\documentclass{amsart}
\usepackage[nohug]{diagrams}
\usepackage[mathscr]{eucal}
\usepackage{amssymb, amsmath, amsthm}
\usepackage{amsfonts}
\usepackage{latexsym}
\usepackage{tabularx,array}
\usepackage{enumerate}
\usepackage[all,arc]{xy}
\usepackage{latexsym}
\usepackage{mathrsfs}
\usepackage{color}
\usepackage{mathtools}
\usepackage{leftidx}
\usepackage{bm}
\usepackage{url}
\usepackage{setspace}
\newfont{\msam}{msam10}

\newtheorem{theorem}[]{Theorem}
\newtheorem{proposition}[]{Proposition}
\newtheorem{corollary}[]{Corollary}
\newtheorem{lemma}[]{Lemma}

\theoremstyle{definition}
\newtheorem{definition}[]{Definition}
\newtheorem{defn}[theorem]{Definition}
\newtheorem{remark}[]{Remark}

\newtheorem{conj}[]{Conjecture}
\newtheorem{question}[]{Question}

\def\remark{\noindent\textbf{Remark.}}

\def\cH{\mathcal H}
\def\ucH{\underline{{\mathcal H}}}

\let\nc\newcommand

\def\bthm{\begin{theorem}}
\def\ethm{\end{theorem}}
\def\blemma{\begin{lemma}}
\def\elemma{\end{lemma}}
\def\bproof{\begin{proof}}
\def\eproof{\end{proof}}
\def\bprop{\begin{proposition}}
\def\eprop{\end{proposition}}
\def\bcor{\begin{corollary}}
\def\ecor{\end{corollary}}
\def\bconj{\begin{conj}}
\def\econj{\end{conj}}
\nc{\la}{\label}
\def\O{\mathcal{O}}
\def\Z{\mathbb{Z}}
\def\N{\mathbb{N}}
\def\Q{\mathbb{Q}}
\def\c{\mathbb{C}}

\def\L {\boldsymbol{L}}
\def\bR{\boldsymbol{R}}
\def\brR{\boldsymbol{\mathrm{R}}}

\def\Com{\mathtt{Com}}
\def\Vect{\mathtt{Vect}}

\def\sAlg{\mathtt{sAlg}}

\def\LAlg{\mathtt{Lie\,Alg}}
\def\DGL{\mathtt{DGLA}}

\def\cDGC{\mathtt{DGCC}}

\def\Mod{\mathtt{Mod}}

\def\cAlg{\mathtt{Comm}}
\def\scAlg{\mathtt{sComm}}

\def\Sets{\mathtt{Sets}}
\def\DGA{\mathtt{DGA}}

\def\cDGA{\mathtt{DGCA}}

\def\C{\mathcal{C}}

\def\Ho{{\mathtt{Ho}}}
\def\mfa{\mathfrak{a}}
\nc{\Ob}{{\rm Ob}}
\nc{\Hom}{{\rm{Hom}}}
\nc{\Homcont}{{\mathcal{H}om}}
\nc{\HOM}{\underline{\rm{Hom}}}
\nc{\DER}{\underline{\rm{Der}}}
\nc{\END}{\underline{\rm{End}}}
\nc{\bSym}{\mathbf{Sym}}
\nc{\Ext}{{\rm{Ext}}}
\nc{\Rep}{{\rm{Rep}}}
\nc{\DRep}{{\rm{DRep}}}
\nc{\NCRep}{\widetilde{\rm{Rep}}}
\nc{\RAct}{{\rm{RAct}}}
\nc{\bs}{\backslash}
\nc{\ob}{{\tt{Obs}}}
\nc{\CE}{\mathcal{C}}
\nc{\TP}{{T\!P}}
\nc{\un}{\underline{n}}
\nc{\um}{\underline{m}}
\nc{\rn}{\langle n \rangle}
\nc{\nn}{{{\natural} {\natural}}}
\nc{\n}{{{\natural}}}
\nc{\A}{\mathbb A}
\nc{\B}{{\mathrm{B}}}
\nc{\Ba}{\overline{\mathrm{B}}}
\nc{\bC}{\overline{C}}
\nc{\bOmega}{\boldsymbol{\Omega}}
\nc{\bB}{\boldsymbol{B}}
\nc{\EXT}{\underline{\rm{Ext}}}
\nc{\TOR}{\underline{\rm{Tor}}}
\def\H{\mathrm H}

\def\HC{\mathrm{HC}}
\def\Loc{\mathrm{Loc}}
\def\HR{\mathrm{HR}}
\def\rHC{\overline{\mathrm{HC}}}

\nc{\End}{{\rm{End}}}
\nc{\GL}{{\rm{GL}}}
\nc{\gl}{{\mathfrak{gl}}}
\nc{\rgl}{\overline{{\mathfrak{gl}}}}
\nc{\g}{{\mathfrak{g}}}
\nc{\h}{{\mathfrak{h}}}
\nc{\PGL}{{\rm{PGL}}}
\nc{\SL}{{\rm{SL}}}
\nc{\sll}{\mathfrak{sl}}
\nc{\cn}{ \mbox{\rm c\^{o}ne} }
\nc{\PSL}{{\rm{PSL}}}
\nc{\ad}{{\rm{ad}}}
\nc{\Ad}{{\rm{Ad}}}
\nc{\dlim}{\varinjlim}
\nc{\plim}{\varprojlim}
\nc{\colim}{{{\rm colim}}}

\newcommand{\ra}{\to}

\newcommand{\HH}{{\rm{HH}}}

\newcommand{\Tor}{{\rm{Tor}}}
\newcommand{\Spec}{{\rm{Spec}}}

\newcommand{\Sym}{\Lambda}

\newcommand{\id}{{\rm{Id}}}

\newcommand{\Tr}{{\rm{Tr}}}

\newcommand{\Ker}{{\rm{Ker}}}

\newcommand{\into}{\,\hookrightarrow\,}

\newcommand{\onto}{\,\twoheadrightarrow\,}
\newcommand{\sonto}{\,\stackrel{\sim}{\twoheadrightarrow}\,}
\def\cb{\boldsymbol{\Omega}}

\def\Top{\mathtt{Top}}

\def\bs{\backslash}

\def\ve{\mathtt{Vect}}

\def\sGr{\mathtt{sGr}}
\def\Gr{\mathtt{Gr}}

\def\ffgr{\mathfrak{G}}
\def\GG{\mathfrak{G}}

\def\sset{\mathtt{sSet}}
\def\lgr{\mathbb{G}}
\def\lin{\mathtt{lin}}

\newcommand{\rar}{\xrightarrow{}}

\newcommand\frgr[1]{\langle #1 \rangle}
\nc{\env}{\mathrm{End}(V)}
\nc{\FT}{\mathcal{C}}

\numberwithin{equation}{section}
\numberwithin{theorem}{section}
\numberwithin{lemma}{section}
\numberwithin{proposition}{section}
\numberwithin{definition}{section}
\numberwithin{corollary}{section}
\numberwithin{example}{section}
\numberwithin{remark}{section}

\newcommand{\sym}{\mathrm{Sym}}

\newcommand{\rH}{\overline{\mathrm{H}}}
\def\cO{\mathcal O}
\newcommand{\Chp}{{\tt Ch}_{\geq 0}}
\newcommand{\sMod}{{\tt sMod}}
\newcommand{\grMod}{{\tt grMod}}
\newcommand{\surj}{{\rm surj}}
\newcommand{\op}{{\rm op}}
\newcommand{\sten}{\, \bar{\otimes} \,}
\newcommand{\sh}{{\rm sh}}
\def\cP{\mathcal{P}}
\newcommand{\dgAlg}{\mathtt{dgAlg}} 
\newcommand{\xra}[2][]{\, \xrightarrow[#1]{#2} \,}

\def\bdf{\begin{defn}}
\def\edf{\end{defn}}
\newcommand{\coeq}{{\rm coeq}}
\def\brm{\begin{remark}}
\def\erm{\end{remark}}
\newcommand{\uten}{\, \underline{\otimes} \,}
\newcommand{\wt}{{\rm wt}}

\theoremstyle{definition}
\def\bdf{\begin{definition}}
\def\edf{\end{definition}}

\newcommand{\bS}{{\mathbb S}}

\newcommand{\cA}{{\mathcal A}}

\newcommand{\mH}{{\mathcal H}}
\newcommand{\umH}{\underline{\mathcal H}}

\newcommand{\LL}{{\mathscr L}}
\DeclareMathOperator*{\Moplus}{\text{\raisebox{0.25ex}{\scalebox{0.75}{$\bigoplus$}}}}
\def\arbreBA{\vcenter{\xymatrix@R=2pt@C=2pt{
&&&&\\
&&&*{}\ar@{-}[ul] & \\
&&*{}\ar@{-}[uurr] \ar@{-}[uull] \ar@{-}[d]     &&\\
&&&&
}}}

\def\arbreAB{\vcenter{\xymatrix@R=2pt@C=2pt{
&&&&\\
&*{}\ar@{-}[ur] &&& \\
&&*{}\ar@{-}[uurr] \ar@{-}[uull] \ar@{-}[d]     &&\\
&&&&
}}}

\def\arbreABC{\vcenter{\xymatrix@R=1pt@C=1pt{
&&&&&&\\
&*{}\ar@{-}[ur] &&&&& \\
&&*{}\ar@{-}[uurr] &&&&\\
&&&*{}\ar@{-}[uuurrr] \ar@{-}[uuulll] \ar@{-}[d] &&&\\
&&&&&&
}}}

\def\arbreBAC{\vcenter{\xymatrix@R=1pt@C=1pt{
&&&&&&\\
&&&*{}\ar@{-}[ul] &&& \\
&&*{}\ar@{-}[uurr] &&&&\\
&&&*{}\ar@{-}[uuurrr] \ar@{-}[uuulll] \ar@{-}[d] &&&\\
&&&&&&
}}}

\def\arbreACB{\vcenter{\xymatrix@R=1pt@C=1pt{
&&&&&&\\
&*{}\ar@{-}[ur] &&&&& \\
&&&&*{}\ar@{-}[uull] &&\\
&&&*{}\ar@{-}[uuurrr] \ar@{-}[uuulll] \ar@{-}[d] &&&\\
&&&&&&
}}}

\def\arbreBCA{\vcenter{\xymatrix@R=1pt@C=1pt{
&&&&&&\\
&&&&&*{}\ar@{-}[ul] & \\
&&*{}\ar@{-}[uurr] &&&&\\
&&&*{}\ar@{-}[uuurrr] \ar@{-}[uuulll] \ar@{-}[d] &&&\\
&&&&&&
}}}

\def\arbreCAB{\vcenter{\xymatrix@R=1pt@C=1pt{
&&&&&&\\
&&&*{}\ar@{-}[ur] &&& \\
&&&&*{}\ar@{-}[uull] &&\\
&&&*{}\ar@{-}[uuurrr] \ar@{-}[uuulll] \ar@{-}[d] &&&\\
&&&&&&
}}}

\def\arbreCBA{\vcenter{\xymatrix@R=1pt@C=1pt{
&&&&&&\\
&&&&&*{}\ar@{-}[ul] & \\
&&&&*{}\ar@{-}[uull] &&\\
&&&*{}\ar@{-}[uuurrr] \ar@{-}[uuulll] \ar@{-}[d] &&&\\
&&&&&&
}}}

\def\arbreACA{\vcenter{\xymatrix@R=1pt@C=1pt{
&&&&&&\\
&*{}\ar@{-}[ur] &&&&*{}\ar@{-}[ul] & \\
&&&&&&\\
&&&*{}\ar@{-}[uuurrr] \ar@{-}[uuulll] \ar@{-}[d] &&&\\
&&&&&&
}}}
\begin{document}

\title{Representation homology of simply connected spaces}
%
\author{Yuri Berest}
\address{Department of Mathematics,
Cornell University, Ithaca, NY 14853-4201, USA}
\email{berest@math.cornell.edu}
\author{Ajay C. Ramadoss}
\address{Department of Mathematics,
Indiana University, Bloomington, IN 47405, USA}
\email{ajcramad@indiana.edu}
\author{Wai-Kit Yeung}
\address{Department of Mathematics,
Indiana University, Bloomington, IN 47405, USA}
\email{yeungw@iu.edu}
\begin{abstract}
Let $G$ be an affine algebraic group defined over field $k$ of characteristic zero. We study the derived moduli space of $G$-local systems on a pointed connected CW complex $X$ trivialized at the basepoint of $X$. This derived moduli space is represented by an affine DG scheme $ \brR \Loc_G(X,\ast) \,$: we call the 
(co)homology of the structure sheaf of $ \brR \Loc_G(X,\ast) $ the {\it representation homology of $X$ in $G$} 
and denote it by $ \HR_\ast(X,G) $. The $ \HR_0(X,G) $ is isomorphic to the coordinate ring of 
the representation variety $ \Rep_G[\pi_1(X)] $ of the fundamental group of $X$ in $G$ --- a well-known
algebro-geometric invariant of $X$ with many applications in topology. The case when $X$ is simply connected 
seems much less studied: in this case, the  $ \HR_0(X,G) $ is trivial but the higher representation 
homology is still an interesting rational invariant of $X$ depending on the algebraic group $G$.
In this paper, we use rational homotopy theory to compute the  $ \HR_\ast(X,G) $ for an arbitrary simply connected
space $X$ (of finite rational type)  in terms of its Quillen and Sullivan algebraic models.
When $G$ is reductive, we also compute the $G$-invariant part of representation homology, $ \HR_\ast(X,G)^G $, and study the question when $ \HR_\ast(X,G)^G $ is free of locally finite type as a graded commutative algebra. This question turns out to be closely  related to the so-called Strong Macdonald Conjecture, a celebrated result in representation theory proposed (as a conjecture) by B. Feigin and P. Hanlon in the 1980's and proved by S. Fishel, I. Grojnowski and C. Teleman in \cite{FGT}. Reformulating the Strong Macdonald Conjecture in topological terms, we give a simple characterization of spaces $X$ for which $ \HR_\ast(X,G)^G $ is a  graded symmetric algebra for any complex reductive group $G$. 
\end{abstract}
\maketitle


\section{Introduction}
The present paper is a sequel to our earlier work, \cite{BRYI, BRYII}, where we study representation 
homology of topological spaces. In \cite{BRYI}, we established basic properties of representation homology, constructed  natural maps and spectral sequences relating it to some well-known homology theories associated with spaces (e.g., higher Hochschild homology and 
homology of based loop spaces) as well as did concrete computations. Further, in \cite{BRYII}, we studied the linearization of representation homology and proved some vanishing theorems for groups, surfaces and certain 3-dimensional
manifolds of interest in geometric topology.

The main aim of this paper is to compute the representation homology of an arbitrary {\it simply connected} space $X$ over a field $k$ of characteristic zero. From \cite{BRYI}, we know that the representation homology of such a space is a rational homotopy invariant (i.e. it depends only on the homotopy type of the rationalization $ X_{\Q} $ of $X$); on the other hand, by a fundamental theorem of D. Sullivan \cite{Su}, the homotopy type of $X_{\Q}$ is completely determined by  its algebraic model: a commutative cochain DG algebra  $ \cA_X $, called the {\it Sullivan model} of $X$. This leads us to the the following natural 
\begin{question}
\la{qst1}
How to express the representation homology of $X$ in terms of $ \cA_X $? 
\end{question}

The representation homology $\HR_\ast(X,G) $ of a space $X$ in an algebraic group $G$ may be thought of as a  multiplicative version of ordinary (co)homology, where the commutative Hopf algebra $\cO(G) $ plays the role of coefficients (see \cite{BRYI}). The $\HR_\ast(X,G) $ is a graded commutative algebra, whose structure depends on the homotopical structure of $X$ and the algebraic structure of $G$. In this regard, representation homology is analogous to higher Hochschild homology, $\HH_*(X,A) $, which can be viewed as a homology of the space $X$ with coefficients in a commutative algebra $A$ (see \cite{P1}). While the two homology theories may be defined in a similar way and are, in fact, closely related\footnote{As shown in \cite{BRYI}, there is a natural isomorphism $\,\HR_*(\Sigma(X_{+}), G) \cong \HH_*(X, \cO(G)) \,$ for any space $X$.}, there is one important difference: unlike $\HH_*(X,A) $, the $\HR_\ast(X,G) $ carries a natural algebraic $G$-action induced by the adjoint action of $G$. Examples show that this action depends on the space $X$ in quite a nontrivial way, which makes representation homology a richer and somewhat more geometric theory than Hochschild homology. When $X$ is simply connected (so that $\HR_0(X,G) = k $) and $G$ is reductive, it is natural to treat $ \HR_\ast(X,G) $ as an object of representation theory --- or even classical invariant theory (in the sense of \cite{Weyl}) --- and
ask basic questions about the structure of  the algebra $\HR_*(X,G) $ as a $G$-module and its subalgebra $\HR_*(X,G)^G $ of $G$-invariants. Perhaps, the first natural question that arises in this direction is

\begin{question}
\la{qst2}
When is the algebra of invariants $\,\HR_*(X,G)^G \,$ free and (locally) finitely generated, i.e. isomorphic to the graded symmetric algebra of a (locally) finite-dimensional graded vector space over $k$?
\end{question}

Besides its intrinsic interest Question~\ref{qst2} turns out to be related to some of the deeper questions in Lie theory and algebraic representation theory. Our second aim in this paper is to  shed some new light on these questions linking them  to topology. To state our results we first recall a few basic facts about representation homology (for more details and proofs the reader is referred to \cite{BRYI}).
\subsection{Three definitions of representation homology}
There are (at least) three different ways to define representation homology. Historically the first and (arguably) most appealing definition  comes from derived algebraic geometry (see, e.g., \cite{K}, \cite{PTVV}, \cite{PT}, \cite{To}). Let $G$ be an affine algebraic group defined over a field $k$ of characteristic zero. Given a pointed connected CW complex $X$ consider the (framed) moduli space $ \Loc_G(X,\ast) $ of $G$-local systems on $X$ with trivialization at the basepoint of $X$. As shown in \cite{K}, this classical moduli space has a natural derived extension which is represented by an affine differential-graded (DG) scheme $\brR \Loc_G(X,\ast) $. The structure sheaf of $ \brR \Loc_G(X,\ast) $ is, by definition, a (negatively graded) commutative cochain DG algebra whose cohomology is a homotopy invariant of $X$. We set
\begin{equation}
\la{locsys}
\HR_*(X,G) := \H^{-\ast}[\cO_{\boldsymbol{\mathrm R} \Loc_G(X,\ast)}]
\end{equation}
and call $ \HR_*(X,G) $ the {\it representation homology of  $X$ in $ \,G$}. This terminology is motivated by the fact that $ \Loc_G(X,\ast) $  can be identified with the classical representation
scheme $ \Rep_G[\pi_1(X)] $, parametrizing the representations of the fundamental group of $X$ in $G$, and  $ \HR_0(X,G) $ is thus naturally isomorphic to $\cO[\Rep_G(\pi_1(X))]$, the affine coordinate ring of  $ \Rep_G[\pi_1(X)] $. 

Another, less geometric but more general and conceptually simpler definition was proposed in \cite{BRYI}. This definition rests on a fundamental result in simplicial homotopy theory, due to D. Kan \cite{Kan1}, that describes the homotopy types of pointed connected spaces in terms of simplicial groups. More precisely, Kan's Theorem asserts that the model category $\sGr $ of simplicial groups is Quillen equivalent to the category $ \sset_0 $ of reduced simplicial sets, which is, in turn, Quillen equivalent to the category $ \Top_{0,*} $ of pointed connected (CGWH)  spaces; thus, there are natural equivalences of homotopy categories
\begin{equation}\la{tsg}
 \Ho(\sGr)  \,\cong\, \Ho(\sset_0)\,\cong\, \Ho(\mathtt{Top}_{0,\ast})\ .
\end{equation}
The starting point for our construction of representation homology is the simple observation that the functor of points $ G: \cAlg_k \to \Gr $ of any affine algebraic group (scheme) $G$ has a left adjoint
\begin{equation}
\la{lgrsc}
(\,\mbox{--}\,)_G:\,\mathtt{Gr} \to \cAlg_k
\end{equation}
which \mbox{---} when applied to a given group $ \Gamma $ \mbox{---} gives the coordinate ring of the affine scheme $ \Rep_G(\Gamma) $: i.e., $\Gamma_G = \O[\Rep_G(\Gamma)]\,$. Thus  the functor \eqref{lgrsc} provides an alternative (dual) description of the representation scheme $ \Rep_G(\Gamma) $ and is called the {\it representation functor} in $G$. Now, to define representation homology in $G$ we simply derive \eqref{lgrsc} using the standard simplicial technique in homological algebra \cite{Q1}. First, we prolong the adjoint functors $\,(\,\mbox{--}\,)_G: \Gr \rightleftarrows \cAlg_k : G\,$ to the simplicial categories $ \sGr $ and $\scAlg_k $ (by applying them degreewise to the corresponding
simplicial objects) and then replace the resulting adjunction
\begin{equation}
\la{sAdj}
(\,\mbox{--}\,)_G\,: \,\sGr\, \rightleftarrows \, \scAlg_k\, :\, G\ .
\end{equation}
with its `universal homotopy approximation' represented by derived functors. To be precise, the main theorem of \cite{BRYI} (see {\it loc. cit.} Theorem~1.1) asserts each of the adjoint functors in \eqref{sAdj} has a total derived functor
(left and right, respectively) and these derived functors form an adjoint pair at the level of homotopy categories\footnote{We should warn the reader that the functors \eqref{sAdj} 
do {\it not} form a Quillen pair between the {\it model} categories $ \sGr $ and $\scAlg_k $ for standard (projective) model structures. The existence of the derived functors and the corresponding derived adjunction \eqref{derAdj}
is a nontrivial fact that does not follow immediately from Quillen's well-known theorem for model categories \cite{Q1}.}:
\begin{equation}
\la{derAdj}
\L(\,\mbox{--}\,)_G\,: \,\Ho(\sGr)\, \rightleftarrows \, \Ho(\scAlg_k)\, :\, \bR G
\end{equation}

We can now make the following definition which
will be our main definition for the present paper ({\it cf.} \cite[Definition~3.1]{BRYI}).

\begin{definition}
\la{DRepX}
For a space $X \in \Top_{0,*} $, we choose a simplicial group model $ \Gamma X $  and define the {\it representation homology of $X$ in $G$}  by
\begin{equation}
\la{hrX}
\HR_\ast(X,G) \,:= \,\pi_\ast \L(\Gamma X)_{G} \,:= \,\H_\ast[N\L(\Gamma X)_{G}]\ ,
\end{equation}
where $ N $ stands for the standard (Dold-Kan) normalization functor (see Appendix~\ref{sdk}). 
\end{definition}
Note that, since $ \L(\,\mbox{--}\,)_G $ is a homotopy functor on simplicial groups, formula \eqref{hrX} does not depend on the choice of a simplicial group model of $X$. In fact,
there are several natural models that can be used in practical computations
(see \cite{BRYII}). In this paper, we will use most exclusively the so-called {\it Kan loop group model}
$\, \Gamma = \lgr{X}$, which is a semi-free simplicial group functorially attached to the space $X$ (see \cite[Chap. V]{GJ} or \cite[Sect. 2]{BRYI} for a brief summary of this construction).
Since semi-free simplicial groups are cofibrant objects in $\sGr$, formula \eqref{hrX} simplifies in this case to
\begin{equation}
\la{KhrX}
\HR_\ast(X,G) \, = \,\pi_\ast(\lgr{X})_{G} \ .
\end{equation}

To compare Definition~\ref{DRepX} with the algebro-geometric construction of representation homology, \eqref{locsys}, we associate to the derived representation functor \eqref{hrX} the {\it derived representation scheme}
$$
\DRep_G(X) := \boldsymbol{\mathrm{R}}\Spec\, [\L(\Gamma X)_{G}]\ .
$$
Here `$\, \boldsymbol{\mathrm{R}}\Spec \,$'  stands for the To\"en-Vezzosi derived Yoneda functor 
\cite{TV05, TV08} that assigns to a simplicial commutative algebra $A$ --- a derived ring in 
terminology of \cite{TV08} --- the simplicial presheaf (prestack)  
$$
\boldsymbol{\mathrm{R}}\Spec(A):\, \mathtt{dAff}_k^{\rm op} := \scAlg_k \ \to \ \sset\ ,\quad
B \,\mapsto\, \underline{\mathrm{Hom}}(Q(A),\,B)\ ,
$$
where $Q(A)$ is a cofibrant model for $A$ and $ \underline{\mathrm{Hom}} $ is the simplicial mapping 
space (function complex) in $ \scAlg_k $. For any $A \in \scAlg_k $, the prestack  $ 
\boldsymbol{\mathrm{R}}\Spec(A) $ satisfies the descent condition
for \'etale hypercoverings and hence defines a derived stack (which is a derived affine scheme in the 
sense of \cite{TV08}). Now, in \cite[Appendix~A.2]{BRYII}, we showed that for any pointed connected CW 
complex $X$, there is an equivalence of derived stacks $\, \DRep_G(X) \simeq \boldsymbol{\mathrm R} \Loc_G(X,\ast)\,$. This implies that the two definitions of representation homology --- \eqref{locsys} 
and \eqref{hrX} --- actually agree.

Our third definition of $\, \HR_\ast(X,G) \,$ --- perhaps the most straightforward and elementary one ---
is given  in terms of functor homology. Let $ \mathfrak{G} $ denote the (small) subcategory of $ \Gr $
whose objects $ \langle n \rangle := \mathbb{F}_n $ are f. g. free groups, one for each $ n \ge 0 $, 
and the morphisms are arbitrary group homomorphisms. This category carries a natural (strict) monoidal 
structure,  with  product $ \ast:\,\mathfrak{G} \times  \mathfrak{G} \to  \mathfrak{G} $ being the free product (coproduct) of free  groups: $ \langle n \rangle \ast 
\langle m \rangle = \langle n + m \rangle $. It is known  that every commutative Hopf algebra defines a (strict) monoidal functor on $ \GG $ with values in $ \cAlg_k $, and conversely,
every such functor corresponds to a commutative Hopf algebra
(see, e.g., \cite{P2}). As in \cite{BRYI}, given a commutative Hopf algebra $ \cH $,  we denote the corresponding
functor by
\begin{equation}
\la{hopff}
\ucH:\, \GG \to \cAlg_k \ ,\quad \rn \to \cH^{\otimes n} \ .
\end{equation}
Note that \eqref{hopff} naturally extends to a functor on all groups: $ \Gr \to \cAlg_k $ by taking the left Kan extension along the inclusion $ \GG \into \Gr$. To avoid  complicated notation we will use the same symbol $ \umH $ to denote the functor \eqref{hopff} and its natural extension to $\Gr $ (moreover, we will often drop the underline in this symbol when there is no danger of confusion).  Now, to define representation homology with coefficients in $ \mH $ we simply precompose the corresponding functor $ \umH: \Gr \to \cAlg_k $ with the Kan loop group of a given space $X$: the result is the simplicial commutative algebra
$$
\cH(\lgr{X}):\ \Delta^{\rm op} \xrightarrow{\lgr{X}} \Gr \xrightarrow{\ucH} \cAlg_k
$$
whose homology we denote by\footnote{As noted in \cite{BRYI}, this simple definition of representation homology is analogous to Pirashvili's definition of higher Hochschild homology in \cite{P1}.}
\begin{equation*}
\HR_\ast(X,\cH) := \pi_*[\cH(\lgr{X})] = \H_\ast[N \cH (\lgr{X})]
\end{equation*}
For $\mH = \cO(G) $, where $G$ is an affine algebraic group scheme over $k$, it is easy to show that there is a natural isomorphism (see \cite[Prop. 4.1]{BRYI}):
\begin{equation}
  \la{def3}
  \HR_\ast(X,\,\cO(G)) \,\cong \,\HR_\ast(X,G)\ .
\end{equation}
Thus, we may think of the representation homology as a homology of a space with coefficients in commutative Hopf algebras in the same way as one thinks of the ordinary homology as a homology with coefficients in abelian groups
or the higher Hochschild homology \cite{P1} as a homology with coefficients in commutative algebras.

Now, for any (discrete) group $ \Gamma \in \Gr $, the group algebra
$k[\Gamma] $ has a natural {\it cocommutative} Hopf algebra structure and therefore defines a contravariant monoidal functor on $\GG$:
$$
\underline{k[\Gamma]}:\ \GG^{\rm op} \to \cAlg_k\ ,\quad \rn \mapsto k[\Gamma]^{\otimes n}\ .
$$
Regarding $ \underline{\cO(G)} $ and $ \underline{k[\Gamma]} $
as linear functors on $\GG $ (with values in $\Vect_k $), we can form their  tensor product $\,k[\Gamma] \otimes_{\GG} \cO(G) \,$. It turns out that there is a natural isomorphism: 
$\,
k[\Gamma] \otimes_{\GG} \cO(G)\,\cong\,\cO[\Rep_G(\Gamma)] 
\,$; more generally, it is shown in \cite{BRYI} that
\begin{equation}
  \la{torG}
  \HR_\ast(\mathrm{B}\Gamma,\,G) \,\cong \,\Tor_\ast^{\GG}(k[\Gamma],\,\cO(G))\ ,
\end{equation}
where $ \Tor^{\GG}_\ast $ is the (homology of the) classical derived tensor product $\,\otimes_{\GG}^{\L}\,$ between covariant and contravariant linear functors on $\GG$.
The `Tor'-formula \eqref{torG} is remarkable for two reasons:
first, it gives a natural interpretation of representation homology in terms of usual (abelian) homological algebra,
placing it in one row with other classical invariants, such as
Hochschild and cyclic homology (see, e.g., \cite{L}). Second --- as we will see in this paper --- it provides an efficient tool for computations\footnote{We should also mention that, in recent years homological algebra in functor categories over $\mathfrak{G}$ 
has been extensively used in computations of stable homology of automorphism groups of free groups and the study of related questions
of $K$-theory and topology (see, e.g., \cite{DV1, DV2, DV3, PV, V}) and also \cite[Sect. 7]{BRYI}).}.

\subsection{Main results} Throughout, $k$ stands for a commutative base field, which is always assumed to be of characteristic zero but (unless specified so) not necessarily algebraically closed. 

Our answer to Question~\ref{qst1} can be encapsulated into the following theorem which is the main result of the present paper.

\bthm
\la{hrsullivan}
Let $X$ be a $1$-connected pointed space of finite rational type with Sullivan model $ \cA_X $. Let $ \bar{\cA}_X $ denote the augmentation ideal of $ \cA_X $  corresponding to the basepoint of $X$. 
 
 $(a)$ For any affine algebraic group $G$ defined over $ k $ with Lie algebra $ \g $, there is an isomorphism of graded commutative algebras
\begin{equation*}
\la{sulltype}
\HR_\ast(X,G) \,\cong\, \H^{- \ast}(\g(\bar{\cA}_X);\, k)\ ,
\end{equation*}
where $ \g(\bar{\cA}_X) $ is the current Lie algebra of $\g $ over 
the commutative DG algebra $ \bar{\cA}_X $.

$(b)$ If $G$ is a {\rm reductive} affine algebraic group over $k$, then
\begin{equation*}
\la{sulltypeG}
\HR_\ast(X,G)^G \,\cong\, \H^{-\ast}(\g(\cA_X), \g; \, k)\ ,
\end{equation*}
where $ \g(\cA_X) $ is the current Lie algebra over 
$ \cA_X $ and $ \g \subseteq \g(\cA_X)  $ is its canonical Lie subalgebra.

\ethm
Theorem~\ref{hrsullivan} needs some explanations. First, 
recall that for a Lie algebra $\g$ and a commutative DG algebra $\cA$,
the current Lie algebra $ \g(\cA) $ is defined to be
the tensor product $\g(\cA) :=\g \otimes \cA $ with Lie bracket $[\xi \otimes a\,,\,\eta \otimes b] \,:=\,[\xi\,,\,\eta] \otimes ab$ and the differential $d(\xi \otimes a) := \xi \otimes da $. If $X$ is a pointed $1$-connected topological space of finite rational type, its Sullivan model $\cA_X$ is an augmented commutative cochain DG algebra, so we can form the current Lie algebras $\g(\cA_X)$ and $\g(\bar{\cA}_X)$, both of which are cohomologically graded. In Part $(a)$ of Theorem~\ref{hrsullivan},  $\H^{-\ast}(\g(\bar{\cA}_X); k)$ stands for the classical (Chevalley-Eilenberg) cohomology of the Lie algebra $ \g(\bar{\cA}_X) $ with trivial coefficients; in Part $(b)$, $\,\H^{-\ast}(\g({\cA}_X),\g; k)\,$ is the relative Lie algebra cohomology of the canonical pair  $ \g \subseteq \g(\cA_X) $. The `minus' sign in the superscript of both cohomologies indicates that they are considered with homological grading.

The proof of Theorem~\ref{hrsullivan} is fairly long and technical: it occupies most of Section~\ref{sect6} and relies heavily on results of Quillen \cite{Q}. For reader's convenience, we  outline the main steps of this proof in Section~\ref{sect6.3.1}. Here we mention only two key results that are of independent interest.
The first is Theorem~\ref{conj1} which expresses the representation homology of a simply connected space $X$ in terms of its Quillen DG Lie algebra model $\mfa_X $:
\begin{equation}
\la{compthm}
\HR_\ast(X,G) \,\cong \, \HR_\ast(\mfa_X, \g)\ .
\end{equation}
We call Theorem~\ref{conj1} the `Comparison Theorem' as it `compares' two representation homology functors: one with coefficients in an algebraic group $G$ and the other with coefficients in its Lie algebra $\g$. The second notable result is Theorem~\ref{liereps}: it provides a functor homology interpretation --- a natural counterpart of  `Tor'-formula \eqref{torG} --- for representation homology of Lie algebras\footnote{We briefly review the definition of representation homology of Lie algebras in Section~\ref{DREP} below.}:
\begin{equation}
\la{toralg}
 \HR_\ast(\mfa, \g) \,\cong\, \Tor_{\ast}^{\GG}(U\mfa,\,\g)\ .
\end{equation}
Both isomorphisms \eqref{compthm} and \eqref{toralg} are deduced from Theorem~\ref{dtpliemodel1}, which is a result in rational homotopy theory --- a natural refinement of one of the main results of \cite{Q}.

We now turn to Question~\ref{qst2}. We will approach this question by constructing topologically  some natural maps with values in $\HR_\ast(X, G)^G $ whose images --- in good cases --- will generate  $\HR_\ast(X, G)^G $ as an algebra. Given a simply connected space $X$, we consider the space $\LL X $ of all continuous maps $S^1 \to X $ from
the topological circle $S^1$ to $X$ equipped with compact open topology. This classical space, called the {\it free loop space of $X$}, carries a natural $S^1$-action induced by the action of $S^1$ on itself by rotations: thus, we can define its (reduced) $S^1$-equivariant homology $\rH^{S^1}_\ast(\LL X,\, k) $. It is well known that, when $k =\Q$ (or more generally, $k$ has characteristic $0$), there is a natural direct sum decomposition
\begin{equation}
\la{hodge}
\rH^{S^1}_\ast(\LL X,\,k) \,=\,\bigoplus_{p=0}^{\infty}\,\rH^{S^1,\, (p)}_\ast(\LL X,\,k)\ ,
\end{equation}
which is usually called the {\it Hodge decomposition of} $\,\rH^{S^1}_\ast(\LL X,\, k) $. The $p$-th direct summand in \eqref{hodge} --- the Hodge component of degree $p$ -- is defined topologically as the common eigenspace of the
degree $p$ Frobenius operations, i.e. the graded endomorphisms of $ \rH^{S^1}_\ast(\LL X,\, k) $
induced by the finite coverings of the circle: $S^1 \to S^1 $, $\,e^{i\theta} \to e^{i n \theta} $, corresponding to the eigenvalues $n^p$, $\, n\ge 0 \,$ (see Section~\ref{secdrintr}). A theorem of Burghelea, Fiedorowicz and Gajda (see \cite{BFG}, Theorem~A) asserts that, if all (rational) Betti numbers of $X$ are finite, each Hodge component of $\, \rH^{S^1}_\ast(\LL X,\, k) $ is locally finite: i.e., 
\begin{equation}\la{finhodge}
\dim_k\, \rH^{S^1,\, (p)}_i(\LL X,\,k) < \infty \quad \mbox{for all}\ i \ge 0\ \mbox{and all}\ p\ge 0 \ .
\end{equation}

Now, assume that $ k = \c $ and $G$ is a complex reductive group of rank $ l \ge 1 $. Let $\g$ be the Lie algebra of $G$ with classical exponents  $\, \{m_1, \, m_2, \, \ldots \,,\, m_l\}\, $, and let $I(\g) := \sym(\g^*)^G $ be the ring of $G$-invariant polynomials on $\g$. It is well known that $ I(\g) $ is generated by $l$ algebraically independent homogeneous polynomials $\,P_1,\, P_2,\, \ldots,\, P_l \, $, such that $ \deg(P_i) = m_i + 1 $ for $\,i=1,2,\ldots, l\,$. In Section~\ref{secdrintr}, for every such generator $ P_i $, we construct a natural linear map
\begin{equation*}
\la{dttr}
\rH^{S^1,\, (m_i)}_\ast(\LL X,\,\c) \,\to\, \HR_*(X, G)^G
\end{equation*}
defined on the $m_i$-th Hodge component of  \eqref{hodge}. Assembling these maps (for all $\,i=1,2,\ldots, l\,$), we get a  graded algebra homomorphism
\begin{equation}
\la{Idrinhom1}
\Sym \big[ \Moplus_{i=1}^l \rH_\ast^{\,S^1,\,(m_i)}({\mathscr L}X,\,\c)\big] \,\to
\, \HR_\ast(X, G)^G\ ,
\end{equation}
which we call a {\it Drinfeld homomorphism} for  $(X,\,G)$.

Note that if $G$ is an algebraic torus, then $ m_i = 0 $ for all $i=1,2,\ldots, l $, and $\,\HR_\ast(X, G)^G = \HR_*(X,G)\,$, because $G$ is commutative. On the other hand, for any simply connected space $X$,  we have 
$$
\rH_\ast^{\,S^1,\,(0)}({\mathscr L}X,\,\c) \cong \H_{*+1}(X,\,\c) \ , 
$$
where the isomorphism is given by composition of the classical
Gysin map $ \rH_\ast^{S^1}({\mathscr L}X,\,\c) \to  \rH_{\ast+1}({\mathscr L}X,\,\c) $ and the natural map $\rH_{\ast+1}({\mathscr L}X,\,\c) \to \rH_{\ast+1}(X,\,\c) $ induced by evaluation of loops at the origin.
Thus, for an algebraic torus, the Drinfeld homomorphism becomes
$$ 
\Sym \big[\H_{\ast+1}(X;k)^{\oplus  l}\big]\,\to\, \HR_\ast(X, G)\ .
$$
A simple calculation with a minimal Quillen model shows that the above map is an {\it isomorphism} for any simply connected space $X$ 
and, in fact, for any commutative --- not necessarily diagonalizable --- algebraic group $G$ (see Theorem~\ref{gagm}). Thus, we get an answer to Question~\ref{qst2}, though in a very special and somewhat trivial case. 

Suppose now that $G$ is an arbitrary complex reductive group. Then we can ask: 

\begin{question}
\la{qst3}
For which spaces $X$ is the Drinfeld homomorphism \eqref{Idrinhom1} an isomorphism?
\end{question}

The following theorem, which is our second main result in this paper, specifies simple conditions on cohomology of the space $X$ that are sufficient for \eqref{Idrinhom1} to be an isomorphism for {\it all} reductive groups $G$.
\bthm[see Theorem~\ref{drinhomcprxs}] 
\la{drinhomcpr}
Assume that the rational cohomology algebra  $\H^\ast(X;\Q)$ of a simply connected space $X$ is either generated by one element $($in any dimension$)$ or {\rm freely} generated by two elements: one in even and one in odd dimensions. Then, the Drinfeld homomorphism \eqref{Idrinhom1} is an isomorphism for $X$ and any complex reductive algebraic group $G$.
\ethm
Note that any space $X$ satisfying the assumptions of Theorem~\ref{drinhomcpr}  obviously satisfies the assumptions of the Burghelea-Fiedorowicz-Gajda Theorem \cite{BFG}, which ensures the (local) finiteness of all Hodge components of $\LL X$: see \eqref{finhodge}. Thus, Theorem~\ref{drinhomcpr} combined with \cite{BFG} provides an answer to Question~\ref{qst2}:
\begin{corollary}
\la{ansq2}
If $X$ satisfies the conditions of Theorem~\ref{drinhomcpr}, then for any complex reductive group $G$,  $\, \HR_*(X,G)^G $ is a free graded commutative algebra of locally finite type over $ \c $.
\end{corollary}
Theorem~\ref{drinhomcpr} relies on (part $(b)$ of) Theorem~\ref{hrsullivan} and a certain (minor) refinement of the main result of the paper \cite{FGT} by S. Fishel, I. Grojnowski and C. Teleman. This last paper settles the so-called {\it Strong Macdonald Conjecture} --- a deep and celebrated result in representation theory proposed as a conjecture by I. Macdonald \cite{Macd}, B. Feigin \cite{Fe}, and P. Hanlon \cite{H1, H2} in the early 80's and proved (in full generality) in \cite{FGT}. The Strong Macdonald Conjecture comprises actually two cases: the first describes the structure of cohomology of the nilpotent Lie algebras $ \g[z]/(z^{r+1}) $ (see \cite{H2} for the case $ \g = \gl_n $ and \cite[Theorem~A]{FGT} for
an arbitrary reductive $\g$) and the second describes the  cohomology of the Lie superalgebra $ \g[z,s] $ (see \cite{Fe} and \cite[Theorem~B]{FGT}).
These two cases roughly correspond to the two cases of Theorem~\ref{drinhomcpr}. Thus, Theorem~\ref{drinhomcpr} gives a topological meaning to the full Strong Macdonald Conjecture. The proof of \cite{FGT} is an algebraic {\it tour de force}. Given the simplicity of our topological reformulation, it is tempting to expect that topology might also lead to a new simpler proof. We leave this as a project for the future.

We would like to conclude this Introduction with a few nice examples illustrating Corollary \ref{ansq2}. Let us consider the  spaces $X$ with rational cohomology algebra $\, \H^*(X, \Q) \cong \Q[z]/(z^{r+1}) \,$, where the generator $z$ is in 
even dimension $ d \ge 2 $. The most well-known examples of such spaces are the even-dimensional spheres $  \bS^{2n}$  ($r=1,\,d=2n$)
and the classical projective spaces: namely, the complex projective spaces $\c\mathbb{P}^r$ ($r \ge 1,\,d=2 $), the quaternionic projective spaces $\mathbb{H}\mathbb{P}^r$ ($r \ge 1,\,d=4$) and the Cayley plane $\mathbb{O}\mathbb{P}^2$ ($r=2,\,d=8$). For these spaces, we have (see Corollary~\ref{rephomcpr}):
$$
\HR_{\ast}(X, G)^{G}\,\cong\, 
\Sym\,[\xi^{(i)}_1,\,\xi^{(i)}_2,\,\ldots,\, \xi^{(i)}_{r}\,:\, i=1,2,\ldots,l]\ ,
$$
where the generators $\xi^{(i)}_{j}$ have homological degree 
$$
\deg\,\xi^{(i)}_j\,=\, \,(d(r+1)-2)m_i+dj-1\ .
$$
%
%
Notice that, in this case, the algebra $ \HR_*(X,G)^G $ is generated by finitely many elements of {\it odd} degrees: hence, it is finite-dimensional (as a vector space) and concentrated in finitely many homological degrees. In fact, knowing the exact degrees of generators,  it is easy to calculate the exact upper bound for the vanishing of $\, \HR_n(X,G)^G $:
$$
\sum_{i=1}^l\,\sum_{j=1}^r \,\deg\,\xi^{(i)}_j\,=\,\frac{1}{2}\,r\,(d(r+1) - 2)\,\dim G\ .
$$
Somewhat miraculously, this allows one to determine the exact upper bound for the {\it full} \,representation homology  of $X$ (see Lemma \ref{vanrephomcpr}):
\begin{equation*}
\la{uppb}
\HR_n(X, G) \,=\,0\quad\, \mbox{for all}\ \ n \,>\,\frac{1}{2}\,r\,(d(r+1) - 2)\,\dim G \ .
\end{equation*}
%
%
Now, the weighted Euler-Poincar\'{e} series of $ \HR_*(X,G)^G $ is given by the polynomial
$$ 
P_{X,G}(q,z)\,=\, \prod_{i=1}^l \prod_{j=1}^r (1\,+\,q^{j+m_i(r+1)}\,z^{\deg\,\xi^{(i)}_j})\ ,
$$
which specializes (at $z=-1$) to the following (weighted) Euler characteristic
\begin{equation*} 
\la{eulercharcpr1} 
\chi_{X,G}(q)\,=\,\prod_{i=1}^l \prod_{j=1}^r (1\,-\,q^{j+m_i(r+1)})\ . 
\end{equation*}
The latter can be also computed --- by Theorem~\ref{hrsullivan}$(b)$ --- as an Euler characteristic of the Chevalley-Eilenberg complex $\C^{-\ast}(\g({\mathcal A}_X), \g; \c) $, where ${\mathcal A}_X $ is the (minimal) Sullivan model of the corresponding space $X$. The resulting equality of Euler characteristics gives the following combinatorial identity
\begin{equation*}  
\mathrm{CT}\left\{ \prod_{j=0}^r\,\prod_{\alpha\,\in\,R} (1\,-\,q^j e^{\alpha})\right\}\,=\, \prod_{i=1}^l \prod_{j=1}^r \frac{1\,-\,q^{j+m_i(r+1)}}{1\,-\,q^j}\ ,
\end{equation*}
which is Macdonald's famous $q$-Constant Term  identity \cite{Macd}. For more examples and explicit calculations we refer the reader to Section~\ref{Sect4.4}.

\subsection*{Appendix} 
In Appendix~\ref{doldkan}, we describe an abstract monoidal version of the classical Dold-Kan correspondence relating the category of (non-negatively graded) DG $\cP$-algebras and the category of simplicial $\cP$-algebras for an arbitrary $k$-linear operad $\cP$. This version of the  Dold-Kan correspondence is needed for our proof of Comparison Theorem in Section~\ref{sect6}. The main result of  Appendix~\ref{doldkan} is Theorem \ref{NNstar_equiv}, which states that when $k$ is a field of characteristic $0$, there is a Quillen equivalence between the category of (non-negatively graded) DG $\cP$-algebras and the category of simplicial $\cP$-algebras. Various special  cases of this theorem have appeared in the literature. First of all, when $\cP$ is the Lie operad, a slightly weaker version (namely, a Quillen equivalence between the category of {\it positively graded} DG Lie algebras and reduced simplicial Lie algebras) was proved in \cite[Part I, Theorem 4.6]{Q}. In {\it loc. cit.} Quillen also outlines a proof  for the commutative operad (which controls commutative unital  $k$-algebras) under the same reducedness assumptions. For general (non-reduced) commutative algebras, the proof of the Dold-Kan correspondence is given in \cite[Proposition~A.1]{TV}. The case of the associative operad is treated in greater generality (for any commutative ring $k$) in \cite{SS}, where the DG associative algebras and simplicial associative algebras are viewed as monoids in the (symmetric) monoidal model categories of chain complexes and simplicial $k$-modules, respectively. In this case, the Dold-Kan correspondence follows from an abstract comparison theorem between monoids in different (symmetric) monoidal model categories. The arguments that establish each of these special cases seem to apply only to the case in hand. To the best of our knowledge, a unified proof for any linear operad is missing in the literature. Our Theorem \ref{NNstar_equiv} fills in this gap\footnote{We should mention, however, that one of the key arguments that we use in our proof of Theorem \ref{NNstar_equiv} is sketched in \cite[Remark 6.4.5]{Fre} in the special case of the commutative operad.}.
Theorem \ref{NNstar_equiv} is crucial for the proof of our Theorem~\ref{hrsullivan}. While Quillen's original result for reduced DG Lie algebras is sufficient for this proof, the full strength of Theorem \ref{NNstar_equiv} is needed to prove Proposition~\ref{liereps}, which is an interesting result on its own.

\subsection*{Outline of the paper} 
The paper is organized as follows. In Section~\ref{DREP}, we recall the definition of representation homology of Lie algebras from \cite{BFPRW} and prove our first result,  Theorem~\ref{liereps}, which gives a realization of this kind of representation homology
as functor homology, see \eqref{toralg}. In Section~\ref{sect6}, we prove our main result, Theorem~\ref{hrsullivan}, answering Question~\ref{qst1} stated in the beginning of the Introduction. We deduce this result from Theorem~\ref{conj1} --- the Comparison Theorem --- which expresses representation homology of a simply connected space in terms of its Quillen Lie model. The Comparison Theorem is technically the most involved result of this paper: its proof occupies the whole of Section~\ref{pfconj1} (with a brief outline given in Section~\ref{sect6.3.1}). We close Section~\ref{sect6} with a conjectural generalization of Theorem~\ref{conj1} to non-simply connected spaces (see Conjectue~\ref{conj5}). Our conjecture is inspired by the recent work \cite{BFMT1,BFMT2,BFMT3}  of Buijs, F\'{e}lix, Murillo and Tanr\'{e} who proposed a natural generalization of Quillen models to non-simply connected. In Section~\ref{S4}, after some necessary preliminaries we construct the Drinfeld homomorphism \eqref{Idrinhom1} and prove our second main result, Theorem~\ref{drinhomcpr}, that gives (partial) answers to Question~\ref{qst2} and Question~\ref{qst3}. We also describe
explicitly the algebra $ \HR_\ast(X,G)^G $ for all spaces $X$ satisfying the conditions of Theorem~\ref{drinhomcpr} and give many
concrete examples of such spaces. Finally, we show how the classical root systems identities --- the original $q$- and $(q,t)$-Macdonald conjectures proposed in \cite{Macd} and proved in \cite{Ch} --- arise from our examples. The last section is an Appendix on the Dold-Kan correspondence (see above) that can be read independently of the rest of the paper.

\subsection*{Acknowledgements} 
The first author would like to thank Giovanni Felder and Ivan Cherednik for interesting questions and many stimulating discussions that occurred during his sabbatical stay at Forschungsinstitut f\"{u}r Mathematik (ETH, Z\"{u}rich) in Fall 2019. He is grateful to ETH for its hospitality and financial support. 
The research of the first author was partially supported by 2019 Simons Fellowship and NSF grant DMS 1702372. Research of the second author was partially supported by NSF grant DMS 1702323.

\section{Representation homology of Lie algebras}
\la{DREP}
The goal of this section is to prove Theorem~\ref{liereps} which gives a functor homology interpretation --- a counterpart of formula \eqref{torG} --- for the representation homology of Lie algebras. This result is a key step in the proof of our main theorem in Section~\ref{sect6}.
We begin by recalling the construction of representation homology in the form it first appeared in \cite{BFPRW}.

\subsection{The representation functor for Lie algebras}
Let $ \g $ be a finite-dimensional Lie algebra over $k$. Given an (arbitrary) Lie algebra $ \mfa \in \LAlg_k $, the moduli scheme $ \Rep_{\g}(\mfa) $   classifying the $k$-linear representations of $ \mfa $ in $ \g $ is defined by the functor on the category of commutative algebras
\begin{equation*}
\la{repg}
\Rep_{\g}(\mfa) :\, \cAlg_k \to \Sets\ , \quad
A \mapsto \Hom_{\tt Lie}(\mfa, \, \g(A))\ ,
\end{equation*}
that assigns to an algebra $A$ the set of families of representations of $ \mfa $
in $ \g $ parametrized by the $k$-scheme $ \Spec(A) $. It is easy to show that  
this functor is representable, and the commutative algebra $ \mfa_\g $ representing $ \Rep_{\g}(\mfa) $ has the following canonical presentation  ({\it cf.} \cite[Prop.~6.4]{BFPRW}):
\begin{equation}
\la{ag}
\mfa_{\g} \, = \,\frac{\Sym_k(\mfa \otimes \g^*)}{
\langle\!\langle\, (x\otimes \xi^*_{1})\cdot (y \otimes \xi^*_{2}) -
(y \otimes \xi^*_{1})\cdot (x \otimes \xi^*_{2}) - [x, y] \otimes \xi^* \,\rangle\!\rangle}\ .
\end{equation}
Here $\, x \otimes \xi^* $ are elements of $ \mfa \otimes \g^* $, where $ \g^* := \Hom_k(\g,k) $ is the vector space dual to $ \g $,  and $\,\xi^* \mapsto \xi^*_{1} \wedge \xi^*_{2}\,$ is the linear map $\,\g^* \to \wedge^2 \g^* $ dual to the Lie bracket on $ \g $.
The tautological (universal) representation $\,\varrho_{\g}: \mfa \to \g(\mfa_\g) \,$ is given by the natural Lie algebra map
\begin{equation} \la{univrep}
\mfa \to \mfa \otimes \g^* \otimes \g \into \Sym_k(\mfa \otimes \g^*) \otimes \g \onto \mfa_{\g} \otimes \g = \g(\mfa_\g)\ ,\quad
x \mapsto \sum_i\, [x \otimes \xi_i^*] \otimes \xi_i\ ,
\end{equation}
where $ \{\xi_i \} $ and $\{\xi_i^*\} $ are dual bases in $ \g $ and $ \g^* $.
The $k$-algebra $ \mfa_{\g} $ has a canonical augmentation
$\,\varepsilon: \mfa_{\g} \to k \,$ induced by the zero map $\, \mfa \otimes \g^* \to 0 \,$. The assignment $\,\mfa \mapsto (\mfa_{\g}, \varepsilon)   \,$ defines a functor with values in the category of augmented commutative algebras
\begin{equation}
\la{repcl}
(\,\mbox{--}\,)_{\g}\,:\ \LAlg_k \to \cAlg_{k/k}\ ,
\end{equation}
which is left adjoint to the current Lie algebra functor $\,
\g: \cAlg_{k/k} \to \LAlg_k\,$,$\ A \mapsto \g(\bar{A})\,$.
We call \eqref{repcl} the {\it representation functor} in $ \g $. Geometrically,  one can think of $ (\mfa_{\g}, \varepsilon)  $ as the coordinate ring $ k[\Rep_{\g}(\mfa)] $ of the based affine scheme $ \Rep_{\g}(\mfa) $, with the basepoint corresponding to the trivial representation.

The adjoint functors $ ((\,\mbox{--}\,)_{\g},\,\g) $ extend naturally to the categories of differential-graded (DG) algebras:
\begin{equation}
\la{drepg}
(\,\mbox{--}\,)_{\g}\,:\ 
\DGL_k \,\rightleftarrows \,\cDGA_{k/k}:\, \g 
\end{equation}
It is well known \cite{Q1} that the categories $ \DGL_k $ and $ \cDGA_{k/k} $ carry natural (projective) model structures, where the weak equivalences (resp., fibrations) are the
quasi-isomorphisms (resp., degreewise surjective maps) of DG algebras. It is shown in \cite{BFPRW} that \eqref{drepg} is a Quillen adjunction with respect to these model structures.  Hence, although the representation functor $ (\,\mbox{--}\,)_{\g}  $ is not homotopy invariant (it does not preserve quasi-isomorphisms), it is left Quillen and therefore has a well-behaved left derived functor 
\begin{equation}
\la{drepg1}
\L(\,\mbox{--}\,)_{\g}\,:\,\Ho(\DGL_k) \rar \Ho(\cDGA_{k/k})\ .
\end{equation}
%
For a given DG Lie algebra $ \mfa $, we now define the {\it representation
homology} of $ \mfa $ in $ \g $ by
$$
\HR_{\ast}(\mathfrak{a},\g)\,:=\, \L_\ast(\mathfrak{a})_{\g} \ ,
$$
where $  \L_\ast(\,\mbox{--}\,)_{\g} := \H_\ast[\L(\,\mbox{--}\,)_{\g}] $ denotes the composition of \eqref{drepg1} with the homology functor on  $\cDGA_{k/k} $. By definition,  $ \HR_\ast(\mfa, \g) $ is a graded commutative $k$-algebra, which depends on $\g$ and (the homotopy type of) the DG Lie algebra $ \mfa $. If $ \mfa \in \LAlg_k $ is an ordinary Lie algebra, there is a natural isomorphism $ \H_0(\mathfrak{a},\g) \cong \mfa_{\g} $ which justifies our definition for the derived representation scheme of $ \Rep_\g(\mfa)$: 
$$ 
\DRep_{\g}(\mfa) := \boldsymbol{\mathrm{R}} \Spec[\L_\ast(\mathfrak{a})_{\g}] 
$$

Now, let $G$ be an affine algebraic group over $k$ associated with the Lie algebra $\g$. Observe that for any
$ \mfa \in \LAlg_k $, $ G $ acts naturally on $ \mfa_{\g} $ by
automorphisms: this action is algebraic and functorial in $ \mfa $. We write $\,(\,\mbox{--}\,)_{\g}^{G}\,:\,
\LAlg_k \to \cAlg_{k/k} $ for the subfunctor of
$ (\,\mbox{--}\,)_{\g} $ defined by taking the $G$-invariants:
$$
\mfa_{\g}^{G} := \{x \in \mfa_{\g}\,:\,g(x) = x\, ,\, \forall\,g\in G\}\ .
$$
The algebra $ \mfa_{\g}^{G} $ represents the affine
quotient scheme $\,\Rep_{\g}(\mfa)/\!/G \,$ parametrizing the closed orbits of $ G $ in $\,\Rep_{\g}(\mfa) \,$.
Although it is not, in general, left Quillen, the functor $\,(\,\mbox{--}\,)_{\g}^{G} $ also admits the (total) left derived functor
$$
\L(\,\mbox{--}\,)_{\g}^{G}\,:\,\Ho(\DGL_k) \rar \Ho(\cDGA_{k/k})\ ,
$$
and we can consider the associated homology functor $  \L_\ast(\,\mbox{--}\,)^G_{\g} := \H_\ast[\L(\,\mbox{--}\,)^G_{\g}]\, $ (cf. \cite[Theorem~2.6]{BKR}). Then, if the algebraic group $G$ is reductive over $k$, there is a natural isomorphism
$$
\L_\ast(\mfa)^G_{\g} \, \cong \, \HR_{\ast}(\mathfrak{a},\g)^{G} .
$$
where $ \HR_{\ast}(\mathfrak{a},\g)^{G} $ denotes the invariant part of the representation homology of $ \mfa $ in $ \g $.
\subsection{A functor homology interpretation} 
\la{proof2}
Recall from the Introduction that $\mathfrak{G}$ denotes the full subcategory of $\mathtt{Gr}$ whose objects are the free groups $\langle n \rangle$ based on the sets $\underline{n}:=\{1,2,\ldots,n\}$, $\,n \geq 0$. We write $ \mathfrak{G}\text{-}\Mod $ (resp., $\Mod\text{-}\mathfrak{G}$) for the categories 
of all covariant (resp., contravariant) functors on $\mathfrak{G}$ with values in the category of $k$-vector spaces. Since $\mathfrak{G} $ is a small category, the categories $ \mathfrak{G}\text{-}\Mod $  and $\Mod\text{-}\mathfrak{G}$ are both abelian with sufficiently many injective and projective objects. We view (and refer to) the objects of $\mathfrak{G}\text{-}\Mod$ and $\Mod\text{-}\mathfrak{G}$ as left and right $\mathfrak{G}$-modules, respectively.

There is a natural bifunctor called the functor tensor product (see, e.g., \cite[Appendix C]{L}):
\[
\mbox{--} \otimes_{\mathfrak{G}} \mbox{--} \,:\, \Mod\text{-}\mathfrak{G} \times \mathfrak{G}\text{-}\Mod \rar \ve_k
\]
This bifunctor is right exact with respect to each argument, preserves sums, and is left balanced. By classical homological algebra \cite{CE}), the derived functors of $\mbox{--} \otimes_{\mathfrak{G}} \mbox{--}$ with respect to each argument are thus isomorphic, and we denote their common value by $\Tor^{\mathfrak{G}}_{\ast}(\,\mbox{--}\,,\,\mbox{--}\,)$. 
Now, as explained in the Introduction, every commutative Hopf algebra $ {\mathcal H} $ defines a left $\mathfrak{G}$-module $\underline{\mathcal H}:  \langle n \rangle \mapsto {\mathcal H}^{\otimes n} $, and dually, every cocommutative Hopf algebra 
$ U $ defines a right $\mathfrak{G}$-module $ \underline{U}: \langle n \rangle \mapsto U^{\otimes n} $. Abusing the notation we will often omit the `underline' in the above formulas, indentifying the $ \mathfrak{G}$-modules $ \underline{\mathcal H} $ and $ \underline{U} $  with the corresponding Hopf algebras $ {\mathcal H} $ and $U$.

\bthm 
\la{liereps} Let $ G $ be an affine algebraic group with coordinate ring
$ \cO(G) $ and the associated Lie algebra $ \g $. Then, for any Lie algebra $\mfa\,\in\,\mathtt{Lie Alg}_k $, there is a natural isomorphism
$$
\HR_\ast(\mfa,\g)\,\cong\,\mathrm{Tor}^{\mathfrak{G}}_\ast(U\mfa,\,
\cO(G))\ ,
$$
where $ \cO(G) $ and $ \, U\mfa  $ are equipped with the standard Hopf algebra
structures $($commutative and cocommutative, respectively$)$.
\ethm
Theorem~\ref{liereps} follows from Lemma~\ref{ualphatensorog},
which is a simple formal result (probably well known to experts: see, for example, \cite{KP}), and Proposition \ref{tors} --- an apparently much deeper result on functor homology ---  whose proof involves topological arguments.

\blemma \la{ualphatensorog}
For any $\mfa\,\in\,\LAlg_k$, there is a natural isomorphism of commutative algebras
$$U\mfa \otimes_{\mathfrak{G}} \cO(G) \,\cong\, \mfa_\g\,,$$
where $\mfa_\g$ is the representation algebra defined in \eqref{ag}.
\elemma

\bproof
Let $B\,\in\, \cAlg_k$. From the left $\mathfrak{G}$-module $\cO(G)$, we form the right $ \mathfrak{G}$-module $\Hom_k({\cO(G)}, B)$, which assigns $\Hom_k(\cO(G)^{\otimes m}, B)$ to $\langle m \rangle$. Since $B$ is a commutative $k$-algebra and since $\cO(G)$ is a strictly monoidal left $\mathfrak{G}$-module, $\Hom_k(\cO(G), B)$ acquires the structure of a lax monoidal right $\mathfrak{G}$-module. This structure is given by the maps
$$\begin{diagram}[small] \Hom_k(\cO(G)^{\otimes m}, B) \otimes \Hom_k(\cO(G)^{\otimes n}, B) & \rTo^{\mu_B \circ (\mbox{--}\,\otimes\,\mbox{--})} & \Hom_k(\cO(G)^{\otimes (m+n)}, B) \end{diagram}\,, $$
where $\mu_B$ is the product on $B$. By the standard $\Hom - \otimes$ adjunction, there is a natural isomorphism of $k$-vector spaces
$$ \Hom_k(U\mfa \otimes_{\mathfrak{G}} \cO(G), B) \,\cong\, \Hom_{\mathtt{Mod}\text{-}\mathfrak{G}}(U\mfa, \Hom_k({\cO(G)}, B)) \ .$$
It is routine to check that under this isomorphism, the  $k$-algebra homomorphisms from $U\mfa \otimes_{\mathfrak{G}} \cO(G)$ to $B$ correspond to
 the right $\mathfrak{G}$-module homomorphisms from $U\mfa$ to $\Hom_k({\cO(G)}, B)$ that respect the (lax) monoidal structure. Since $\cO(G)$ is a coalgebra and $B$ is an algebra,
 $\Hom_k(\cO(G),B)$ has an algebra structure (with product given by convolution). Another routine verification shows that the set of right $\mathfrak{G}$-module homomorphisms from
 $U\mfa$ to $\Hom_k({\cO(G)}, B)$ that respect the (lax) monoidal structure is in (natural) bijection with the set of $k$-algebra homomorphisms $\varphi$ from
 $U\mfa$ to $\Hom_k(\cO(G),B)$ that satisfy the following additional conditions:
$$\varphi(x)(fg) \,=\, \varphi(x^{(1)})(f)\varphi(x^{(2)})(g)\,,\,\, \varphi(x)(1_{\cO(G)})\,=\,\varepsilon(x)1_B\,,\,\,\varphi(Sx)(f)\,=\,\varphi(x)(Sf)\,  $$
for all $x \,\in\,U\mfa$ and $f,g\,\in\,\cO(G)$. Here, $\varepsilon$ and $S$ stand for the counit and antipode of $U\mfa$ respectively the coproduct in $U\mfa$ is given by  $x \mapsto x^{(1)} \otimes x^{(2)}$ in Sweedler notation. It is not difficult to verify that the third condition above follows from the first two. As shown in~\cite[Example 3.4]{MT}, the algebra homomorphisms from $U\mfa$ to $\Hom_k(\cO(G),B)$ satisfying the above conditions are in natural bijection with Lie algebra homomorphisms from $\mfa$ to $\g({B})$. Indeed, $\varphi$ satisfies these conditions for all $x$ in $U\mfa$ iff it satisfies these conditions for $x\,\in\,\mfa$. For $x\,\in \,\mfa$, these conditions are equivalent to the assertion that $\varphi(x)$ is a $k$-linear derivation on $\cO(G)$ with respect to the homomorphism $1_B \circ \varepsilon_{\cO(G)}$, where $\varepsilon_{\cO(G)}$ denotes the canonical augmentation on $\cO(G)$. Such derivations are indeed in bijection with elements of $\Hom_k(\g^{\ast},B) \cong \g(B)$.
We thus, have a natural bijection
$$ \Hom_{\cAlg_k}(U\mfa \otimes_{\mathfrak{G}} \cO(G), B) \,\cong\, \Hom_{\LAlg_k}(\mfa,\g(B))\ .$$
The desired lemma now follows from the Yoneda lemma.
\eproof

\bprop
\la{tors} Let $V$ be a $k$-vector space, and let $ \underline{TV} $ be the right $\mathfrak{G}$-module associated to the tensor algebra of $V$ equipped with 
the standard cocommutative Hopf algebra structure. Then
$$
\Tor_i^{\mathfrak{G}}(\underline{TV},\, \cO(G)) \cong
\left\{
\begin{array}{lll}
\Sym(\g^{\ast} \otimes V) & \mbox{\rm if}\ i=0\\*[1ex]
0 & \mbox{\rm if}\ i > 0
\end{array}
\right.
$$
In particular,  $\, \Tor_i^{\mathfrak{G}}[(\underline{TV})_q, \cO(G)] =0\,$ for all $ i> 0$ and $q \geq 0$.
\eprop
Our proof of Proposition~\ref{tors} is based on topological arguments: specifically, it uses Theorem~\ref{dtpliemodel1} (and its Corollary \ref{dtpliemodel}) as well as our earlier computations of the representation homology of wedges of spheres (\cite[Proposition 5.3]{BRYI}). We do not know a completely algebraic proof of this result.

\bproof[Proof of Proposition \ref{tors}]
Note that the cocommutative Hopf algebra $TV$ can be viewed as the universal enveloping algebra $U(LV)$ of the free Lie algebra generated by $V$.  The corresponding module $\underline{TV}$ has a weight grading induced by the weight grading on $TV$ in which $V$ has weight $1$. Let $(\underline{TV})_q$ denote the component of $\underline{TV}$ of weight $q$. For example, $\underline{V}:= \underline{TV}_1$ is the $\mathfrak{G}$-module defined by $ \lin_{k}^* \otimes V $, where 
$ \lin_\Q $ is the linearization functor ({\it cf.} \cite[Ex. 3.1]{BRYI})
\begin{equation*}
\la{linfun}
\lin_k:\,  \ffgr \to \ve_k\ ,\quad \frgr{n} \mapsto  \
\frgr{n}_{\rm ab} \otimes_{\Z} k = k^{n}\ .
\end{equation*}
Since $\underline{TV} \otimes_{\mathfrak{G}} \cO(G)\,\cong\, U(LV) \otimes_{\mathfrak{G}} \cO(G)$, the required isomorphism for $ i = 0$ follows from Lemma \ref{ualphatensorog}.
To prove the vanishing of $ \Tor_i^{\mathfrak{G}}(\underline{TV},\, \cO(G)) $ for $ i>0$, we assign  $V$ (homological) degree $2$. Then $\underline{TV}$ is a graded right $\mathfrak{G}$-module, whose component in degree $2q$ is $(\underline{TV})_q$. Thus,
$$ \H_n[\underline{TV} \otimes^{\L}_\mathfrak{G} \cO(G)] \,\cong\, \Moplus_{2q+i=n} \Tor^i_{\mathfrak{G}}[(\underline{TV})_q, \cO(G)]\,\text{.}  $$
The desired proposition will follow once we show that
\begin{equation} \la{dreptv}\H_\ast[\underline{TV} \otimes^{\L}_\mathfrak{G} \cO(G)] \,\cong\, \Sym(\g^\ast \otimes V) \,\text{.}\end{equation}
By Theorem~\ref{NNstar_equiv} there are Quillen equivalences refining the Dold-Kan correspondence
$$N^{\ast}\,:\, \DGL^{+}_{k} \rightleftarrows \mathtt{sLie}_{k}\,:N\,, \qquad N^{\ast}\,:\, \DGA^{+}_{k} \rightleftarrows \mathtt{sAlg}_{k}\,:N\, , $$
where $ \mathtt{s}{\mathscr C}$ denotes the category of simplicial objects in a category $ {\mathscr C} $. Equip $N^{\ast}TV \,\cong\, T(N^{-1}V)$ (see formula \ref{Nstar_semifree_formula} in Appendix \ref{doldkan}) with the simplicial cocommutative Hopf algebra structure given by its identification with $ U L(N^{-1}V)$. This gives $N^{\ast}TV$ the structure of a simplicial right $\mathfrak{G}$-module (which we denote by $\underline{N^{\ast}TV}$). This module assigns to the free group $\langle m \rangle$ the simplicial vector space $N^{\ast}TV^{\otimes m}$. Since $V$ has degree $2$, $N^{\ast}LV \,\cong\,L(N^{-1}V)$ is a semi-free simplicial Lie model for the space $X$ given by the wedge of $(\dim_k V)$ copies of the $3$-spheres $ \bS^3 $. By Theorem~\ref{dtpliemodel1} (in particular, Corollary \ref{dtpliemodel} thereof) and \cite[Proposition~5.3]{BRYI}, we then conclude
$$ \H_\ast[N(\underline{N^{\ast}TV}) \otimes^{\L}_\mathfrak{G} \cO(G)] \,\cong\, \HR_\ast(X,G) \,\cong\, \Sym(\g^{\ast} \otimes V) \,\text{.}$$
To complete the proof it remains to note that the natural map $\,\underline{\varepsilon}\,:\,\underline{TV} \rar N(\underline{N^{\ast}TV})\,$ (induced by the unit of the adjunction between the functors $N$ and $N^*$) is a quasi-isomorphism of right $\mathfrak{G}$-modules.
Indeed, $\, \underline{\varepsilon}\,$ is defined by the family of maps
$$
\underline{\varepsilon}(\langle m \rangle):\ TV^{\otimes m}  \xrightarrow{\epsilon^{\otimes m}}  N(N^{\ast}TV)^{\otimes m} \to  N(N^{\ast}TV^{\otimes m})\ ,
$$
where the last arrow is the Eilenberg-Zilber map (which is well-defined for $m>2$ because of the associativity of the Eilenberg-Zilber map for $m=2$). That this is a quasi-isomorphism follows from the K\"{u}nneth Theorem and the fact that $\varepsilon\,:\, TV \rar N(N^{\ast}TV)$  is a quasi-isomorphism of algebras. The associativity of the Eilenberg-Zilber map implies that the maps $\underline{\varepsilon}(\langle m\rangle)$ indeed assemble into a morphism of right $\mathfrak{G}$-modules.
\eproof
\begin{remark}
The result of Proposition \ref{tors} extends to   (homologically) graded vector spaces. To be precise, if $ V = V_\ast $ be a DG $k$-module with trivial differential, such that $V_i = 0$ for all $ i \ll 0 $, then there is an isomorphism in the derived category of $k$-modules:
$$
\underline{TV} \otimes^{\L}_\mathfrak{G} \cO(G) \cong  \Sym(\g^\ast \otimes V) \ .
$$
As a result, there is a homology spectral sequence of the form
\begin{equation} \la{sstv} E^2_{pq}\,=\,\Tor_p^{\mathfrak{G}}(\underline{\H}_q(TV)\,,\,\cO(G)) \,\implies\, \Sym(\g^{\ast} \otimes V)\ , 
\end{equation}
where $\underline{\H}_q(TV)$ stands for the component of the right $\mathfrak{G}$-module $\underline{TV}$ in homological degree $q$ (note that this module is in general {\it different} from $(\underline{TV})_q$). 
Now, if we take $ V = \rH_\ast(X,k) $, the (reduced) homology
of some pointed space $X$, then the spectral sequence \eqref{sstv} has a topological meaning: it is isomorphic to  the fundamental spectral sequence of \cite[Theorem 4.3]{BRYI} for the reduced suspension $\Sigma X$:
$$ E^2_{pq}\,=\,\Tor_p^{\mathfrak{G}}(\underline{\H}_q(\Omega \Sigma X)\,,\,\cO(G))\,\implies\, \HR_\ast(\Sigma X, G)\ .
$$
Indeed, by \cite[Prop. 5.3]{BRYI}, $\HR_\ast(\Sigma X,G)\,\cong\,\Sym(\g^{\ast} \otimes V)$. On the other hand, by the classical Bott-Samelson Theorem \cite{BS}, we have
an isomorphism of graded Hopf algebras: $\,\H_\ast(\Omega \Sigma X,k)\,\cong\, TV$; the latter gives isomorphisms of right  $\mathfrak{G}$-modules: $\, \underline{\H}_q(\Omega \Sigma X,k)\,\cong\,\underline{\H}_q(TV)\,$ for all $\,q \ge 0 \,$.
\end{remark}

\bproof[Proof of Theorem~\ref{liereps}]
By Theorem~\ref{NNstar_equiv} (see Appendix), there are Quillen equivalences refining the classical Dold-Kan correspondence
$$N^{\ast}\,:\, \DGL^{+}_{k} \rightleftarrows \mathtt{sLie}_{k}\,:N\,, \qquad N^{\ast}\,:\, \cDGA^{+}_{k} \rightleftarrows \mathtt{scAlg}_{k}\,:N\, , $$
where $ \mathtt{s}{\mathscr C}$ denotes the category of simplicial objects in a category $ {\mathscr C} $. Let $\mathcal{L} \stackrel{\sim}{\rar} \mfa$ be a semi-free DG resolution of $\mfa$. Let $L:=N^{\ast}(\mathcal L)$. By Theorem~\ref{NNstar_equiv}, $L \stackrel{\sim}{\rar} \mfa$ is a cofibrant resolution in $\mathtt{sLie}_k$. Since the representation functor $(\,\mbox{--}\,)_\g$ is  left adjoint, it commutes with $N^{\ast}$,  i.e, there is a commutative diagram of functors
$$
\begin{diagram}[small]
\DGL_{k}  & \rTo^{N^{\ast}} & \mathtt{sLie}_k\\
  \dTo^{(\,\mbox{--}\,)_\g} & & \dTo_{(\,\mbox{--}\,)_\g}\\
 \cDGA_k & \rTo^{N^{\ast}} & \mathtt{scAlg}_k\\
 \end{diagram}\ .
 $$
Thus, $\HR_\ast(\mfa,\g) \,\cong\,\H_\ast[\mathcal{L}_\g]\,\cong\, \pi_\ast[L_\g]$.
By Lemma \ref{ualphatensorog}, $L_\g \,\cong\, \underline{U\mathcal L} \otimes_{\mathfrak{G}} \cO(G)$. Since $L$ is semi-free by Proposition \ref{Nstar_semifree_prop}, the right $\mathfrak{G}$-module of $n$-simplices in the simplicial right $\mathfrak{G}$-module $\underline{U\mathcal L}$ is of the form $\underline{TV}$ for some vector space $V$. It follows from Proposition~\ref{tors} that the map $\mathrm{C}(\underline{U\mathcal L}) \otimes^{\L}_{\mathfrak{G}} \cO(G) \rar  \mathrm{C}(\underline{U\mathcal L}) \otimes_{\mathfrak{G}} \cO(G)$ is a quasi-isomorphism, where $\mathrm{C}(\mbox{--})$ stands for associated chain complex. The desired result then follows once we establish that $\underline{U\mathcal L}$ is a simplicial resolution of $\underline{U\mfa}$. For this, we need to check that for any $m$, $U\mathcal L^{\otimes m}$ resolves $U\mfa^{\otimes m}$. This follows from the Eilenberg-Zilber and K\"{u}nneth Theorems.
\eproof

\section{The Main Theorem} 
\la{sect6}
In this section, we prove Theorem \ref{hrsullivan} stated in the Introduction. We deduce this result from Theorem~\ref{conj1} which we call the Comparison Theorem. Despite its simple appearance, this theorem is a nontrivial result, the proof of which relies heavily on Quillen's theory \cite{Q} and requires a number of technical refinements
thereof. As these refinements may be useful for other applications, we
state them carefully and prove in a detailed manner.

\subsection{Comparison Theorem}
In this section for simplicity, we assume that $k=\Q$ to use directly results from Quillen's rational homotopy paper~\cite{Q}. However, as explained in Remark~\ref{rem636}, the results of this section extend to an arbitrary field of characteristic $0$ by a universal coefficient argument.

Let $X$ be a 1-connected 
topological space of finite rational type.
Recall ({\it cf.} \cite{FHT1}) that one can associate to $ X $ a commutative
cochain DG algebra $ {\mathcal A}_X $, called a {\it Sullivan model} of $X$, and a connected (chain) DG Lie algebra $ \mfa_X $, called a {\it Quillen model} of $X$. Each of these algebras is uniquely determined up to homotopy and each encodes the rational homotopy type of $X$. The relation between them is given by a DG algebra quasi-isomorphism
\begin{equation}
\la{SQ}
\C^{\ast}(\mfa_X;\mathbb Q) \stackrel{\sim}{\to} {\mathcal A}_X\ ,
\end{equation}
where $\, \C^{\ast}(\mfa_X;\Q) \,$ is the Chevalley-Eilenberg cochain complex of $ \mfa_X $. The homology of $ \mfa_X
$ is the homotopy Lie algebra $ L_X = \pi_*(\Omega X)_\Q $, while the cohomology of  $ {\mathcal A}_X   $ is the rational cohomology algebra $ \H^*(X; \Q) $ of $X$. Among Quillen models of $X$, there is a {\it minimal} one  given by a semi-free DG Lie algebra $  (\mathcal{L}(V), d) $ generated by a graded $\Q$-vector space $V$
with differential $d$ satisfying $d(V) \subset [\,\mathcal L(V), \,\mathcal L(V)]\,$. Such a minimal model  is determined uniquely up to (noncanonical) isomorphism. In particular, $ V \cong \rH_{\ast}(X;\Q)[-1]$  (see  \cite{FHT1}, p. 326).

Now, given an algebraic group $G$, one can associate to a 1-connected space $X$  two representation homologies: the representation homology  $\mathrm{HR}_{\ast}(X,G) $ of $ X $ with coefficients in $G$ and the representation homology  $ \HR_{\ast}(\mfa_X,\g) $ of a Lie model $ \mfa_X $ of $X$
with coefficients in the Lie algebra of $G$ (in sense of Section~\ref{DREP}). A priori, these two representation homologies are defined in a very different way, but the following Comparison Theorem shows that they actually agree.

\bthm[Comparison Theorem] \la{conj1}
For any affine algebraic group $G$ with Lie algebra $\g$, there is an isomorphism of graded commutative $\Q$-algebras
$$ \mathrm{HR}_{\ast}(X,G)\,\cong\, \HR_{\ast}(\mfa_X,\g)\text{.}$$
\ethm

Theorem \ref{hrsullivan} follows readily from Theorem~\ref{conj1} modulo some general algebraic results on representation homology of Lie algebras proved in \cite{BFPRW}.
\bproof[Proof of Theorem \ref{hrsullivan}]
Since the Sullivan model of $X$ is uniquely determined up-to homotopy, it suffices to prove the desired theorem for a particular choice of Sullivan model of $X$. Let $\mfa_X:= (\mathcal{L}(V),d)$ be the minimal Quillen model of $X$. Then, $\mfa_X$ is connected, i.e, concentrated in positive homological degree and finite dimensional in each homological degree. Hence, $C\,:=\,\C_\ast(\mfa_X;\Q)$ is $2$-connected (i.e, its coaugmentation coideal is concentrated in degrees $\geq 2$) and finite dimensional in each homological degree.  The graded $\Q$-linear dual of $C$ is $\cA_X\,:= \,\C^\ast(\mfa_X;\Q)$, which is a Sullivan model of $X$.  Moreover, $C$ is Koszul dual to $\mfa_X$\footnote{It is well known that if $\mathrm{char}(k)=0$, there is a Quillen equivalence $ \cb_{\mathtt{comm}}\,:\, \cDGC_{\Q/\Q} \rightleftarrows \DGL_{\Q}\,:\, \C_\ast(\,\mbox{--}\,;\Q)$ between the category of $\cDGC_{k/k}$ of (coaugmented, conilpotent) DG Lie coalgebras and the category $\DGL_k$ of DG Lie algebras. 
Thus, there is a quasi-isomorphism of DG Lie algebras $\cb_{\mathtt{comm}}(C) \stackrel{\sim}{\rar} \mfa_X$, which shows that $C$ is Koszul dual to $\mfa_X$.
}.  It follows from Theorem \ref{conj1} and \cite[Theorem 6.7 (b)]{BFPRW} (also see {\it loc. cit.}, Theorem 6.3 and the subsequent remark) that
$$ \mathrm{HR}_{\ast}(X,G)\,\cong\, \HR_\ast(\mfa_X,\g) \,\cong\, \H^{-\ast}(\g(\bar{\cA}_X);\Q)\ . $$
If, moreover,  $G$ is reductive, we have
$$  \mathrm{HR}_{\ast}(X,G)^G\,\cong\, \H^{-\ast}(\g(\bar{\cA}_X);\Q)^G \,\cong\, \H^{-\ast}(\g(\bar{\cA}_X);\Q)^{\ad\,\g} \,=\,  \H^{-\ast}(\g({\cA_X}), \g;\Q) \ . $$
The first isomorphism above follows from the fact that all (quasi-)isomorphisms in the proof of Theorem~\ref{conj1} are $G$-equivariant. Indeed, every $G$-action involved is induced by the $G$-action on the left $\mathfrak{G}$-module $\cO(G)$ coming from the conjugation action of $G$ on itself. This finishes the proof of the theorem.
\eproof

Before proving Theorem~\ref{conj1}, we record one useful consequence that gives an explicit DG algebra model for the representation homology of $X$ in terms of the minimal Quillen model of $X$.

\bcor \la{cancomp}
Let $\mfa_X\,=\,(\mathcal L(V),d)$ be the minimal Quillen model of $X$. Then,  $ (\mfa_X)_\g $  is a canonical DG $\Q$-algebra whose homology is isomorphic to $\HR_\ast(X,G)$. Thus, as graded algebras,
$$ \HR_\ast(X,G) \,\cong\, \H_\ast[\Sym_k(\g^{\ast} \otimes V), \partial]\,,$$
where the differential $ \partial $ is given on generators by
$$\partial(\xi^{\ast} \otimes v) \,=\, \langle \xi^{\ast}, \varrho(dv) \rangle\,,\,\,\forall \,\xi \in \g^{\ast}, \,v \in V\,, $$
where $\varrho\,:\, \mathcal L(V) \rar \Sym(\g^{\ast} \otimes V) \otimes \g$ is the universal representation~\eqref{univrep}.
\ecor

\bproof
Since $\mfa_X$ is a semi-free (hence, cofibrant) DG Lie algebra, $\HR_\ast(\mfa_X,\g) \,\cong\, \H_\ast[(\mfa_X)_\g]$. The first assertion is then immediate from Theorem~\ref{conj1}. The algebra isomorphism $ (\mfa_X)_\g \,\cong\,\Sym_k(\g^{\ast} \otimes V)$ follows easily from formula \eqref{ag}. The formula for the differential $\partial$ can follows easily from the fact that the universal representation $\rho\,:\,\mfa_X \rar (\mfa_X)_\g \otimes \g$ is a {\it differential graded} Lie algebra homomorphism.
\eproof

\noindent
\textbf{Example 1.} Recall (see Example 5, \cite[Ch. 24]{FHT1}) that the minimal Lie model for the complex projective space $\c\mathbb P^r, r \geq 1$ is given by the free Lie algebra $\mfa_{r}\,:=\,\mathcal L(v_1,v_2,\ldots,v_r)$ generated by $v_1,\ldots,v_r$, where the degree of $v_i$ is $2i-1$, and the differential is defined on generators by $dv_1=0$, $dv_i = \frac{1}{2}\,\sum_{j+k=i} [v_j,v_k]$ for all $i \geq 2$.
By Corollary~\ref{cancomp}, we have
$$\HR_{\ast}(\c\mathbb P^r,G) \,\cong\, \H_\ast[(\mfa_r)_\g] \,\cong\, \H_\ast \big[ \Sym\big( \Moplus_{i=1}^r \g^{\ast} \cdot v_i \big), \partial \big]\,, $$
where $\g^{\ast} \cdot v_i$ denotes a copy of $\g^\ast$ in degree $2i-1$ indexed by $v_i$ and where the differential $d$ is given on generators by
$$ \partial(\xi^{\ast} \cdot v_i) \,=\, \sum_{j+k=i} (\xi^{\ast}_1\cdot v_j)(\xi^{\ast}_2 \cdot v_k) \,\text{.}$$
Here, the cobracket on $\g^{\ast}$ is given by $ \xi^\ast \mapsto \xi^{\ast}_1 \wedge \xi^{\ast}_2$ in Sweedler notation.

\vspace{1.5ex}

\noindent
\textbf{Example 2.}
As another application of Corollary~\ref{cancomp}, we can compute explicitly the representation homology of higher connected spaces in low homological degrees. To be precise, let $ X $ be an $n$-connected space for some $ n \ge 1$. Consider the minimal Quillen model $ \mfa_X\,=\,(\mathcal{L}(V),\partial)$ of $X$. Then $\, V_i \,\cong\,\H_{i+1}(X;\Q) $ for all $i \ge 0 $.  By the Rational Hurewicz Theorem, $  \H_i(X,\Q) \cong  \pi_{i}(X)_{\Q} $ for all $ 1 \leq i \leq 2 n $. Hence, $V_i=0$ for $i  \leq n-1$. Then the (nonzero) elements of $ [V,V] $ must have homological degree $\geq 2n$, and therefore, by minimality of $ \mfa_X $,  $d(V_i)=0$ for  $n \leq i \leq 2n$. The differential $\partial$ on $(\mfa_X)_\g \,=\,\Sym(\g^{\ast} \otimes V)$ then vanishes on chains of degree $\leq 2n$, and Corollary~\ref{cancomp} implies
\begin{equation*}
\la{hi}
\HR_i(X, G) =
\left\{
\begin{array}{cll}
k & \mbox{\rm for} & i = 0\\*[1ex]
0 & \mbox{\rm for} & 1 \leq i < n\\*[1ex]
\H_{i+1}(X; \g^*) & \mbox{\rm for} & n \leq i \leq 2n-1
\end{array}
\right.
\end{equation*}
The above isomorphisms were found by a different method in \cite{BRYI} (see {\it loc. cit.}, Proposition 4.3).

\subsection{Proof of Comparison Theorem} 
\la{pfconj1}
\subsubsection{Outline of the proof} 
\la{sect6.3.1} The proof of Theorem~\ref{conj1} is based on several technical results. Let $L$ be a semi-free simplicial Lie model of $X$ and let $\underline{R}$ is the simplicial right $\mathfrak{G}$-module associated with the simplicial cocommutative Hopf algebra $R:=U L$.

 Recall that Quillen's rational homotopy theory gives a zig-zag of maps
$$\begin{diagram}[small] \Q\lgr{X} & \rTo & \widehat{\Q}\lgr{X} & \rTo^g & \widehat{R} &\lTo & R \end{diagram}\,,$$
where $\widehat{(\,\mbox{--}\,)}$ stands for completion with respect to the canonical augmentation. Here, the map $g$ (which is by no means unique) is a weak equivalence in the model category of simplicial complete cocommutative Hopf algebras ($\mathrm{sCHA}$'s). The first and third arrows in the above zig-zag induce isomorphisms on all homotopy groups (see~\cite[Sec.3, Part I]{Q}). Our first step is to prove Theorem \ref{dtpliemodel1}, which states that the above zig-zag of maps enriches to a zig-zag of weak equivalences of the associated simplicial right $\mathfrak{G}$-modules:
$$\begin{diagram}[small] \underline{\Q\lgr{X}} & \rTo & \underline{\widehat{\Q}\lgr{X}} & \rTo^{\underline{g}} & \underline{\widehat{R}} &\lTo & \underline{R} \end{diagram} \ . $$
This is verified in a series of Propositions in Section~\ref{proof1} using a relatively straightforward extension of the arguments of {\it loc. cit.} The only subtlety here is that the notion of  weak-equivalence in $\mathrm{sCHA}$ is {\it a priori} different from that of a map inducing an isomorphism on all homotopy groups (see~\cite[Sec.4, Part II]{Q}). This makes it necessary to argue that the map on simplicial right $\mathfrak{G}$-modules induced by $g$ indeed  induces isomorphisms on all homotopy groups. We conclude Section \ref{proof1} by noting that Theorem~\ref{dtpliemodel1} and \cite[Theorem~4.2]{BRYI} together imply that $\HR_\ast(X,G)$ is isomorphic to the homology of the {\it derived} tensor product $N(\underline{R}) \otimes^{\L}_{\mathfrak{G}} \cO(G)$ (Corollary \ref{dtpliemodel}). 

In our second step, starting from Corollary~\ref{dtpliemodel}, we proceed to argue in Section~\ref{6.3.4} that $\HR_\ast(X,G)\,\cong\,\HR_\ast(\mfa_X,\g)$ as graded vector spaces. Our argument is a minor modification of the proof of Theorem \ref{liereps}. For this, we first observe that $\HR_\ast(\mfa_X,\g) \,\cong\, \H_\ast[N(L_\g)]\,\cong\,\H_\ast[\mathrm{C}(L_\g)]$, where $\mathrm{C}$ stands for associated chain complex. Now, $L_\g= \underline{R} \otimes_{\mathfrak{G}} \cO(G)$. Hence, $\mathrm{C}(L_\g)\,\cong\, \mathrm{C}(\underline{R}) \otimes_{\mathfrak{G}} \cO(G)$. By Corollary \ref{dtpliemodel}, it suffices to verify that the map $\mathrm{C}(\underline{R}) \otimes^{\L}_{\mathfrak{G}} \cO(G) \rar \mathrm{C}(\underline{R}) \otimes_{\mathfrak{G}} \cO(G)$ induces an isomorphism on homology. The crucial ingredient in this verification is Proposition~\ref{tors} (stated and proved in Section \ref{proof2}), which implies that the $n$-simplices of $\underline{R}$ are right $\mathfrak{G}$-modules whose higher Tor's with $\cO(G)$ vanish. 

Our third and final step is carried out in Section~\ref{6.3.5}, where we show that the isomorphism $\HR_\ast(X,G) \cong \HR_\ast(\mfa_X,\g)$ of graded vector spaces constructed in Section~\ref{6.3.4} is indeed an isomorphism of graded commutative algebras. We do this by exhibiting for any $q\in \N$ a morphism of simplicial commutative algebras inducing the isomorphism $\HR_i(X,G)\,\cong\,\H_i[N(L_\g)]$ for $i \leq q$. To show this, we first note that the canonical filtration (by powers of the augmentation ideal) on $R$ induces a filtration on the right $\mathfrak{G}$-module $\underline{R}$. Then we use a generic connectivity argument due to Curtis \cite[Sec. 4]{Cu} to show that $\pi_q(F^r\underline{R})=0$ for $r>q$ (Proposition~\ref{quotientr}). This allows us to replace $\underline{R}$ with $\underline{R}/F^r\underline{R}\,, r \gg 0$ when computing homologies in degree $\leq q$ of $N(\underline{R}) \otimes^{\L}_{\mathfrak{G}} \cO(G)$ (i.e, $\HR_i(X,G)$ for $i \leq q$). Again as a consequence of Proposition~\ref{tors}, the $n$-simplices of $\underline{R}/F^r\underline{R}$ are right $\mathfrak{G}$-modules whose higher Tor's with $\cO(G)$ vanish. It follows from these facts that the composite map
$$ \begin{diagram}[small] \underline{\Q\lgr{X}} \otimes_{\mathfrak{G}} \cO(G) & \rTo & \underline{\widehat{\Q}\lgr{X}} \otimes_{\mathfrak{G}} \cO(G) & \rTo & \underline{\widehat{R}}  \otimes_{\mathfrak{G}} \cO(G) & \rTo &  \underline{R}/F^r\underline{R} \otimes_{\mathfrak{G}} \cO(G) \end{diagram}$$
induces the isomorphism $\HR_i(X,G) \,\cong\, \pi_i[L_\g]$ for $i \leq q$ (on $\pi_i$'s). It is not difficult to check that the maps above are morphisms of simplicial commutative algebras. This concludes our argument.

\subsubsection{Step I} \la{proof1} 
By Theorem~\ref{NNstar_equiv} there are Quillen equivalences refining the Dold-Kan correspondence
$$N^{\ast}\,:\, \DGL^{+}_{\Q} \rightleftarrows \mathtt{sLie}_{\Q}\,:N\,,\qquad N^{\ast}\,:\, \DGA^{+}_{\Q} \rightleftarrows \mathtt{sAlg}_{\Q}\,:N\, \qquad N^{\ast}\,:\, \cDGA^{+}_{\Q} \rightleftarrows \mathtt{scAlg}_{\Q}\,:N\, , $$
where $ \mathtt{s}{\mathscr C}$ denotes the category of simplicial objects in a category $ {\mathscr C} $. 
By Proposition~\ref{Nstar_semifree_prop}, applying the functor $N^{\ast}$ to a semi-free Quillen model of $X$ gives a reduced semi-free simplicial Lie model of $X$. Let $L:=L_X$ be a reduced semi-free simplicial Lie model of $X$. 
Consider the simplicial cocommutative Hopf algebra $ R :=  U(L) $ as well as the simplicial complete cocommutative Hopf algebra $\widehat{R}\,\cong\, \widehat{U}(L)$ (where the completion is with respect to the canonical augmentation). 
These correspond to the right $\mathfrak{G}$-modules $\underline{R}$ and $\underline{\widehat{R}}$, which assign to $\langle m \rangle$ the the tensor product $R^{\otimes m}$ and the completed tensor product $\widehat{R}^{\widehat{\otimes} m}$, respectively. 
Similarly, the simplicial cocommutative Hopf algebra $\Q\lgr{X}$ and the simplicial complete cocommutative Hopf algebra 
$\widehat{\Q}\lgr{X}$ correspond to the right $\GG$-modules $\underline{\Q \lgr{X}}$ and $\underline{\widehat{\Q}\lgr{X}}$ 
which assign to $\langle m \rangle$ the tensor product $\Q \lgr{X}^{\otimes m}$ and the completed tensor product $\widehat{\Q}\lgr{X}^{\widehat{\otimes} m}$, respectively. 
Recall that the main result from Quillen's rational homotopy theory \cite{Q} is about the existence of a zig-zag of maps
$$\begin{diagram}[small] \Q\lgr{X} &\rTo& \widehat{\Q}\lgr{X} &\rTo^g & \widehat{R} & \lTo & R\end{diagram} $$
of simplicial commutative Hopf algebras such that the first and last arrows induce isomorphisms on all homotopy groups while the map $g$ is a weak-equivalence in the model category of simplicial complete cocommutative Hopf algebras. First, we prove the following extension of this result.
\bthm \la{dtpliemodel1}
There is a zig-zag of weak-equivalences of simplicial right $\mathfrak{G}$-modules
$$  \begin{diagram}[small]\underline{\Q\lgr{X}} &\rTo & \underline{\widehat{\Q}\lgr{X}} & \rTo^{\underline{g}} & \underline{\widehat{R}} &\lTo & \underline{R} \end{diagram}\ .$$
\ethm
\bproof
The desired result is an immediate consequence of Proposition \ref{comp1}, Proposition \ref{comp2} and Proposition \ref{comp3} which we state and prove below. The proofs of these propositions are exercises in Quillen's rational homotopy theory. 
\eproof 

The propositions leading to Theorem \ref{dtpliemodel1}, are in turn based on the following lemma. Let $V$ be a filtered reduced simplicial vector space. Let $\widehat{V}$ denote the completion of $V$ with respect to the given filtration. More generally, for any $m \in \N$, one has the simplicial vector spaces $V^{\otimes m}, \widehat{V}^{\otimes m}$ and $V^{\widehat{\otimes} m} = \widehat{V}^{\widehat{\otimes} m}$, where $V^{\widehat{\otimes} m}$ denotes the completed tensor product $\varprojlim_r (V/F^rV)^{\otimes m}$. Let $\widehat{\Sym}^m(V)$ denote the image in $V^{\widehat{\otimes} m}$ of the symmetrization idempotent $e_m:= \frac{1}{n!} \sum_{\sigma \in S_m} \sigma$. Let $\widehat{\Sym}(V):=\prod_{m=0}^{\infty} \widehat{\Sym}^m(V)$. Recall that a $\pi_\ast$-equivalence (see~\cite{Maj}) is a morphism inducing isomorphisms on all homotopy groups.

\blemma \la{completions}
Suppose that for each $q>0$, $\pi_q(F^rV)=0$ for $r$ sufficiently large. Then,

\vspace{1ex}

\noindent
$(i)$ For each $q>0$, $\pi_q(F^r\widehat{V})=0$ for $r$ sufficiently large.\\
$(ii)$ The map $V^{\otimes m} \rar \widehat{V}^{\widehat{\otimes} m}$ is a $\pi_\ast$-equivalences for all $m$.\\
$(iii)$ The map $\Sym(V) \rar \widehat{\Sym}(V)$ is a $\pi_\ast$-equivalence.
\elemma

\bproof
By a long exact sequence of homotopy groups (LESH) argument, the natural map $\pi_q(V) \rar \pi_q(V/F^rV)$ is an isomorphism for $r$ sufficiently large. Thus, the inverse system $\{\pi_q(V/F^rV)\}$ is eventually constant. Thus, $\lim^1\{\pi_q(V/F^rV)\}=0$. It follows from~\cite[Part I, Prop. 3.8]{Q} that $\pi_q(\widehat{V})\,\cong\,\pi_q(V/F^rV)$ for $r$ sufficiently large. Since $V/F^rV \,\cong\,\widehat{V}/F^k\widehat{V}$, we see that $\pi_q(\widehat{V})\,\cong\, \pi_q(\widehat{V}/F^r\widehat{V})$ for $r$ sufficiently large. Again by a LESH argument, $\pi_q(F^r\widehat{V})=0$ for $r$ sufficiently large. This proves $(i)$.

Moreover, by the Eilenberg-Zilber and K\"{u}nneth Theorems, $\pi_q(V^{\otimes m}) \,\cong\, \pi_q[(V/F^rV)^{\otimes m}]$ for $r$ sufficiently large (since the same is true for $m=1$). It follows that the inverse system $\{\pi_q[(V/F^rV)^{\otimes m}]\}$ is eventually constant. Arguing as for the case when $m=1$, we see that $\pi_q(\widehat{V}^{\widehat{\otimes m}}) \,\cong\, \pi_q[(V/F^rV)^{\otimes m}]$ for $r$ sufficiently large. This proves that the natural map $V^{\otimes m} \rar \widehat{V}^{\widehat{\otimes} m}$ induces an isomorphism on $\pi_q$ for any fixed $q$. This proves $(ii)$.

Since the map $V^{\otimes m} \rar V^{\widehat{\otimes} m}$ is $S_m$-equivariant and since $\Sym^m(V)$ and $\widehat{\Sym}^m(V)$ are the images of the symmetrization idempotent $e_m$ acting on $V^{\otimes m}$ and $V^{\widehat{\otimes} m}$ respectively, the natural map $\Sym^m(V) \rar \widehat{\Sym}^m(V)$ is a $\pi_\ast$-equivalence. Thus, the map $\Sym(V) \rar \oplus_m \widehat{\Sym}^m(V)$ is a $\pi_\ast$-equivalence. Since $V$ is reduced and by $(ii)$, $\pi_q(V^{\widehat{\otimes} r})=0$ for $r >q$ (by the Eilenberg-Zilber and K\"{u}nneth Theorems). It follows that $\pi_q(\oplus_{m \geq r} \widehat{\Sym}^m(V) =0$ for $r>q$. Applying $(ii)$ to $W:=\oplus_m \widehat{\Sym}^m(V)$ with filtration given by $F^rW:=\oplus_{m \geq r} \widehat{\Sym}^m(V)$, we see that the map $\oplus_{m } \widehat{\Sym}^m(V) \rar \widehat{\Sym}(V)$ is a $\pi_\ast$-equivalence. This proves $(iii)$.
\eproof

\bprop
\la{comp1}
The canonical map of $\mathfrak{G}$-modules $\underline{R} \rar \underline{\widehat{R}}$ is a $\pi_\ast$-equivalence.
\eprop

\bproof
It needs to be shown that the map $R^{\otimes m} \rar \widehat{R}^{\widehat{\otimes} m}$ is a $\pi_\ast$-equivalence. By~\cite[Part I, Thm. 3.7]{Q}, for any fixed $q$, $\pi_q(F^rR)$ vanishes for $r$ sufficiently large. Lemma~\ref{completions} $(ii)$ then implies that the map $R^{\otimes m} \rar \widehat{R}^{\widehat{\otimes} m}$ is a $\pi_\ast$-equivalence, as desired.
\eproof

Recall that $\lgr{X}$ denotes the Kan loop group functor applied to a reduced simplicial/cellular model of $X$. Then, $\Q\lgr{X}$ is a simplicial cocommutative Hopf algebra equipped with a canonical augmentation.  The completion $\widehat{\Q}\lgr{X}$ of $\Q\lgr{X}$ with respect to its canonical augmentation is a simplicial complete cocommutative Hopf algebra  ($\mathrm{sCHA}$). $\Q\lgr{X}$ as well as $\widehat{\Q}\lgr{X}$ correspond to simplicial right $\mathfrak{G}$-modules, which we denote by $\underline{\Q\lgr{X}}$ and $\underline{\widehat{\Q}\lgr{X}}$ respectively.

\bprop
\la{comp2}
The map $\underline{\Q\lgr{X}} \rar \underline{\widehat{\Q}\lgr{X}}$ is a $\pi_\ast$-equivalence.
\eprop

\bproof
We need to show that for each $m$, the map $\Q\lgr{X}^{\otimes m} \rar \widehat{\Q}\lgr{X}^{\widehat{\otimes} m}$ is a $\pi_\ast$-equivalence. By Lemma~\ref{completions} $(ii)$, this follows one we verify that for any fixed $q$, $\pi_q(F^r\Q\lgr{X})=0$ for $r$ sufficiently large. This is immediate from~\cite[Thm. 4.72]{Maj}.
\eproof

We recall that the category $\mathtt{sCHA}$ of reduced $\mathrm{sCHA}$'s is a model category, whose cofibrant objects are retracts of semi-free $\mathrm{sCHA}$'s. The definition of semi-free $\mathrm{sCHA}$ is the obvious extension to the simplicial setting of the definition of a free complete cocommutative Hopf algebra: the free comple cocommutative Hopf algebra generated by a vector space $V$ is $\widehat{TV}$, where $V$ is primitive. We now apply Quillen's rational homotopy theory: in~\cite{Q}, Quillen proves several equivalences of homotopy categories (see {\it loc. cit}, pg. 211, Fig. 2) from which it follows that there is an isomorphism in $\Ho(\mathtt{sCHA})$ $\widehat{\Q}\lgr{X}\,\cong\,\widehat{R}$. By Theorem 4.7 of {\it loc. cit}, there is a morphism $g\,:\, \widehat{\Q}\lgr{X} \rar \widehat{R}$ that is a simplicial homotopy equivalence. Denote the corresponding map of right $\mathfrak{G}$-modules by $\underline{g}\,:\, \underline{\widehat{\Q}\lgr{X}} \rar \underline{\widehat{R}} $.

 \bprop \la{comp3}
 $\underline{g}$ is a $\pi_\ast$-equivalence.
 \eprop

\bproof
By~\cite[Part I, Theorem 3.7]{Q} and Lemma~\ref{completions}, the completion map $\widehat{R}^{\otimes m} \rar \widehat{R}^{\widehat{\otimes} m}$ is a $\pi_\ast$-equivalence. Similarly, it can be shown that the map $\widehat{\Q}\lgr{X}^{\otimes m} \rar \widehat{\Q}\lgr{X}^{\widehat{\otimes} m}$ is a $\pi_\ast$ equivalence. To prove the desired lemma, we need to show that $\underline{g}^{\widehat{\otimes} m} \,:\, \widehat{\Q}\lgr{X}^{\widehat{\otimes} m} \rar \widehat{R}^{\otimes m}$ is a $\pi_\ast$ equivalence for each $m$. Since the diagram
$$
\begin{diagram}[small]
 \widehat{\Q}\lgr{X}^{\otimes m} & \rTo &  \widehat{\Q}\lgr{X}^{\widehat{\otimes} m}\\
 \dTo_{g^{\otimes m}} & & \dTo_{g^{\widehat{\otimes} m}}\\
\widehat{R}^{\otimes m} &\rTo& \widehat{R}^{\widehat{\otimes} m}\\
 \end{diagram}
$$
commutes, it suffices (by the Eilenberg-Zilber and K\"{u}nneth Theorems) to show that $g$ is a $\pi_\ast$-equivalence. Let $\mathcal P$ denote the functor of primitive elements. By~\cite[Appendix A, Cor. 2.16]{Q}, there is an isomorphism of simplicial vector spaces $\widehat{\Sym}(\mathcal P \widehat{R}) \stackrel{\sim}{\rar} \widehat{R}$. For the same reason, $\widehat{\Q}\lgr{X}$ is isomorphic to $\widehat{\Sym}(\mathcal P \widehat{\Q}\lgr{X})$ as simplicial vector spaces. Since $\mathcal P \widehat{R}$ is a canonical retract of $\widehat{R}$, $\pi_q(F^r \mathcal P\widehat{R})=0$ for $r$ large enough (since the same holds for $\widehat{R}$) and for the same reason, $\pi_q(F^r \mathcal P \widehat{\Q} \lgr{X})=0$ for $r$ sufficiently large. By Lemma~\ref{completions} $(iii)$, the horizontal arrows in the commutative diagram below are $\pi_\ast$-equivalences.
$$
\begin{diagram}[small]
 \Sym(\mathcal P \widehat{\Q}\lgr{X}) & \rTo &  \widehat{\Sym}(\mathcal P \widehat{\Q}\lgr{X}) \, \cong\, \widehat{\Q}\lgr{X}\\
 \dTo_{\Sym(\mathcal P g)}  & & \dTo_{g}\\
 \Sym(\mathcal P \widehat{R}) &\rTo& \widehat{\Sym}(\mathcal P \widehat{R}) \,\cong\, \widehat{R}\\
 \end{diagram}
$$
By~\cite[Part II, Theorem 4.7]{Q}, $\mathcal P g$ is a $\pi_\ast$-equivalence. Thus, the left vertical arrow in the above diagram is a $\pi_\ast$-equivalence. It follows that $g$ is a $\pi_\ast$-equivalence, as desired.
\eproof

The following corollary of Theorem \ref{dtpliemodel1} completes the first step towards proving Theorem \ref{conj1}.

\bcor \la{dtpliemodel}
There is an isomorphism of graded vector spaces $\HR_\ast(X,G)\,\cong\,\H_\ast[N(\underline{R}) \otimes^{\L}_{\mathfrak{G}} \cO(G)]$.
\ecor
\bproof
By Theorem \ref{dtpliemodel1}, $N(\underline{\Q\lgr{X}})\,\cong\, N(\underline{R})$ in the derived category of right $\mathfrak{G}$-modules. Hence, there is an isomorphism in the derived category $\mathcal D(\Q)$ of complexes of $\Q$-vector spaces
$$ N(\underline{\Q\lgr{X}}) \otimes^{\L}_{\mathfrak{G}} \cO(G) \,\cong\, N(\underline{R}) \otimes^{\L}_{\mathfrak{G}} \cO(G) \,\text{.}$$
The desired result now follows from \cite[Theorem~4.2]{BRYI}.
\eproof

\subsubsection{Step II}\label{6.3.4}

We now show that  $\HR_\ast(X,G) \,\cong\,\HR_\ast(\mfa_X,\g)$ as graded vector spaces. This step is a minor modification of the proof of Theorem \ref{liereps}. Without loss of generality, we may assume that $\mfa_X$ is semi-free. Since the representation functor $(\,\mbox{--}\,)_\g$ is  left adjoint, it commutes with $N^{\ast}$,  i.e, there is a commutative diagram of functors
$$
\begin{diagram}[small]
\DGL_{\Q}  & \rTo^{N^{\ast}} & \mathtt{sLie}_\Q\\
  \dTo^{(\,\mbox{--}\,)_\g} & & \dTo_{(\,\mbox{--}\,)_\g}\\
 \cDGA_\Q & \rTo^{N^{\ast}} & \mathtt{scAlg}_\Q\\
 \end{diagram}
 $$
Since $N^{\ast} \,:\,\cDGA_{\Q} \rightleftarrows \mathtt{scAlg}_{\Q}\,:\,N$ is a Quillen equivalence, the above commutative diagram implies that $\HR_\ast(\mfa_X,\g) \,\cong\, \H_\ast[N(L_\g)]$ as graded algebras, where $ L:= N^{\ast}\mfa_X$. By Lemma \ref{ualphatensorog}, $L_\g \,\cong\, \underline{R} \otimes_{\mathfrak{G}} \cO(G)$, where $R:= UL$. Thus, $\HR_\ast(\mfa_X,\g) \,\cong\, \H_\ast[\mathrm{C}(\underline{R} \otimes_{\mathfrak{G}} \cO(G))]\,\cong\,\H_\ast[\mathrm{C}(\underline{R}) \otimes_{\mathfrak{G}} \cO(G)] $, where $\mathrm{C}$ stands for associated chain complex (indeed, the inclusion $N(L_\g) \hookrightarrow C(L_\g)$ is a quasi-isomorphism). Since $L$ is a semi-free simplicial Lie algebra by Proposition \ref{Nstar_semifree_prop}, the right $\mathfrak{G}$-module of $n$-simplices in the simplicial right $\mathfrak{G}$-module $\underline{R}$ is of the form $\underline{TV}$ for some vector space $V$. It follows from  Proposition~\ref{tors} that $\mathrm{C}(\underline{R})$ is a complex of right $\mathfrak{G}$-modules whose higher Tor's with $\cO(G)$ vanish. Thus, the map $\mathrm{C}(\underline{R}) \otimes^{\L}_{\mathfrak{G}} \cO(G) \rar \mathrm{C}(\underline{R}) \otimes_{\mathfrak{G}} \cO(G)$ is a quasi-isomorphism. Since there is a quasi-isomorphism of complexes of right $\mathfrak{G}$-modules $N(\underline{R}) \hookrightarrow \mathrm{C}(\underline{R})$, there are isomorphisms of graded vector spaces
$$\HR_\ast(X,G)\,\cong\, \H_\ast[N(\underline{R}) \otimes^{\L}_{\mathfrak{G}} \cO(G)] \,\cong\, \H_\ast[\mathrm{C}(\underline{R}) \otimes^{\L}_{\mathfrak{G}} \cO(G)]\,\cong\, \H_\ast[\mathrm{C}(\underline{R}) \otimes_{\mathfrak{G}} \cO(G)] \,\cong\,\HR_\ast(\mfa_X,\g)\ , $$
where the first isomorphism above is by  Corollary \ref{dtpliemodel}. This completes the second step in the proof of Theorme \ref{conj1}. However, we do not see a resolution of $P \stackrel{\sim}{\rar} \mathrm{C}(\underline{R})$ by projective right $\mathfrak{G}$-modules such that the functor $P\,:\,\mathfrak{G}^{\mathrm{op}} \rar \Com_k$ is monoidal. As a result, we are unable to see the algebra structure on $\H_\ast[\mathrm{C}(\underline{R}) \otimes^{\L}_{\mathfrak{G}} \cO(G)]$ independently of Corollary \ref{dtpliemodel}. We therefore require further work in Section \ref{6.3.5} below to show that the isomorphism $\HR_\ast(X,G)\,\cong\,\HR_\ast(\mfa_X,\g)$ of graded vector spaces is indeed an isomorphism of graded algebras.

\subsubsection{Step III}\label{6.3.5}
To complete the proof of Theorem~\ref{conj1}, it remains to show that $\HR_\ast(X,G) \,\cong\,\HR_\ast(\mfa_X,\g)$ as graded $\Q$-algebras. For this, given any $r\,\in\,\mathbb N$, we shall produce a morphism of simplicial commutative algebras that induces the isomorphism
$ \HR_q(X,G) \,\cong\, \H_q[N(L_\g)]$ for $q<r$.

Recall that $R:= UL$ is a semi-free simplicial associative algebra filtered by powers of its augmentation ideal. This filtration induces a filtration on the simplicial right $\mathfrak{G}$-module $\underline{R}$: if the algebra of $n$-simplices of $R$ is $TV$ for some vector space $V$, then the right $\mathfrak{G}$-module of $n$-simplices of $F^r\underline{R}$ is $\oplus_{ q \geq r} (\underline{TV})_q$. The following connectivity result holds for the filtered right $\mathfrak{G}$-module $\underline{R}$.
\bprop \la{quotientr}
For $r>q$, we have $\pi_q(F^r\underline{R})=0$.
\eprop
\bproof
It needs to be shown that for all $\langle m \rangle$, $\pi_q(F^r\underline{R}(\langle m \rangle))=0$ for $r>q$. For $m=0$, this is obvious. For $m=1$, this is~\cite[Part I, Thm. 3.7]{Q}. For arbitrary $m$, we generalize the argument in {\it loc. cit.}. The functor $\mathtt{Lie}_{\Q} \rar \ve_\Q\,,\,\,L \mapsto F^r\underline{U L}(\langle m \rangle)$ takes $0$ to $0$ and commutes with direct limits.  By~\cite[Remark 4.10]{Cu}, the arguments in Section 4 of {\it loc. cit.} proving Lemma (2.5) therein apply to this functor as well. It therefore, suffices to verify the desired proposition for $R=U\mathfrak{l}$, where $\mathfrak{l}$ is the free simplicial Lie algebra generated by $V:=\overline{\Q}K$, where $K$ is a finite wedge sum of simplicial circles. Note that in this case, $R=TV$, and $V$ is a {\it connected} simplicial vector space. In this case, $F^r\underline{R}(\langle m \rangle) = \oplus_{r_1+ \ldots+r_m \geq r} V^{\otimes r_1} \otimes \ldots \otimes V^{\otimes r_m}$. That $\pi_q$ of each summand vanishes for $q<r$ follows from the Eilenberg-Zilber and K\"{u}nneth Theorems. This proves the desired proposition.
\eproof

\bprop \la{goingdown}
For $r$ sufficiently large, all arrows in the following commutative diagram induce isomorphisms on the homology groups $\H_i[\,\mbox{--}\,]\,,\,i \leq q$.
$$
\begin{diagram}[small]
\mathrm{C}(\underline{R}) \otimes^{\L}_{\mathfrak{G}} \cO(G) & \rTo & \mathrm{C}(\underline{R}/F^r\underline{R}) \otimes^{\L}_{\mathfrak{G}} \cO(G)\\
  \dTo & & \dTo \\
 \mathrm{C}(\underline{R}) \otimes_{\mathfrak{G}} \cO(G) & \rTo & \mathrm{C}(\underline{R}/F^r\underline{R}) \otimes_{\mathfrak{G}} \cO(G)\\
  \end{diagram}$$
\eprop

\bproof
Both $\mathrm{C}(\underline{R})$ and $\mathrm{C}(\underline{R}/F^r\underline{R})$ are complexes of right $\mathfrak{G}$-modules whose higher Tors with $\cO(G)$ vanish by Proposition~\ref{tors}. It follows that the vertical arrows induce isomorphisms on all homology groups. It therefore, suffices to show that the horizontal arrow on top of the above diagram induces isomorphisms on $\H_i[\,\mbox{--}\,]\,,\,i \leq q$ for $r$ sufficiently large.

Consider the good truncation $\tau_{\geq q+1}\mathrm{C}$ (see~\cite[Sec. 1.2.7]{W}) of a chain complex $\mathrm{C}$ of right $\mathfrak{G}$-modules. The exact sequence $0 \rar \tau_{\geq q+1}\mathrm{C} \rar \mathrm{C} \rar \tau_{<q+1}\mathrm{C} \rar 0$ of complexes of right $\mathfrak{G}$-modules gives a distinguished triangle in $\mathcal D(\Q)$ for any right $\mathfrak{G}$-module $N$.

$$\tau_{\geq q+1}\mathrm{C} \otimes^{\L}_{\mathfrak{G}} N \rar \mathrm{C} \otimes^{\L}_{\mathfrak{G}} N \rar \tau_{<q+1}\mathrm{C} \otimes^{\L}_{\mathfrak{G}} N \rar \tau_{\geq q+1}\mathrm{C} \otimes^{\L}_{\mathfrak{G}} N [1]$$

It is easy to see that $\H_i(\tau_{\geq q+1}\mathrm{C} \otimes^{\L}_{\mathfrak{G}} N)=0$ for $i<q+1$. The long exact sequence of homologies associated with the above distinguished triangle then implies that

 \begin{equation} \la{homologies} \H_i(\mathrm{C} \otimes^{\L}_{\mathfrak{G}} N) \cong \H_i(\tau_{<q+1}\mathrm{C} \otimes^{\L}_{\mathfrak{G}} N)\text{ for } i \leq q\,\text{.}\end{equation}

By Proposition~\ref{quotientr}, the map $\tau_{<q+1}\mathrm{C}(\underline{R}) \rar \tau_{<q+1} \mathrm{C}(\underline{R}/F^r\underline{R})$ is a quasi-isomorphism for $r>q$. Thus, the map
$\H_\ast[\tau_{<q+1}\mathrm{C}(\underline{R}) \otimes^L_\mathfrak{G} \cO(G)] \rar \H_\ast[\tau_{<q+1} \mathrm{C}(\underline{R}/F^r\underline{R}) \otimes^L_\mathfrak{G} \cO(G)]$ is an isomorphism of graded $\Q$-vector spaces. The desired proposition now follows from~\eqref{homologies}.
\eproof

Note that the filtration on $\widehat{R}$ induces a filtration on the right $\mathfrak{G}$-module $\underline{\widehat{R}}$. Clearly, $\underline{R}/F^r\underline{R} \,\cong\, \underline{\widehat{R}}/F^r{\underline{\widehat{R}}}$. The following corollary is immediate from Proposition~\ref{comp1} and Proposition~\ref{goingdown}.

\bcor \la{gdcor}
For $r$ sufficiently large, all arrows in the following commutative diagram induce isomorphisms on the homology groups $\H_i[\,\mbox{--}\,]\,,\,i \leq q$.
$$
\begin{diagram}[small]
\mathrm{C}(\underline{R}) \otimes^{\L}_{\mathfrak{G}} \cO(G) & \rTo & \mathrm{C}(\underline{\widehat{R}})\otimes^{\L}_{\mathfrak{G}} \cO(G)\\
  \dTo & \ldTo &\\
 \mathrm{C}(\underline{R}/F^r\underline{R}) \otimes^{\L}_{\mathfrak{G}} \cO(G) & &\\
\end{diagram}
$$
\ecor

Recall that there is a weak equivalence between cofibrant objects in $\mathtt{sCHA}$ $g\,:\, \widehat{\Q}\lgr{X} \rar \widehat{R}$ inducing a map of simplicial right $\mathfrak{G}$-modules $\underline{g}$ (see Proposition~\ref{comp3}). Consider the following commutative diagram, where the second arrow on the top and bottom rows is induced by $\underline{g}$.

\begin{equation} \la{maindiag}
\begin{diagram}[small]
\mathrm{C}(\underline{\Q\lgr{X}}) \otimes^{\L}_{\mathfrak{G}} \cO(G) & \rTo & \mathrm{C}(\underline{\widehat{\Q}\lgr{X}}) \otimes^{\L}_{\mathfrak{G}} \cO(G) & \rTo &  \mathrm{C}(\underline{\widehat{R}}) \otimes^{\L}_{\mathfrak{G}} \cO(G) & &\\
 \dTo & & \dTo & & \dTo & \rdTo &\\
\mathrm{C}[\underline{\Q\lgr{X}}  \otimes_{\mathfrak{G}} \cO(G)] & \rTo &  \mathrm{C}[\underline{\widehat{\Q}\lgr{X}} \otimes_{\mathfrak{G}} \cO(G)] & \rTo & \mathrm{C}[ \underline{\widehat{R}} \otimes_{\mathfrak{G}} \cO(G)] & \rTo &  \mathrm{C}[\underline{R}/F^r\underline{R} \otimes_{\mathfrak{G}} \cO(G)]\\
\end{diagram}
\end{equation}

By \cite[Theorem~4.2]{BRYI}, the left vertical arrow in~\eqref{maindiag} induces isomorphisms on all homologies. The two arrows on the top row of~\eqref{maindiag} induce isomorphisms on all homologies by Propositions~\ref{comp2} and~\ref{comp3} respectively. The diagonal arrow induces isomorphisms on $\H_i[\,\mbox{--}\,]\,,\,i \leq q$ for $r$ sufficiently large by Proposition~\ref{goingdown} and Corollary~\ref{gdcor}. An isomorphism $\HR_i(X,G) \,\cong\, \H_i[N(L_\g)]\,,\,i\leq q$ is thus induced on homologies (for sufficiently large $r$) by the composition of the maps on the bottom row of~\eqref{maindiag}. That the composition of maps in the bottom row is a map of DG commutative algebras follows from the fact that each of the maps
$$\underline{\Q\lgr{X}}  \otimes_{\mathfrak{G}} \cO(G) \rar \underline{\widehat{\Q}\lgr{X}} \otimes_{\mathfrak{G}} \cO(G) \rar  \underline{\widehat{R}} \otimes_{\mathfrak{G}} \cO(G) \rar \underline{R}/F^r\underline{R} \otimes_{\mathfrak{G}} \cO(G) $$
is a morphism of simplicial commutative algebras. Indeed, this last fact follows from~\cite[Prop. 3.4]{KP} and the facts that $\cO(G)$ is a lax-monoidal left $\mathfrak{G}$-module, the $n$-simplices of the right $\mathfrak{G}$-modules $\underline{\Q\lgr{X}}, \underline{\widehat{\Q}\lgr{X}}, \underline{\widehat{R}}$ and $\underline{R}/F^r\underline{R}$ are lax-monoidal for each $n$, and the morphisms
$$ \underline{\Q\lgr{X}} \rar \underline{\widehat{\Q}\lgr{X}} \rar \underline{\widehat{R}} \rar \underline{R}/F^r\underline{R} $$
are natural transformations of lax-monoidal functors on $n$-simplices for each $n$. This completes the proof of Theorem~\ref{conj1}.

\subsubsection{Remark.}\la{rem636} The results of this section go through with $\Q$ replaced by any field $k$ of characteristic $0$. Indeed, the proofs of Propositions~\ref{comp1} and~\ref{comp2} work for any such field $k$. For Proposition~\ref{comp3}, we work with a semi-free simplicial Lie model $L$ of $X$ over $\Q$. The corresponding Lie model over $k$ is $L \otimes_\Q k$. The corresponding $\mathrm{sCHA}$ over $k$ is $\widehat{R \otimes_\Q k}$. The $\pi_\ast$-equivalence of $\mathrm{sCHA}$'s (over $\Q$) $f\,:\, \widehat{R} \rar \widehat{\Q}\lgr{X}$ extends to a $\pi_\ast$-equivalence of $\mathrm{sCHA}$'s over $k$ $f\,:\, \widehat{R \otimes_\Q k} \rar \widehat{k}\lgr{X}$. This proves Proposition~\ref{comp3} over $k$. Theorem~\ref{dtpliemodel1} (and hence, Corollary \ref{dtpliemodel}), Proposition~\ref{tors}, Theorem~\ref{conj1} and Proposition~\ref{liereps} can then be proven over $k$ as done above (over $\Q$).

\subsection{A conjecture for non-simply connected spaces} In a series of recent papers \cite{BFMT1, BFMT2, BFMT3}, Buijs, F\'{e}lix, Murillo and Tanr\'{e} associated a free DG Lie algebra model $ (\mathfrak L_X, d) $ to {\it any} finite simplicial complex $X$.
Unlike Quillen models, the DG Lie algebras $ \mathfrak L_X $ are assumed, in general, to be not connected but {\it complete}
with respect to the canonical decreasing filtrarion $\,
{\mathfrak L}^1 \supseteq {\mathfrak L}^2 \supseteq \ldots \,$ defined by $\mathfrak L^1:=\mathfrak L$ and $\mathfrak L^n:= [\mathfrak L, \mathfrak L^{n-1}]$. The $0$-simplices of $X$ correspond to the degree $-1$ generators of $\mathfrak L_X$ that satisfy the Maurer-Cartan equation, while the $n$-simplices of $X$ correspond to generators in degree $n-1$.
For any connected, finite simplicial complex $X$,
the DG Lie algebra $ \mathfrak L_X $ itself is acyclic (i.e., $ \H_*(\mathfrak L_X, d) = 0 $).
The topological information about $X$ is contained in a DG Lie algebra $(\mathfrak L_X, d_v)$ obtained from $ \mathfrak L_X $
by twisting its differential by Maurer-Cartan elements corresponding to the vertices of $X$, i.e.
$\, d_v:=d +[v\,,\,\mbox{--}\,]\,$, where $v$ denotes (the Maurer-Cartan element corresponding to) a vertex of $X$. Now,
the main result of \cite{BFMT1} (see {\it loc. cit.}, Theorem~A) says that, if $X$ is simply-connected, then
$\,(\mathfrak L_X, d_v)\,$ is quasi-isomorphic to a Quillen model of $X$.  This motivates the following conjectural
generalization of our Theorem~\ref{conj1}.

Let $(\mathfrak L_X,d)$ be a complete free DG Lie algebra model associated to a reduced simplicial set $X$. Let $d_v:=d+[v\,,\,\mbox{--}]$ be the twisted differential on $ \mathfrak L_X $ corresponding to the (unique) basepoint of $X$.
Note that $\HR_0[(\mathfrak L_X,d_v),\g]$ has a canonical augmentation $\varepsilon$ corresponding to the trivial (zero) representation.
Let $\widehat{\HR}_\ast[(\mathfrak L_X,d_v),\g]$ denote the adic completion of  ${\HR}_\ast[(\mathfrak L_X,d_v),\g]$ with respect to the augmentation ideal of $\varepsilon$. Similarly, $\HR_0(X,G)$ has a canonical augmentation corresponding to the trivial (identity) representation of $\pi_1(X, v)$. Let $\widehat{\HR}_\ast(X,G)$ denote the corresponding completion of $\HR_\ast(X,G)$.
\bconj \la{conj5}
{\it There is an isomorphism of completed graded $\Q$-algebras}
$$\widehat{\HR}_{\ast}(X,G)\,\cong\, \widehat{\HR}_\ast[(\mathfrak L_X,d_v),\g] \,\text{.}
$$
\econj
\noindent
Note that Conjecture~\ref{conj5} holds for $X$ simply-connected: indeed, in this case, $(\mathfrak L_X,d_v)$ is quasi-isomorphic to a Quillen model $\mfa_X$ of $X$ and $\HR_0[(\mathfrak L_X,d_v),\g]\,\cong\,\Q$. Thus, the right-hand side of the conjectured isomorphism is $\HR_\ast(\mfa_X,\g)$. Similarly, $\HR_0(X,G)=\Q$, which implies that $\widehat{\HR}_{\ast}(X,G) \,\cong\, {\HR}_{\ast}(X,G)$. Thus, Conjecture~\ref{conj5} is equivalent to Theorem~\ref{conj1} for simply-connected spaces.

\section{Spaces with polynomial representation homology and the Strong  Macdonald Conjecture}
\la{S4}
In this section we address Question \ref{qst2} and prove our second main result --- Theorem \ref{drinhomcprxs} --- stated as Theorem~\ref{drinhomcpr} in the Introduction. We will also work out a number of explicit examples illustrating this theorem
and linking it to the Strong Macdonald Conjecture.

\subsection{Lie-Hodge decompositions} 
\la{secdrintr}
Given a Lie algebra $ \mfa \in \LAlg_k $, we consider
the symmetric ad-invariant multilinear forms on $ \mfa \,$ of a (fixed) degree $ d \ge 1 $. Every such form is induced from the universal one: $\,\mfa \times \mfa \times \ldots \times \mfa \to \lambda^{(d)}(\mfa) \,$, which takes its values in  $\,\lambda^{(d)}(\mfa) := \sym^d(\mfa)/[\mfa, \sym^d(\mfa)]\,$, the space of coinvariants of the adjoint representation of $ \mfa $ in the $d$-th symmetric power of $\mfa$.
The assignment $\,\mfa \mapsto \lambda^{(d)}(\mfa)\,$ defines a (non-additive) functor that naturally extends to the category of DG Lie algebras:
\begin{equation}
\la{lam}
\lambda^{(d)}:\,\DGL_k \rar \Com_k \ .
\end{equation}
The functor \eqref{lam} is not homotopy invariant (it does not preserve quasi-isomorphisms); however, as shown in \cite{BFPRW}, it has a well-defined left derived functor
\begin{equation}
\la{Llam}
\L\lambda^{(d)}:\,\Ho(\DGL_k) \to \Ho(\Com_k)\ ,
\end{equation}
We write  $\,\HC^{(d)}_{\ast}(\mfa)\,$ for the homology of $\, \L\lambda^{(d)}(\mfa) \,$ and call it the 
{\it Lie-Hodge homology}\footnote{Observe that $ \lambda^{(1)} $ is just the abelianization functor on
Lie algebras; hence, for $d=1$, the existence of \eqref{Llam}  follows from general results of \cite{Q1}, and 
$\, \HC^{(1)}_{\ast}(\mfa) \,$ coincides with the Quillen homology of $ \mfa $, which is known to be isomorphic
(up to shift in degree) to the classical Chevalley-Eilenberg homology of $ \mfa $. For $d = 2$, the functor $ \lambda^{(2)} $ was introduced by Drinfeld in \cite{Dr}; the existence of $ \L\lambda^{(2)} $ was established by Getzler and Kapranov \cite{GK} who suggested to view $ \HC^{(2)}_{\ast}(\mfa) $ as an analogue of cyclic homology for Lie algebras. For arbitrary $ d \ge 1 $,  the existence of  \eqref{Llam} was proven in \cite[Sect. 7]{BFPRW} using some earlier general results of \cite{BKR}.} of $ \mfa $.

Next, we consider the (reduced) cyclic functor on  associative DG algebras
$$ 
(\mbox{--})_{\n}\,:\,\DGA_{k/k} \rar \Com_k\,\qquad R \mapsto R/(k+[R,R])\,,
$$
Observe that each  functor $ \lambda^{(d)} $ comes together with a natural transformation $\,\lambda^{(d)} \to U_\n \,$ induced by the symmetrization maps
\begin{equation*}
\la{symfun}
\sym^d(\mfa) \to U\mfa\ ,\quad x_1 x_2 \ldots x_d\, \mapsto\,
\frac{1}{d!}\,\sum_{\sigma \in {\Sigma}_d}\, \pm \,x_{\sigma(1)} \cdot x_{\sigma(2)} \cdot
\ldots \cdot x_{\sigma(d)}\ ,
\end{equation*}
where $ U\mfa $ is the universal enveloping algebra of $\mfa $, and by the Poincar\'e-Birkhoff-Witt Theorem, these natural transformations assemble to an isomorphism of functors
\begin{equation}
\la{eqv1}
\bigoplus_{d=1}^{\infty} \lambda^{(d)} \,\cong \, U_\n \ .
\end{equation}
On the other hand, by a well-known theorem of Feigin and Tsygan \cite{FT}, the functor $ (\,\mbox{--}\,)_\n $ has a left derived functor $ \L(\,\mbox{--}\,)_\n:\, \Ho(\DGA_{k/k}) \to \Ho(\Com_k) $ that computes
the reduced cyclic homology $\,\rHC_\ast(R)\,$ of an associative algebra
$ R \in \DGA_{k/k} $. Since $  U $ preserves quasi-isomorphisms and maps cofibrant
DG Lie algebras to cofibrant DG associative algebras, the isomorphism \eqref{eqv1}
induces an isomorphism of derived functors:
\begin{equation}
\la{eqv2}
\bigoplus_{d=1}^{\infty}\, \L\lambda^{(d)}\, \cong\, \L(\,\mbox{--}\,)_\n \circ \,U \ ,
\end{equation}
At the level of homology, \eqref{eqv2} yields a direct sum decomposition
\begin{equation}
\la{hodgeds}
\rHC_{\ast}(U\mfa) \,\cong\,\bigoplus_{d=1}^{\infty}\, \HC^{(d)}_{\ast}(\mfa)\ \text{.}
\end{equation}
which we call the {\it Lie-Hodge decomposition} for $ U\mfa $ ({\it cf.}~\cite[Theorem 7.4]{BFPRW}).

Now, let $X$ be a simply connected topological space, and let $\LL X $ denote the free loop space over $X$, i.e. the space of all continuous maps $ S^1 \to X $ equipped with compact open topology. This space carries a natural
$ S^1$-action (induced by rotations of $ S^1 $), hence we  can consider its equivariant homology
$$
\H^{S^1}_{\ast}(\LL X, k) := \H_\ast(ES^1 \times_{S^1} \LL X, k)\ .
$$
We will actually work with a  reduced version of $S^1$-equivariant homology of $\LL X $ defined by
$$
\rH^{S^1}_{\ast}(\LL X, k) := \Ker[\,\H^{S^1}_{\ast}(\LL X, k) \xrightarrow{\pi_*} \H_\ast(BS^1, k)\,]\ ,
$$
where the map $ \pi_* $ comes from the natural (Borel) fibration
\begin{equation}
\la{fibr}
\LL X \to ES^1 \times_{S^1} \LL X \xrightarrow{\pi} BS^1 \ .
\end{equation}
The following theorem is a well-known result due to Goodwillie \cite{Go} and Jones \cite{J}.
\bthm[\cite{Go, J}]
\la{top1}
Assume that $X$ is a simply connected space of finite rational type, and let $ \mfa_X $ be
a Quillen model of $ X $. Then there is a natural isomorphism of graded vector spaces
\begin{equation}
\la{ax}
\rHC_{\ast}(U\mathfrak{a}_X) \xrightarrow{\sim} {\rH}^{S^1}_{\ast}(\LL X, \Q)
\end{equation}
\ethm
Now, for each integer $n \ge 0 $, consider the $n$-fold covering of the circle:
$$
\omega^n\,:\,S^1 \rar S^1\ ,\quad e^{i\theta} \mapsto e^{in\theta}\ ,
$$
and denote by $\,\varphi_X^n\,;\,\mathcal L X \rar \mathcal L X\,$ the induced map on $ \LL X$. While for $n\ge 1 $, the maps $\varphi_X^n$ are not equivariant with respect to the $S^1$-action on $ \LL X$, they give naturally a commutative diagram in the homotopy category
\begin{equation}\la{diag}
\begin{diagram}[small, tight]
\LL X & \rTo^{\varphi_X^n} & \LL X\\
 \dTo & & \dTo\\
(ES^1 \times_{S^1} \LL X)_\Q & \rTo^{\tilde{\varphi}_X^n} &  (ES^1 \times_{S^1} \LL X)_\Q \\
  \dTo^{} & & \dTo_{}\\
  (BS^1)_\Q & \rTo^{B\omega^n} & (BS^1)_\Q\\
  \end{diagram}\ ,
\end{equation}
where columns are obtained by taking the rationalization of the Borel fibration \eqref{fibr} (see \cite{BFG}). The maps $\tilde{\varphi}_X^n $ in \eqref{diag}  induce graded endomorphisms
\begin{equation*}
\label{powerop}
\tilde{\Phi}_X^n :\,  {\rH}^{S^1}_{\ast}(\LL X, k)  \to   {\rH}^{S^1}_{\ast}(\LL X, k)\ , \quad n \geq 0\ ,
\end{equation*}
defined over any field $ k $ containing $\Q$. We call these endomorphisms the  {\it power} or {\it  Frobenius operations}  on $ {\rH}^{S^1}_{\ast}(\LL X, k) $ and write $ {\rH}^{S^1,\, (p)}_{\ast}(\LL X, k) $  for their (common) eigenspaces with eigenvalues $ n^p $: i.e.,
\begin{equation}\la{phi}
{\rH}^{S^1,\, (p)}_{\ast}(\LL X, k)\, :=
 \bigcap_{n \ge 0}\,\Ker(\tilde{\Phi}_X^n - n^p\,\id)\ .
\end{equation}

The next result proven in \cite{BRZ} provides a topological interpretation of the Lie-Hodge homology.
\bthm[\cite{BRZ}, Theorem~4.2]
\la{top2}
The Goodwillie-Jones isomorphism \eqref{ax} restricts to isomorphisms 
\begin{equation*}
\HC^{(d)}_{\ast}(\mathfrak{a}_X)  \cong  {\rH}^{S^1, \,(d-1)}_{\ast}(\LL X, \Q)\ ,\quad \forall\, d\ge 1 \ .
\end{equation*}
\ethm
It follows from Theorem~\ref{top1} and Theorem~\ref{top2} that, for a Quillen model $\mfa_X$ of a simply connected space $X$,
the Lie-Hodge decomposition \eqref{hodgeds} coincides with the topological Hodge decomposition
$$
\rH^{S^1}_\ast(\LL X,k) \,=\,\bigoplus_{p=0}^{\infty}\,\rH^{S^1,\, (p)}_\ast(\LL X,k)
$$

\subsection{The Drinfeld homomorphism}
\label{Dtracemaps}
Our next goal is to describe certain natural trace maps with values in representation homology. These maps were originally constructed in \cite{BFPRW, BFPRW2} as (derived) characters of finite-dimensional Lie representations. We will give a topological interpretation of these characteris in terms of free loop spaces. From now on, we assume that $ G $ is a {\it reductive} affine algebraic group over $k$. We denote by  $ I(\g)\,:=\,\sym(\g^{\ast})^{G}\,$ the space of invariant polynomials on the Lie algebra $ \g $ of $G$, and for $d \geq 0$, write $  I^d(\g) \subset  I(\g) $ for the subspace of homogeneous polynomials of degree $ d $.

For any commutative algebra $ A $, there is a natural
symmetric invariant $d$-linear form $\,\mfa(A) \times \mfa(A)
\times \ldots \times \mfa(A) \to \lambda^{(d)}(\mfa) \otimes A \,$ on the
current Lie algebra $ \mfa(A) $. Hence, by the universal property of $ \lambda^{(d)} $,
we have a canonical map
\begin{equation}
\la{canB}
\lambda^{(d)}[\mathfrak{a}(A)] \rar \lambda^{(d)}(\mathfrak{a}) \otimes A \,\text{.}
\end{equation}
Applying $ \lambda^{(d)} $ to the universal representation \eqref{univrep} and composing  with \eqref{canB}, we define
\begin{equation} \la{maplambdad}
\begin{diagram}
\lambda^{(d)}(\mathfrak{a}) & \rTo
& \lambda^{(d)}[\g(\mathfrak{a}_{\g})]
& \rTo & \lambda^{(d)}(\g) \otimes \mathfrak{a}_{\g} \end{diagram}\,\text{.}
\end{equation}
On the other hand, for the Lie algebra $ \g $, we have a canonical (nondegenerate) pairing
\begin{equation}
\la{pair}
I^d(\g) \times \lambda^{(d)}(\g) \to k
\end{equation}
induced by the linear pairing between $ \g^* $ and $\g$.
Replacing the Lie algebra $ \mfa $ in \eqref{maplambdad} by its cofibrant resolution
$ \mathcal L \sonto \mfa $ and using \eqref{pair}, we define the morphism of complexes
\begin{equation}
\la{dtrm}
\begin{diagram}
I^{d}(\g) \otimes \lambda^{(d)}(\mathcal L) & \rTo^{\eqref{maplambdad}} &
I^{d}(\g) \otimes \lambda^{(d)}(\g) \otimes \mathcal L_{\g}
& \rTo^{\eqref{pair}} & \mathcal L_{\g} \end{diagram} \ .
\end{equation}
For a fixed polynomial $ P \in I^{d}(\g)$, this morphism induces a map
on homology
$\,
\Tr^P_{\g}(\mathfrak{a})\,:\,
\HC_{\ast}^{(d)}(\mathfrak{a}) \rar \HR_{\ast}(\mathfrak{a},\g)\,
$,
which we call the {\it Drinfeld trace} associated to $P$. It is easy to check that
the image of \eqref{dtrm} is contained in the invariant subalgebra
$ \mathcal L_{\g}^{G} $ of $ {\mathcal L}_{\g}$, hence the Drinfeld trace is actually
a map
\begin{equation} \la{DTr}
\Tr^P_{\g}(\mathfrak{a}):\, \HC_{\ast}^{(d)}(\mathfrak{a}) \rar \HR_{\ast}(\mathfrak{a},\g)^{G} \,\text{.}
\end{equation}

Now, assume that $k=\c$ and $G$ is a complex reductive group of rank $l$. 
In this case, the algebra $ I(\g) = \sym(\g^*)^G $  
is freely generated  by a set of homogeneous polynomials $ \{P_1,\ldots,P_l\} $ whose degrees
$ d_i := \deg(P_i) $ are called the {\it fundamental degrees}   of $\g$. Fixing such a  set $ \{P_1,\ldots,P_l\} $ of generators in $I(\g) $, we assemble the associated trace maps \eqref{DTr} into a single homomorphism of graded commutative algebras
\begin{equation}\la{drinhom}
\Sym \Tr_\g(\mfa):\, \Sym_k\big[\Moplus_{i=1}^l \HC^{(d_i)}_\ast(\mfa)\big]  \rar \HR_\ast(\mfa,\g)^G\ .
\end{equation}
Following \cite{BFPRW,BFPRW2}, we call \eqref{drinhom} the {\it Drinfeld homomorphism} for $(\mfa,\g)$. We note that the Drinfeld homomorphism \eqref{drinhom} depends on the choice of polynomials $\{P_1,\ldots,P_l\} \subset I(\g)$, but for simplicity we suppress this in our notation. 

If $ \mfa = \mfa_X $ is a Lie model of a simply connected space $X$, by Theorem~\ref{top2}, $\,\HC^{(d)}_\ast(\mfa) \cong \rH^{S^1,(d-1)}_\ast(\LL X, \c)$.
On the other hand, by Theorem \ref{conj1}, $\HR_\ast(\mfa,\g)\,\cong\,\HR_\ast(X,G)$. Hence, the Drinfeld homomorphism for $X$ may be rewritten in the following topological form:
\begin{equation} 
\la{drinhomtop} 
\Sym_k \left[\bigoplus_{i=1}^l \rH^{S^1,(m_i)}_\ast(\LL X, k) \right] \rar \HR_\ast(X,G)^G\,,
\end{equation}
where the $m_i=d_i-1$ are the exponents of the Lie algebra of the group $G$. 

Our next goal is to express the Drinfeld homomorphism \eqref{drinhomtop}  in terms of the (minimal) Sullivan model of $\mathcal A_X$ of $X$. Recall that $\mathcal A:=\mathcal A_X$ is Koszul dual to the Lie algebra $\mfa:=\mfa_X$ in the sense that $\mathcal{A}\,\cong\,\C^{\ast}(\mfa;k)$. Then, by \cite[Prop. 7.8]{BFPRW}, there is an isomorphism of graded vector spaces for any $m \geq 0$,
\begin{equation} \la{cychomdual}  \HC^{(m+1)}_\ast(\mfa)\,\cong\, (\rHC^{(m)}_\ast(\mathcal{A}))^{\ast}[-1]\,,\end{equation}
where the superscript $(\mbox{--})^{\ast}$ stands for the graded $k$-linear dual. In particular, we have ({\it cf.} \cite[Thm. B]{BFG}) 
\begin{equation} \la{bfg} \rH^{S^1,(m)}_\ast(\LL X)\,\cong\, (\rHC^{(m)}_\ast(\mathcal{A}))^{\ast}[-1]\ .\end{equation}
On the other hand, by \cite[Thm. 6.7(b)]{BFPRW},
\begin{equation} \la{rephomdual} 
\HR_\ast(\mfa,\g)^G\,\cong\,\H^{-\ast}(\g(\mathcal{A}),\g;\c)\ .
\end{equation}
Now, for each $m \geqslant 0$, for each $P\,\in\,I^{m+1}(\g)$, define a linear map
\begin{equation} \la{phip0} \Phi_P\,:\, \C_\ast(\g(\bar{\mathcal A});\c) \rar  \Omega^m(\mathcal A)/d\Omega^{m-1}(\mathcal{A})[m+1]\,,\end{equation}
by the following explicit formula
\begin{equation} \la{phip}
\Phi_P\big((\xi_0 \otimes a_0) \wedge \ldots \wedge (\xi_m \otimes a_m)\big)\,=\, \frac{1}{(m+1)!} \sum_{\sigma \in \Sigma_{m+1}} \pm a_{\sigma(0)}da_{\sigma(1)}\ldots da_{\sigma(m)}P(\xi_{\sigma(0)},\ldots,\xi_{\sigma(m)})\,,\end{equation}
where $a_0,\ldots,a_m\,\in\,\bar{\mathcal{A}}$ and $\xi_0,\ldots,\xi_m\,\in\,\g$, and let $\Psi_P$ denote the composition
$$\Psi_P\,:\,\begin{diagram} \C_\ast(\g(\mathcal{A}),\g;\c) &\rInto & \C_\ast(\g(\bar{\mathcal{A}});\c) & \rTo^{\Phi_P} & \Omega^m(\mathcal A)/d\Omega^{m-1}(\mathcal A)[m+1] \end{diagram}\ . $$
\blemma \la{drintrduals}
$(i)$ The map $\Psi_P$ is a well-defined chain map whose graded linear dual induces on cohomology
$$\Psi_{P}^{\ast}\,:\, \big(\rHC^{(m)}_\ast(\mathcal A)\big)^{\ast}[-1] \rar \H^\ast(\g(\mathcal A),\g;\c) \ . $$
$(ii)$ The following diagram commutes:
$$ 
\begin{diagram}
\HC^{(m+1)}(\mfa) & \rTo^{\Tr^P_\g} & \HR_\ast(\mfa,\g)^G\\
\dTo^{\eqref{cychomdual}}_{\cong} & & \dTo^{\eqref{rephomdual}}_{\cong}\\
\big(\rHC^{(m)}_\ast(\mathcal A)\big)^{\ast}[-1]  & \rTo^{\Psi_P^{\ast}} & \H^{-\ast}(\g(\mathcal A), \g;\c)\\
\end{diagram}
$$
\elemma 
\bproof
We first recall from \cite{BFPRW2} a construction of the Drinfeld traces via the Chern-Simons formalism. Let $\mathrm{DR}(\mathcal A):=\Sym_{\mathcal A}(\Omega^1 \mathcal{A}[-1])$ equipped with the differential $d+\partial$, where $d$ is the de Rham differential and $\partial$ is the internal differential induced by the differential on $\mathcal A$. Let $\mathrm{DR}^{\geqslant n}(\mathcal A)$ denote the two sided DG ideal in $\mathrm{DR}(\mathcal A)$ generated by $\Omega^n \mathcal{A}[-n]$, and let $\tau^n \mathrm{DR}(\mathcal A)$ denote the quotient $\mathrm{DR}(\mathcal A)/\mathrm{DR}^{\geqslant (n+1)}(\mathcal A)$. Note that since $\mathcal{A}$ is augmented, so is $\tau^n\mathrm{DR}(\mathcal A)$ for each $n$. Let $\tau^n\overline{\mathrm{DR}}(\mathcal A)$ denote the corresponding augmentation ideal. Since $\mathcal A$ is smooth as a graded commutative algebra, \cite[Thm. 5.4.7]{L}, there is a canonical isomorphism for each $m \geq 0$
$$ \rHC^{(m)}_\ast(\mathcal A)\,\cong\, \H_\ast\left(\tau^{m}\overline{\mathrm{DR}}(\mathcal A)[2m]\right) .  $$
Further, since $\mathcal A$ is a graded symmetric algebra equipped with an extra differential, the canonical projection 
\begin{equation} \la{projtoomegam} \tau^m\overline{\mathrm{DR}}(\mathcal A) \twoheadrightarrow \left(\Omega^m(\mathcal A)/d\Omega^{m-1}(\mathcal{A})\right)[-m] \end{equation}
is a quasi-isomorphism (see \cite[Thm. 5.4.12]{L}). Hence, 
$$
\rHC^{(m)}_\ast(\mathcal A)\,\cong\, \H_\ast \left[ \left(\Omega^m(\mathcal A)/d\Omega^{m-1}(\mathcal{A})\right)[m] \right]\ . $$
 Next, note that the Chevalley-Eilenberg chain complex $\C_\ast(\g(\bar{\mathcal A});\c)$ is a cocommutative DG coalgebra. Hence, the Hom complex $\mathcal{B}:=\Hom(\C_\ast(\g(\bar{\mathcal A});k), \tau^m\overline{\mathrm{DR}}(\mathcal A))$ has the structure of a commutative DG algebra with convolution product.  There is a $\g$-valued one form $\theta\,\in\,\mathcal{B}^1 \otimes \g$ on $\mathcal B$ such that the restriction of $\theta$ to $\wedge^k(\g(\bar{\mathcal A}))$ vanishes for $k \neq 1$ and $\theta|_{\g(\bar{\mathcal A})}$ coincides with the composite map 
$$ \begin{diagram} \g \otimes \bar{\mathcal A} & \rTo & \bar{\mathcal A} \otimes \g & \rInto & \tau^m\overline{\mathrm{DR}}(\mathcal A) \end{diagram}\,, $$
where the first arrow is the obvious swap map. For $P\,\in\,I^{m+1}(\g)$, the Chern-Simons form $TP(\theta)\,\in\,\mathcal{B}^{2m+1}$ satisfies $\delta(TP(\theta))\,=\,P(\Omega^{m+1})$, where $\Omega\,\in\,\mathcal{B}^2 \otimes \g$ is the curvature of $\theta$. Since $\Omega^{m+1}=0$ by \cite[Prop. A.2]{BFPRW2}, $TP(\theta)\,\in\,\mathcal{B}^{2m+1}$ is a cocycle. It follows that 
$s^{2m}TP(\theta)$ defines a map of complexes
\begin{equation} \la{chernsimons} \frac{1}{(m+1)!}s^{m}TP(\theta)\,:\,\C_\ast(\g(\bar{\mathcal A});\c) \rar \tau^m\overline{\mathrm{DR}}(\mathcal A)[1]\ .\end{equation}
An explicit formula for the map \eqref{chernsimons} has been given in \cite{Fe} (also see \cite[Eqn. 2.2]{Te} and \cite[Prop. A.3]{BFPRW2}). By \cite[Prop. A.5]{BFPRW2}, the composition of the canonical projection \eqref{projtoomegam} with \eqref{chernsimons} coincides with $\Phi_P$ (see \eqref{phip0}). This implies $(i)$. $(ii)$ is then an immediate consequence of the main result of \cite{BFPRW2} (see {\it loc. cit.,} Theorem 3.2).
\eproof

As a consequence of Lemma \ref{drintrduals}, we obtain the following description of the Drinfeld homomorphism in terms of the (minimal) Sullivan model.
\bthm \la{TdrinhomS}
 For a simply connected space $X$ with minimal Sullivan model $\mathcal{A}_X$, the Drinfeld homomorphism \eqref{drinhomtop} is given by the map
\begin{equation} \la{drinhomsullivan}  \Psi^{\ast}(\mathcal A_X)\,:\,\Sym \left(\bigoplus_{i=1}^l\big(\rHC^{(m_i)}_\ast(\mathcal A_X)\big)^{\ast}[-1] \right) \rar \H^{-\ast}(\g(\mathcal A_X),\g;\c) \end{equation}
obtained by assembling the maps $\Psi_{P_i}^{\ast}$ for a set  $\{P_1,\ldots,P_l\}$ of homogeneous generators of $I(\g)$.
\ethm
\subsection{Spaces with polynomial representation homology}
We now address Question \ref{qst2} stated in the Introduction. Recall that this question is asking for a characterization of spaces $X$ and groups $G$ for which the algebra $\HR_\ast(X,G)^G $ is free of locally finite type over $k$. At the moment, a complete characterization of such pairs $(X,G)$ seems to be out of reach. In what follows, we will consider two --- in some sense extreme --- cases: we first describe a class of algebraic groups $G$ such that $\HR_\ast(X,G)^G $ is free for {\it all}
spaces $X$ (see Theorem~\ref{gagm}) and then characterize a class of spaces $X$ such that $\HR_\ast(X,G)^G $ is free for {\it all} complex reductive groups $G$ (see Theorem~\ref{drinhomcprxs}).

\bthm 
\la{gagm}
If $G$ is a commutative affine algebraic group of dimension $ l$, then, for any simply connected space $X$ of finite rational type, there is an algebra isomorphism  
$$ 
\HR_\ast(X,G)\,\cong\, \Sym_k \big[\H_{\ast+1}(X;k)^{\oplus  l}\big]\,.
$$
\ethm
\bproof
We first prove the desired result in the case when $\dim_k\,G=1$. Let $\mfa=(\mfa(V),\partial)$ be a minimal Quillen model of $X$ freely generated by a graded vector space $V$ with differential $\partial$. Then, $\mfa_\g=\Sym(V)$ with $0$ differential. On the other hand, $\HC^{(1)}_\ast(\mfa)\,\cong\,\mfa/[\mfa,\mfa]\,\cong\,\, V\,,$ with $0$ differential. It is easy to see that the Drinfeld trace\footnote{We remark that the construction of the Drinfeld trace \eqref{DTr} goes through even when $G$ is not reductive for $P\,\in\,\sym(\g^{\ast})^{\ad\,\g}$. Hence, when $G$ (and therefore, $\g$) is abelian, one has the Drinfeld trace $\Tr^{P}_\g(\mfa)\,:\,\HC^{(1)}_\ast(\mfa) \rar \HR_\ast(\mfa,\g)$ for every $P\,\in\,\g^{\ast} \hookrightarrow \sym(\g^{\ast})$. Fixing a basis of $\g^{\ast}$, we assemble the associated traces into the Drinfeld homomorphism as in \eqref{drinhom}.} corresponding to the generator of $\Sym(\g^{\ast})$ is the map 
$$ \mfa/[\mfa,\mfa]\,\cong\,V \hookrightarrow \Sym(V)\,\cong\,\mfa_\g\ .$$
The corresponding Drinfeld homomorphism is therefore identified with the identity on $\Sym(V)$. Finally, note that if $\g$ is abelian of dimension $l$, the Drinfeld homomorphism for $\g$ becomes the map
$$\Sym\big(\mfa/[\mfa,\mfa]\big)^{\otimes l} \stackrel{\tau^{\otimes l}}{\longrightarrow} \Sym(V)^{\otimes l}\,,  $$
where $\tau$ is the Drinfeld homomorphism for a one-dimensional Lie algebra. Hence, it is an isomorphism. Since 
$\H_\ast[\mfa/[\mfa,\mfa]]\,\cong\,\H_{\ast+1}(X;k)$, and since $G$ is abelian, the desired formula for $\HR_\ast(X,G)$ follows as well.
\eproof 
\begin{remark}
It is well known (see, e.g., \cite[Chap.~IV]{DG}) that over an algebraically closed  field of characteristic $0$, any finite-dimensional commutative affine algebraic group is isomorphic to the product of an algebraic torus and a vector group over $k$: i.e,  $ G \cong \mathbb{G}_m^r \times \mathbb{G}_a^s$. If $ G\,\cong\,\mathbb{G}_a^s$, then the result of Theorem \ref{gagm} actually holds for an arbitrary --- not necessarily simply connected ---  space $X$ (see \cite[Example 3.1]{BRYI}). 

\end{remark}

\vspace*{2ex}

 The next theorem (stated as Theorem~\ref{drinhomcpr} in the Introduction) provides a (partial) answer to Question \ref{qst3}, characterizing in simple cohomological terms spaces for which the Drinfeld homomorphism is an isomorphism for all reductive groups $G$. As explained in the introduction, the proof of this theorem relies on a theorem of Fishel, Grojnowski and Teleman \cite{FGT} (formerly known as the Strong Macdonald Conjecture).
\bthm \la{drinhomcprxs}
Let $X$ be a simply connected space such that its rational cohomology algebra  $\H^\ast(X;\Q)$  is either generated by one element $($in any dimension$)$ or {\rm freely} generated by two elements: one in even and one in odd dimensions.
Then, the Drinfeld homomorphism \eqref{drinhomtop} is an isomorphism for every complex reductive algebraic group $G$.
\ethm

The proof of Theorem \ref{drinhomcprxs} is based on the following  refinement of Theorem B of \cite{FGT}.

\bprop \la{polysin2var}
Let $\mathcal{A}\,\cong\,\c[z,s]$ with $0$ differential, where $|z| \ge 2 $ is even and $|s| \ge 3 $ is odd. Then, the map $\Psi^{\ast}(\mathcal A)$ (see \eqref{drinhomsullivan}) is an isomorphism.
\eprop
\bproof
 Viewing all (DG) algebras as homologically graded by inverting degrees, we note that 
 $$
 \Sym_{\mathcal{A}} \Omega^1(\mathcal{A})[1]\,=\,\c[z,s,dz,ds]\ , 
 $$
 where $\deg\,dz\,=\,1-d$ and $\deg\,ds\,=\,1-l$. Here, $d\,:=\,|z|$ and $l\,:=\,|s|$ denote the {\it cohomological} degrees of $z$ and $s$ respectively, whence $\deg\,z\,=\,-d$ and $\deg\,s\,=\,-l$. Hence, for $m \geq 1$,
$$ \Omega^m(\mathcal A)\,=\, \c[z]dz(ds)^{m-1} \oplus \c[z](ds)^m \oplus \c[z]dz s(ds)^{m-1} \oplus \c[z]s(ds)^m\,,$$
and it is easy to verify that for $f(z)\,\in\,\c[z]$,
\begin{equation} \la{eqcycles} f(z)(ds)^m\,\equiv\, -f'(z)dz s(ds)^{m-1}\,,\qquad f(z)dz(ds)^{m-1} \,\equiv\, 0\end{equation}
modulo $\Omega^{m-1}(\mathcal A)$. For $m=0$, 
$$ \overline{\mathrm{DR}}^0(\mathcal A)\,=\,\bar{\mathcal{A}}\,=\,  z\c[z] \oplus \c[z]s\ .$$
Since the differential on $\mathcal A$ is trivial, there are isomorphisms of graded vector spaces for $m \geq 0$ 
\begin{equation} \la{cycminsullivan}  \rHC^{(m)}_\ast(\mathcal{A}) \,\cong\,\Omega^m(\mathcal A)/d\Omega^{m-1}(\mathcal A)[m] \,\cong\, \c[z]s(ds)^m \oplus \c[z]dz\cdot s(ds)^{m-1} \,,\end{equation}
where for $m=0$, the formal summand $\c[z]dzs(ds)^{-1}$ of $\rHC_0(\mathcal A)$ is identified with the summand $z\c[z]$ of $\bar{\mathcal A}$ by the isomorphism 
$f(z) \mapsto df(z)=f'(z)dz$. The restriction of the inverse of the isomorphism \eqref{cycminsullivan} to each summand is given by the obvious inclusion into $\Omega^m (\mathcal A)[m]$ followed by the canonical projection. Composing the isomorphism \eqref{cycminsullivan} with projection to each factor on the right hand side yields two linear maps 
\begin{equation} \la{se} S\,:\,\rHC^{(m)}_\ast(\mathcal{A}) \rar \c[z]\,,\qquad E\,:\,\rHC^{(m)}_\ast(\mathcal{A}) \rar \c[z]dz\ .\end{equation}
As in \cite[Sec. 1.8]{FGT}, there is an isomorphism of DG coalgebras
$$ \C_\ast(\g[z],\g;\Sym^c(s\g[z][1]))\,\cong\,\C_\ast(\g(\mathcal A),\g;\c)\,,$$ $$ \bigwedge_{i=1}^p \xi_i(f_i) \otimes \bigwedge_{j=1}^q \xi_{p+j}(sf_{p+j})  \mapsto \bigwedge_{i=1}^p \xi_i(f_i) \wedge \bigwedge_{j=1}^q \xi_{p+j}(sf_{p+j})\,, $$
where $\xi_i\,\in\,g$, $f_i\,\in\,\c[z]$ for all $i$ and for $f\,\in\,\mathcal A$, $\xi(f):=\xi \otimes f$ for $\xi \,\in\,\g$. 
Identifying $\C_\ast(\g(A),\g;\c)$ with $\C_\ast(\g[z],\g;\Sym^c(s\g[z][1]))$ via the above isomorphism, we note that for $P\,\in\,I^{m+1}(\g)$, the restriction of $\Psi_P$ to $\Lambda^r(\g[z]/\g) \otimes \Sym^q(s\g[z][1])$ vanishes for $r \geq 2$. Indeed, for $r \geq 3$ this vanishing is obvious since every summand contributing to the right hand side of \eqref{phip} has two factors of the form $f'(z)dz$. For $r=2$ the only summands on the right hand side of \eqref{phip} not having two factors of the form $f'(z)dz$ are of the form $f(z)dz(ds)^{m-1}$, which lies in $d\Omega^{m-1}(\mathcal A)$. 

Now, note that for $f_0,\ldots,f_m\,\in\,\c[z]$, 
\begin{equation} \la{sp1} \Psi_P(\bigwedge_{i=0}^m \xi_i(sf_i))\,=\,\frac{1}{(m+1)!} \sum_{\sigma\,\in\,\Sigma_{m+1}}\pm sf_{\sigma(0)}d(sf_{\sigma(1)}) \ldots d(sf_{\sigma(m)})P(\xi_{\sigma(0)},\ldots,\xi_{\sigma(m)}) \end{equation}
Since $d(sf_i)=(ds)f_i-sf'_i(z)dz$ and since $s^2=0$,  the right hand side of \eqref{sp1} equals
\begin{equation} \la{sp2} \frac{1}{(m+1)!} \sum_{\sigma\,\in\,\Sigma_{m+1}} \pm \prod_{i=0}^m f_{\sigma(i)} \cdot s(ds)^mP(\xi_{\sigma(0)},\ldots,\xi_{\sigma(m)})\,=\, P(\xi_0(f_0),\ldots,\xi_m(f_m))s(ds)^m\ .\end{equation}

Next, note that 
\begin{equation} \la{ep1} 
\Psi_P\big(\xi_0(f_0) \otimes \bigwedge_{i=1}^m \xi_i(sf_i)\big) \,=\, \frac{1}{(m+1)!}\sum_{\substack{\sigma\,\in\,\Sigma_{2m+1}\\\sigma(0)=0}} \pm  f_0 \prod_{i=1}^m d(sf_{\sigma(i)})P(\xi_0,\ldots,\xi_m) \end{equation}
$$ 
+ \frac{1}{(m+1)!}
\sum_{\substack{\sigma\,\in\,\Sigma_{m+1}\\\sigma(0) \neq 0}} \pm  sf_{\sigma(0)}\prod_{\substack{i=1\\\sigma(i) \neq 0}}^n d(sf_{\sigma(i)})f'_0(z)dzP(\xi_{\sigma(0)},\ldots,\xi_{\sigma(m)})
\ .$$
Since $d(sf_i)=(ds)f_i-sf'_i(z)dz$ and since $s^2=0$,  the second summand on the right hand side of \eqref{ep1} equals
$$-\frac{m}{(m+1)} f'_0(z)\prod_{i=1}^m f_i(z)dzs(ds)^{m-1}\ . $$
On the other hand, the first summand coincides with 
$$ \frac{1}{(m+1)} \left(f_0(z)\big(\prod_{i=1}^m f_i\big)'(z)dzs(ds)^{m-1}+ \prod_{i=0}^m f_i(z)(ds)^m\right)\,=\, -\frac{1}{(m+1)}f'_0(z)\prod_{i=1}^mf_i(z)dzs(ds)^{m-1} \ .$$
The last equality above is by \eqref{eqcycles}. Hence, 
\begin{equation} \la{ep2} \Psi_P\big(\xi_0(f_0) \otimes \bigwedge_{i=1}^m \xi_i(sf_i)\big)\,=\,-\left[f_0(z)\prod_{i=1}^m f_i(z) dzs(ds)^{m-1}\right]\,, \end{equation}
where $[\mbox{--}]$ stands for the class in $\rHC^{(m)}_\ast(\mathcal A)[1]$. It follows from \eqref{sp1}, \eqref{sp2} and \eqref{ep2} that $S \circ \Psi_P$ (resp., $E \circ \Psi_P$), viewed as a map of complexes $\C_\ast(\g[z],\g,\Sym^c(s\g[z][1])) \rar \Omega^m(\mathcal A)/d\Omega^{m-1}(\mathcal A)[m+1]$ coincides with the map 
$S_P$ (resp., $-E_P$) defined in \cite[Thm. B]{FGT} as a map of $\Z_2$-graded vector spaces (though they differ as maps of $\Z$-graded vector spaces). It follows from {\it loc. cit.} that the map $\Psi^{\ast}(\mathcal A)$ (see \eqref{drinhomsullivan}) and hence, the Drinfeld homomorphism \eqref{drinhom}, is an isomorphism of $\Z_2$-graded vector spaces (and therefore, of $\Z$-graded vector spaces) as desired.  Note that in our case, the restricted dual of $\C_\ast(\g(\mathcal{A}),\g;\c)$ in the sense of \cite{FGT} coincides with all of $\C^{\ast}(\g(\mathcal{A}),\g;\c)$ since the fact that $\bar{\mathcal A}$ is concentrated in cohomological degree $\geqslant 2$ ensures that $\C_\ast(\g(\mathcal{A}),\g;\c)$ is finite dimensional in each homological degree.      
\eproof
\bproof[Proof of Theorem \ref{drinhomcprxs}]
First, we consider the case when $\H^\ast(X;\Q)\,\cong\,\Q[z]$, where $z$ is a generator of even dimension $\geqslant 2$. By \cite[Prop. 5.1]{Me}, the (complexified) minimal Sullivan model of $X$ is given by $\mathcal{A}=\c[z]$ (with zero differential). Since the Drinfeld homomorphism is identified with the map $\Psi^\ast(\mathcal A)$ (see \eqref{drinhomsullivan}) by Theorem \ref{TdrinhomS}, the desired result follows in this case from the classical fact that $\Psi^{\ast}(\mathcal A)$ is an isomorphism for $\mathcal{A}=\c[z]$ (see \cite[Sec. 3]{Te}; also see \cite{Fe}). Next, we consider the case when $\H^\ast(X;\Q)\,\cong\,\Q[s]$, where $s$ is a generator of odd cohomological degree $r \geq 3$. Thus, $X$ has the rational homotopy type of an odd sphere. It follows that the Quillen model $\mfa$ of $X$ is a free Lie algebra on a single generator $u$ of (even) homological degree $r-1$. The Drinfeld homomorphism for $X$ becomes the map
$$\Sym \big(\bigoplus_{i=1}^l \c \cdot u^{d_i}\big) \rar \Sym \big(\g^{\ast}[r-1]\big)^G\,, \qquad u^{d_i} \mapsto P_i\ .  $$
That this is an isomorphism then amounts to the classical fact that $I(\g)$ is generated by the set of homogeneous polynomials $\{P_1,\ldots,P_l\}$. 

It therefore, remains to consider the possibilities that $\H^\ast(X;\Q)$ is a truncated polynomial algebra on a single generator of even dimension, or that $\H^\ast(X;\Q)$ is a polynomial algebra in two homogeneous generators, one of even dimension. In the latter case, by \cite[Prop. 5.1]{Me}, the (complexified) minimal Sullivan model of $X$ is $\mathcal{A}=\c[z,s]$ with zero differential; in the former case, the (complexified) minimal Sullivan model is given by  $\mathcal{A}_r=\c[z,s],\,\partial s=z^{r+1}$  where $z$ and $s$ are of {\it cohomological} degree $d$ and $d(r+1)-1$ respectively, where $d$ is even (see \cite[Sec.5.3]{Me}). If $\mathcal{A}$ (with zero differential) is the minimal Sullivan model of $X$, the desired result is immediate from Proposition \ref{polysin2var}, since the Drinfeld homomorphism is identified with the map $\Psi^{\ast}(\mathcal A)$ by Theorem \ref{TdrinhomS}. Now, assume that the minimal Sullivan model of $X$ is $\mathcal{A}_r$. Let $\mathrm{C}_r$ (resp., $\mathrm{C}$) denote the graded linear dual of $\mathcal{A}_r$ (resp., $\mathcal{A})$). Note that $\mfa_r:=\cb_{\mathtt{Comm}}(\mathrm{C}_r)$ is a Quillen model of $X$, where $\cb_{\mathtt{Comm}}\,:\,\cDGC_{k/k} \rar \DGL_k$ is the cobar functor (see \cite[Sec. 6.2.1]{BFPRW}). Explicitly, $\mfa_r=(\mfa(\bar{\mathrm{C}_r}[-1]),d_1+d_2)$, the graded free Lie algebra generated by $\bar{\mathrm{C}_r}[-1]$ with differential given by the sum of two derivations $d_1$ (induced by the differential on $\mathrm{C}_r)$ and $d_2$ (induced by the coproduct on $\mathrm{C}_r$). Equip $\mfa_r$ with an (increasing) filtration by `internal degree' by letting 
$$F_p\mathcal{L}_r \,=\,\sum_{s \geqslant 1} \sum_{d_1+\ldots+d_s \leqslant p+s} [(\bar{C_r}[-1])_{d_1},[(\bar{C_r}[-1])_{d_2},[\ldots,[(\bar{C_r}[-1])_{d_{s-1}},(\bar{C_r}[-1])_{d_s}]\ldots]]] \ . $$
Then $F_\ast\mfa_r$ is a bounded below exhaustive filtration on $\mfa_r$, and induces (bounded below, exhaustive) filtrations on $\lambda^{(p)}(\mfa_r)$ for all $p$ as well as on $(\mfa_r)_\g$ for any (reductive) $\g$. For a set $\{P_1,\ldots,P_l\}$ of homogeneous generators of $I(\g)$, the Drinfeld homomorphism \eqref{drinhom} is induced on homologies by the homomorphism of commutative DG algebras 
\begin{equation} \la{drhomfil} \Sym[\bigoplus_{i=1}^l \Tr_g^{P_i}]\,:\,\Sym\big[\bigoplus_{i=1}^l \lambda^{(d_i)}(\mfa_r)\big] \rar (\mfa_r)_\g^G\ .\end{equation}
The filtrations induced by $F_\ast$ make \eqref{drhomfil} a homomorphism of {\it filtered} commutative DG algebras, where the filtrations are bounded below and exhaustive. Let $\mfa :=\cb_{\mathtt{Comm}}(\mathrm{C})$. Since $\mathrm{gr}_{F_\ast}(\mfa_r)=\mfa$, the induced map on the $E^1$-page of the corresponding spectral sequences is the map induced on homologies by the DG algebra homomorphism
$$\Sym[\bigoplus_{i=1}^l \Tr_g^{P_i}]\,:\,\Sym\big[\bigoplus_{i=1}^l \lambda^{(d_i)}(\mfa)\big] \rar (\mfa)_\g^G\ . $$
By Theorem \ref{TdrinhomS}, the above map is identified with the map $\Psi^{\ast}(\mathcal A)$. Therefore, it is a quasi-isomorphism. The desired theorem is now immediate from Proposition \ref{polysin2var} and the classical convergence theorem (\cite[Thm. 5.5.1]{W}).
\eproof
\subsection{Examples}
\la{Sect4.4}
We will now illustrate Theorem \ref{drinhomcprxs} with explicit examples. Because of simplicity of cohomological conditions of Theorem \ref{drinhomcprxs}, the spaces satisfying these conditions are easy to construct (in fact, many of these spaces appear as basic examples in classical textbooks in algebraic topology, see e.g. \cite{Hat}). We divide them into three natural classes depending on the structure of their cohomology ring:

\vspace*{2ex}

\begin{enumerate}
\item[(I)]
\ $\H^\ast(X,\Q)\,\cong\,\Q[z]$, where $|z|$ is either odd or even.\\
\item[(II)]
\ $\H^\ast(X,\Q)\,\cong\, \Q[z,s]$, where $|z|$ is even and $|s|$ is odd.\\
\item[(III)]
\ $\H^\ast(X,\Q)\,\cong\,\Q[z]/(z^{r+1})$, where $|z|$ is even.
\end{enumerate}

\vspace*{2ex}

\noindent
{\it Throughout this section, as in Theorem \ref{drinhomcprxs}, $G$ stands for a complex reductive Lie group of rank $l \ge 1 \,$, $\,\g\,$ is the Lie algebra $ G$, and $\, \{m_1, m_2, \ldots, m_l\} \,$ are the classical exponents of $\g$.}

\subsubsection{Case I}
First, as already observed at the beginning of our proof of Theorem \ref{drinhomcprxs}, 
\begin{equation} \la{oddsphere} 
\HR_\ast({\bS}^{2r+1},G)^G\,\cong\, \Sym(\g^{\ast}[2r])^G\,\cong\,\c[P_1,\ldots,P_l]\,,\qquad \deg\,P_i\,=\,2r(m_i+1)\,,\quad 1 \leqslant i \leqslant l\ .\end{equation}
Here, $P_1,\ldots,P_l$ are the homogeneous generators of $I(\g)$ viewed as elements of $\Sym(\g^{\ast}[2r])^G$, whence $\deg\,P_i=2r(m_i+1)$. This computes the ($G$-invariant part of the) reprentation homology of $X$ for the case when $\H^{\ast}(X,\Q)\,\cong\,\Q[z]$, where $\deg\,z\,=\,2r+1$ (in which case $X$ is rationally equivalent to ${S}^{2r+1}$). 
Next, we have
\blemma \la{loopcpinfty}
If $\H^{\ast}(X,\Q)\,\cong\,\Q[z]$, where $z$ is of even dimension $d \geqslant 2$, then
$$\rH^{S^1,(0)}_\ast(\LL X)\,\cong\,\bigoplus_{j=1}^{\infty} \c \cdot \xi_j\ , \qquad \rH^{S^1,(i)}_\ast(\LL X)=0\ ,\quad i >0 \,,$$
where $\xi_j $ has homological degree $ dj-1\,$. 
\elemma 
\bproof 
 By \cite[Prop. 5.1]{Me}, the (complexified) minimal Sullivan model of $X$ is given by $\mathcal{A}=\c[z]$, where $z$ is of (even) cohomological degree $d$.  Hence,  $\rHC^{(0)}_\ast(\mathcal{A})\,\cong\,z\c[z]$ and $\rHC^{(i)}_\ast(\mathcal{A})=0$ for $i>0$. The desired lemma now follows from \cite[Thm. 4.2]{BRZ} and \cite[Prop. 7.8]{BFPRW}, which together imply
 $$ \rH_\ast^{S^1,(m)}(\LL X)\,\cong\,\big(\rHC^{(m)}_\ast(\mathcal{A})\big)^{\ast}[-1]\ .
 $$
 %
\eproof 
The following result is a consequence of Theorem \ref{drinhomcprxs} and Lemma \ref{loopcpinfty}
\bcor \la{rephomcase1}
 Let $X$ be a simply connected space with  $\H^\ast(X,\Q)\,\cong\,\Q[z]$, where $d\,:=\,|z|$ is even. Then, there is an isomorphism of graded commutative algebras
$$ \HR_\ast(X,G)^G\,\cong\, \Sym \left[\rH_{\ast+1}(X,\c)^{\oplus l_0} \right]\,,$$
where $l_0$ is the number vanishing exponent of $G$. More explicitly, 
$$\HR_\ast(X,G)^G\,\cong\, \Sym \left( \xi^{(i)}_j\,:\,1 \leqslant i \leqslant l_0\,,\,j \,\in\,\mathbb{N} \right)\ , $$
where the generators $\xi^{(i)}_j$ have homological degree $\,dj-1\,$ for all $\,i =1,2,\ldots, l_0 \,$. In particular, if $\,l_0=0\,$ (for example, if $G$ is complex semisimple), then $\,\HR_\ast(X,G)^G\,\cong\,\c\,$.
\ecor
The condition $\H^{\ast}(X,\Q)\,\cong\,\Q[z]$ holds, for example, for the following spaces (see \cite{Hat}):
\begin{itemize}
\item The spheres ${\bS}^{2n+1}
$,$\,n\geqslant 1$ ($|z|\,=\,2n+1$).
\item The Eilenberg-Maclane spaces $K(\Z,d)$, for even $d \geqslant 2$ ($|z|\,=\,d$).
\item $\c\mathbb{P}^{\infty}$ (rationally equivalent to $K(\Z,2)$).
\item $\mathbb{H}\mathbb{P}^{\infty}$ (rationally equivalent to $K(\Z,4)$).
\end{itemize}

\vspace{1ex}
\noindent
Hence, by Corollary \ref{rephomcase1}, we have 
\begin{eqnarray*}
\la{cpinfty} \HR_\ast(\c\mathbb{P}^{\infty},G)^G & \cong & \Sym \left( \xi^{(i)}_j\,:\,1 \leqslant i \leqslant l_0\,,\,j \,\in\,\mathbb{N} \right)\,,\quad \deg\,\xi^{(i)}_j\,=\,2j-1\ ,\\
\la{hpinfty} \HR_\ast(\mathbb{H}\mathbb{P}^{\infty},G)^G & \cong & \Sym \left( \xi^{(i)}_j\,:\,1 \leqslant i \leqslant l_0\,,\,j \,\in\,\mathbb{N} \right)\,,\quad \deg\,\xi^{(i)}_j\,=\,4j-1\ ,
\end{eqnarray*}
where $l_0$ is the number of vanishing exponents of $G$. 

\subsubsection{Case II} \la{s4.4.2} 
In this case, we have $d:=|z|$ is even and $p:=|s|$ is odd.
\blemma \la{loopcprxs}
If $X$ is a simply connected space such that $\H^\ast(X;\Q)\,\cong\,\Q[z,s]$, then there is an isomorphism of graded vector spaces
$$\rH^{S^1,(m)}_\ast(\LL X)\,\cong\, \bigoplus_{j=1}^{\infty}(\c \cdot \nu_j \oplus \c \cdot \eta_j)\,, $$
where the homological degrees of the basis elements are given by
$$\deg\,\nu_j\,=\,(p-1)m+dj-1\,,\qquad \deg\,\eta_j\,=\,(p-1)(m+1)+d(j-1)\ .$$
\elemma 
\bproof
Let $\mathcal{A}=\c[z,s]$ denote the (complexified) minimal Sullivan model of $X$ (see \cite[Prop. 5.1]{Me}).  Recall the computation of $\rHC^{(m)}_\ast(\mathcal A)$ in \eqref{cycminsullivan} and \eqref{se}. For $j\,\in\,\mathbb{N}$ and $\omega \,\in\,\rHC^{(m)}_\ast(\mathcal A)$, the coefficient of $z^{j-1}$ (resp., $z^{j-1}dz$) in $S(\omega)$ (resp., $E(\omega)$) determines a homogeneous linear functional $\eta_j$ (resp., $\nu_j$) on $\rHC^{(m)}_\ast(\mathcal A)$ of homological degree $p+m(p-1)+d(j-1)$ (resp., $m(p-1)+dj$). The desired lemma now follows from \cite[Thm. 4.2]{BRZ} and \cite[Prop. 7.8]{BFPRW}, which together imply 
$$\rH^{S^1,(m)}_\ast(\LL X)\,\cong\,\big(\rHC^{(m)}_\ast(\mathcal A)\big)^{\ast}[-1]\ .  $$
\eproof
Observe that if $X$ is as in Lemma \ref{loopcprxs} above, then the (complexified) minimal Sullivan model of $X$, which is given by $\mathcal{A}=\c[z,s]$ has a $\Z^2$-weight grading, with $z$ having weight $(1,0)$ and $s$ having weight $(0,1)$. The $\Z^2$-grading on $\mathcal{A}$ induces a $\Z^2$-grading on the Chevalley-Eilenberg cochain complex $\C^{-\ast}(\g(\bar{\mathcal{A}});k)$ (where the graded linear dual of a space of weight $(p,q)$ in homological degree $i$ has weight $(p,q)$ in homological degree $-i$) that is compatible with its homological grading, differential, as well as with the natural $\g$-action. By Theorem \ref{hrsullivan}, $\HR_\ast(X,G)$ acquires a $\Z^2$-grading compatible with the $G$-action. Let $P_{X,G}(q,t,z)$ denote the Euler-Poincar\'{e} series of the ($G$-invariant part of the) representation homology of $X$:
$$P_{X,G}(q,t,z):= \sum_{n=0}^{\infty}\ \sum_{(r,s)\,\in\,\mathbb{Z}^2} \dim_k[\,\HR_n(X,G)^G_{(r,s)}\,]\, q^r\, t^s\, z^n \ .
$$
Here $\HR_n(X,G)^G_{(r,s)}$ denotes the component of $\HR_n(X,G)^G$ with $\Z^2$-weight $(r,s)$. Note that the specialization of $P_{X,G}(q,t,z)$  at $z=-1$ gives the weighted Euler characteristic $\chi_{X,G}(q,t)$   of $\HR_\ast(X,G)^G$: $P_{X,G}(q,t,-1)=\chi_{X,G}(q,t)$. The following result is a consequence of Theorem \ref{drinhomcprxs} and Lemma \ref{loopcprxs}.
\bcor \la{rephomcprxs}
Let $X$ be a simply connected space such that $\H^\ast(X;\Q)\,\cong\,\Q[z,s]$, where $d=|z|$ is even and $p=|s|$ is odd. Then, there is an isomorphism of graded commutative algebras
$$\HR_\ast(X,G)^G\,\cong\,\Sym\, \big(\nu^{(i)}_j,\eta^{(i)}_j\,:\, i=1,2,\ldots,l\,,\,j \,\in\,\mathbb{N} \big)\,, $$
where the generators $\nu^{(i)}_j$ have homological degree $\deg\,\nu^{(i)}_j\,=\,(p-1)m_i+dj-1$ and the generators $\eta^{(i)}_j$ have homological degree $\deg\,\eta^{(i)}_j\,=\,(p-1)(m_i+1)+d(j-1)$. Further,
$$
P_{X,G}(q,t,z)\,=\, \prod_{i=1}^l\prod_{j=1}^{\infty} \, \frac{1\,+\,q^j\,t^{m_i}\,z^{\deg\,\nu^{(i)}_j} }{1\,-\,q^{j-1}\,t^{m_i+1}\,z^{\deg\, \eta^{(i)}_j}}\ .$$
\ecor
In particular, by letting $z=-1$ in the above formula, we obtain:
\begin{equation} \la{eulercharcprxs} \chi_{X,G}(q,t)\,=\, \prod_{i=1}^l\prod_{j=1}^{\infty} \frac{1\,-\,q^j\,t^{m_i}}{1\,-\,q^{j-1}\,t^{m_i+1}} \end{equation}
since $\deg\,\nu^{(i)}_j$ are always odd numbers and $\deg\,\eta^{(i)}_j$ are always even numbers. On the other hand, with Theorem \ref{hrsullivan}, the Euler characteristic $\chi_{X,G}(q,t)$ can be computed in a different way, from the chain complex $\C^{-\ast}(\g(\bar{\mathcal A});k)^G$, using standard Lie theoretic methods (see, e.g. \cite[Sec. 9.3]{BFPRW}, in paricular Corollary 9.8 therein): 
\begin{equation} \la{eulerchevalleycprxs} \chi_{X,G}(q,t)\,=\, \frac{1}{|W|} \prod_{j=1}^{\infty} \,\left(\frac{1\,-\,q^j}{1\,-\,q^{j-1}\,t}\right)^l\,\mathrm{CT}\left\{ \prod_{j=1}^{\infty} \,\prod_{\alpha\,\in\, R}\, \frac{1\,-
\,q^{j-1}\,e^{\alpha}}{1\,-\,q^{j-1}\,t\,e^{\alpha}} \right\} \ .\end{equation}
Here, $W$ is the Weyl group, $R$ the associated root system of $\g$, and $\mathrm{CT}\,:\,\Z[Q] \rar \Z$ is the classical constant term map defined on the group ring of the root lattice $Q=Q(R)$ of $R$. Comparing the right hand sides of \eqref{eulercharcprxs} and \eqref{eulerchevalleycprxs}, we obtain the celebrated Macdonald's $(q,t)$-constant term identity (see \cite{Macd}):
\begin{equation}\la{macdonaldqt} \frac{1}{|W|}\,\mathrm{CT}\left\{ \prod_{j=1}^{\infty} \,\prod_{\alpha\,\in\, R}\, \frac{1\,-
\,q^{j-1}\,e^{\alpha}}{1\,-\,q^{j-1}\,t\,e^{\alpha}} \right\}\,=\, \prod_{i=1}^l\prod_{j=1}^{\infty} \frac{(1\,-\,q^{j-1}\,t)\,(1\,-\,q^j\,t^{m_i})}{(1-q^j)\,(1\,-\,q^{j-1}\,t^{m_i+1})}\ . \end{equation}

We close this section by listing some spaces to which Corollary \ref{rephomcprxs} applies:
\begin{itemize}
\item  $K(\Z,d) \times \mathbb{S}^p$, where $d \geqslant 2$ is even and $p \geqslant 3$ is odd ($|z|\,=\,d\,,\,|s|\,=\,p$)
\item $\c\mathbb{P}^{\infty} \times \mathbb{S}^{2r+1}$ (rationally equivalent to $K(\Z,2) \times \mathbb{S}^{2r+1}$)
\item $\mathbb{H}\mathbb{P}^{\infty} \times \mathbb{S}^{4r+3}$ (rationally equivalent to $K(\Z,4) \times \mathbb{S}^{4r+3}$)
\item $\Omega(\mathbb{H}\mathbb{P}^r)$ (rationally equivalent to $K(\Z,4r+2) \times \mathbb{S}^{3}$)
\end{itemize} 
In particular, by Corollary \ref{rephomcprxs}, we have
\begin{equation*}
\la{cpsphere} 
\HR_\ast(\c\mathbb{P}^{\infty} \times \mathbb{S}^{2r+1},G)^G  
\cong  \Sym\, \big(\nu^{(i)}_j,\,\eta^{(i)}_j \big)\ ,
\quad \deg \nu^{(i)}_j = 2rm_i+2j-1 \ ,
\quad  \deg \eta^{(i)}_j = 2r(m_i+1)+2(j-1)\ ,
\end{equation*}
%
%
%
\begin{equation*} 
\la{hpsphere} \HR_\ast(\mathbb{H}\mathbb{P}^{\infty} \times \mathbb{S}^{4r+3},G)^G  \cong  \Sym\, \big(\nu^{(i)}_j,\,\eta^{(i)}_j\big)\ , \ \deg\,\nu^{(i)}_j\,=\,(4r+2)m_i+4j-1\ , \ \deg\,\eta^{(i)}_j\,=\,(4r+2)(m_i+1)+4(j-1)\ ,
\end{equation*}
%
%
\begin{equation*}
\la{omegacpr} \HR_\ast(\Omega(\mathbb{H}\mathbb{P}^r),G)^G  \cong  \Sym\, \big(\nu^{(i)}_j,\,\eta^{(i)}_j \big)\ , \quad \deg\,\nu^{(i)}_j\,=\,2m_i+(4r+2)j-1\ ,\quad \deg\,\eta^{(i)}_j\,=\,2(m_i+1)+(4r+2)(j-1),
\end{equation*}

\vspace{0.5ex}

\noindent
where $i\,\in\,\{1,2,.\ldots,l\}$ and $j\,\in\,\mathbb{N}$.
\subsubsection{Case III}
In this case, we have $d:=|z|$ is even.
\blemma \la{loopcpr}
If $X$ is a simply connected space such that $\H^\ast(X;\Q)\,\cong\,\Q[z]/(z^{r+1})$, then there is an isomorphism of graded vector spaces
$$ \rH^{S^1,(m)}_\ast(\LL X)\,\cong\,\bigoplus_{j=1}^r \c \cdot \xi_j\,,$$
where the homological degrees of the basis elements are given by
$$
\deg\,\xi_j\,=\,(d(r+1)-2)m+dj-1\ .
$$
\elemma
\bproof
Recall that the minimal Sullivan model of $X$ is given by $\mathcal{A}_r=\c[z,s]\,,\,\partial s=z^{r+1}$ (see \cite[Prop. 5.1]{Me}) where $z$ and $s$ have dimensions  $d$ and $d(r+1)-1$ respectively. As in \eqref{cycminsullivan}, for any $m \geq 0$,
$$\rHC_\ast^{(m)}(\mathcal{A}_r)\,\cong\, \H_\ast\left[ \Omega^m(\mathcal{A}_r)/d\Omega^{m-1}(\mathcal{A}_r\right]\,\cong\, \H_\ast[\c[z]s(ds)^m \oplus \c[z]dz\cdot s(ds)^{m-1},\partial]\,, $$
where the differential $\partial$ is induced by the differential on $\mathcal{A}_r$. A direct computation shows that 
$$\partial([z^ks(ds)^m])=-\big(k+(m+1)(r+1)\big)[z^{k+r}dz\cdot s(ds)^{m-1}]\,,\qquad \partial([z^kdz \cdot s(ds)^{m-1}])=0\ .  $$
Hence, 
$$ \rHC^{(m)}_\ast(\mathcal{A}_r)\,=\, \mathrm{Span}_\c\{ [z^kdz \cdot s(ds)^{m-1}]\,,\,0 \leqslant k < r \}\ .$$
Choose a basis $\{\xi_{j}\,,\, 1\leqslant i \leqslant l\,,\,1 \leqslant j \leqslant r\}$ of 
$\rHC^{(m)}(\mathcal{A}_r)^{\ast}[-1]$ dual to the basis 
$\{[z^{j-1}dz \cdot s(ds)^{m-1}]$ of $\rHC_\ast^{(m)}(\mathcal{A}_r)[1]$. Clearly, 
$|\xi_{j}|=(d(r+1)-2)m+dj-1$.  
Since 
$$\rH^{S^1,(m)}_\ast(\LL X)\,\cong\,\big(\rHC^{(m)}_\ast(\mathcal A)\big)^{\ast}[-1] $$
by \cite[Thm. 4.2]{BRZ} and \cite[Prop. 7.8]{BFPRW}, the desired lemma follows.
\eproof 
Observe that if $X$ is as in Lemma \ref{loopcpr} above, then (complexified) minimal Sullivan model of $X$, which is given by $\mathcal{A}=\c[z,s], \partial\,s=z^{r+1}$, has a $\Z$-weight grading with $z$ having weight $1$ and $s$ having weight $r+1$. The $\Z$-grading on $\mathcal{A}$ induces a $\Z$-grading on the Chevalley-Eilenberg cochain complex $\C^{-\ast}(\g(\bar{\mathcal{A}});k)$ (where the graded linear dual of a space of weight $p$ in homological degree $i$ has weight $p$ in homological degree $-i$) that is compatible with its homological grading, differential, as well as with the natural $\g$-action. By Theorem \ref{hrsullivan}, $\HR_\ast(X,G)$ acquires a $\Z$-grading compatible with the $G$-action. Let $P_{X,G}(q,z)$ denote the Euler-Poincar\'{e} series of the ($G$-invariant part of the) representation homology of $X$:
$$P_{X,G}(q,z):= \sum_{n=0}^{\infty}\ \sum_{ p \,\in\,\mathbb{Z}} \,\dim_k [\,\HR_n(X,G)^G_{p}\,]\, q^p \, z^n \ ,
$$
where $\HR_n(X,G)^G_{p}$ denotes the component of $\HR_n(X,G)^G$ with $\Z$-weight $p$. Note that $P_{X,G}(q,-1)$ is the weighted Euler characteristic $\chi_{X,G}(q)$   of $\HR_\ast(X,G)^G$. The following result is a consequence of Theorem \ref{drinhomcprxs} and Lemma \ref{loopcpr}.
\bcor \la{rephomcpr}
 Let $X$ be a simply connected space such that $\H^\ast(X;\Q)\,\cong\,\Q[z]/(z^{r+1})$, where $z$ is of (even) dimension $d$. Then there is an isomorphism of graded commutative algebras
$$
\HR_{\ast}(X, G)^{G}\,\cong\, 
\Sym\, (\xi^{(i)}_1,\,\xi^{(i)}_2,\,\ldots,\, \xi^{(i)}_{r}\,:\, i=1,2,\ldots,l)\ ,$$
where the generators $\xi^{(i)}_{j}$ have homological degree $\deg\,\xi^{(i)}_j\,=\, \,(d(r+1)-2)m_i+dj-1\,$. Further,
$$ P_{X,G}(q,z)\,=\, \prod_{i=1}^l \prod_{j=1}^r (1\,+\,q^{j+m_i(r+1)}\,z^{\deg\,\xi^{(i)}_j})\ .$$
\ecor
In particular, specializing the above formula at $z=-1$, we obtain:
\begin{equation} \la{eulercharcpr} \chi_{X,G}(q)\,=\,\prod_{i=1}^l \prod_{j=1}^r (1\,-\,q^{j+m_i(r+1)})\ . \end{equation}
As in Section \ref{s4.4.2}, by Theorem \ref{hrsullivan}, $\chi_{X,G}(q,t)$ can also be computed as the weighted Euler characteristic of $\C^{-\ast}(\g(\bar{\mathcal A});k)^G$, using standard Lie theoretic methods. The calculations similar to those in \cite[Sec. 9.3.1]{BFPRW} give:
\begin{equation} \la{eulerchevalleycpr} \chi_{X,G}(q)= \frac{1}{|W|} \prod_{j=1}^r(1\,-\,q^j)^l\,\mathrm{CT}\left\{ \prod_{j=0}^r\,\prod_{\alpha\,\in\,R} (1\,-\,q^j\,e^{\alpha})\right\}\ .\end{equation}
Equating the right hand sides of equations \eqref{eulercharcpr} and \eqref{eulerchevalleycpr} (and dividing both expressions by $\prod_{j=1}^r (1-q^j)^l$), we obtain the Macdonald's $q$-constant term identity (see \cite{Macd})
\begin{equation} \la{macdonaldq} \mathrm{CT}\left\{ \prod_{j=0}^r\,\prod_{\alpha\,\in\,R} (1\,-\,q^j\,e^{\alpha})\right\}\,=\, \prod_{i=1}^l \prod_{j=1}^r \frac{1\,-\,q^{j+m_i(r+1)}}{1\,-\,q^j}\ .\end{equation}
The classical spaces satisfying the conditions of Corollary \ref{rephomcpr} are:

\begin{itemize}
\item the even-dimensional spheres $\bS^{2n}$  ($r=1,\,d=2n$)\,,
 \item  the complex projective spaces $\c\mathbb{P}^r$ ($r \ge 1,\,d=2 $)\,,
 \item the quaternionic projective spaces $\mathbb{H}\mathbb{P}^r$ ($r \ge 1,\,d=4$)\,.
\item  the Cayley projective plane $\mathbb{O}\mathbb{P}^2$ ($r=2,\,d=8$)\,.
\end{itemize}

For these spaces, Corollary \ref{rephomcpr} gives
\begin{eqnarray*}
\HR_{\ast}(\bS^{d},G)^G  &\cong & \Sym\, [\,\xi^{(i)}\, :\, 1 \leq i \leq l\,]\ ,\qquad\qquad\qquad\quad\,\, \deg\,\xi^{(i)}\,=\,(d-1) (2m_i+1)\ , \\
\HR_{\ast}(\c\mathbb{P}^{r},G)^G &\cong & \Sym\, [\,\xi^{(i)}_j\, : \, 1 \leq i \leq l\ ,\ 1 \leq j \leq r\,]\ ,\qquad \deg\,\xi^{(i)}_j\,=\,2rm_i+2j-1\ , \\
\HR_{\ast}(\mathbb{H}\mathbb{P}^{r},G)^G & \cong & \Sym\, [\,\xi^{(i)}_j\, : \, 1 \leq i \leq l\ ,\ 1 \leq j \leq r\,]\ ,\qquad \deg\,\xi^{(i)}_j\,=\,(4r+2)m_i+4j-1 \ , \\
\HR_{\ast}(\mathbb{O}\mathbb{P}^{2},G)^G &\cong & \Sym\, [\,\xi^{(i)}_1,\,\xi^{(i)}_2\, :\, 1 \leq i \leq l\,]\,,\qquad\qquad\quad\,\, \deg\,\xi^{(i)}_j\,=\,22m_i+8j-1\ .
\end{eqnarray*}
%
%

\vspace*{2ex}

We close this section with one curious consequence of Corollary \ref{rephomcpr}: it shows how
knowing the exact structure of the $G$-invariant part of representation homology allows one   
(sometimes) to get information about the full representation homology.
\blemma \la{vanrephomcpr}
Let $X$ be a simply connected space such that $\,\H^{\ast}(X;\Q)\,\cong\,\Q[z]/(z^{r+1})\,$, where $|z| = d \ge 2 $ is even. Put 
$$ 
N := \frac{1}{2}\, r\, (d(r+1) - 2)\, \dim G\ .
$$
Then, $\,\HR_n(X,G) = 0\,$ for all $\,n > N\,$. Moreover, $\,\HR_N(X,G)\,\cong\, \c \,$.
\elemma
\bproof
By a classical theorem of Kostant \cite{Ko}, we have $\,\sum_{i=1}^l (2 m_i + 1) = \dim\, G\,$
for any complex reductive group $G$. This implies that
$$
\sum_{i=1}^l\,\sum_{j=1}^r \,\deg\,\xi^{(i)}_j\,=\,\frac{1}{2}\, r\, (d(r+1) - 2)\, \dim G \, =: N\ .
$$
where $\, \deg\,\xi^{(i)}_j\,$ are the degrees of the free generators of $ \HR_{\ast}(X,G)^G $ given in Corollary \ref{rephomcpr}. By Corollary \ref{rephomcpr}, we then conclude that 
$\, \HR_N(X,G)^G = 0\,$ for all $ n > N $, while $\,\dim_\c\,\HR_N(X,G)^G = 1 $.

On the other hand, the (complexified) minimal Sullivan model $\mathcal{A}_X$ of $X$ is formal (indeed, the map of DG algebras $(\c[z,s], \partial s = z^{r+1}) \rar \c[z]/(z^{r+1})$ given on generators by $z \mapsto z$, $s \mapsto 0$ is obviously  a quasi-isomorphism). Hence, by part $(a)$ of Theorem \ref{hrsullivan}, the $\,\HR_\ast(X,G)$ is isomorphic the homology of the Chevalley-Eilenberg complex $\C^{-\ast}(\g(\rH^{\ast}(X;\c));\,\c)$. By definition, this last complex is a graded exterior algebra on $\, r \cdot \dim \g\, $ generators of homological degree $\,dj-1\,$, where $\,j =1,2,\ldots, r\,$. Therefore, its homology $\,\HR_n(X,G)\,$ is {\it a fortiori}\, concentrated in homological degrees $\, n \leq N' \,$, where
$$
N' := \sum_{j=1}^r (dj-1)\,\dim \g
$$
Moreover, $\, \dim_{\c}\,\HR_{N'}(X,G) \leq 1 \,$. A trivial calculation shows that
$\, N' = N \,$. Since $ \HR_{\ast}(X,G)^G \subseteq \HR_*(X,G) $, this numerical 
coincidence implies the result of the lemma.
\eproof

\appendix
\section{Monoidal Dold-Kan correspondence} 
\la{doldkan}

The Dold-Kan correspondence is a classical result that establishes an equivalence between
the category $\Chp(\mathscr{A})$ of non-negatively graded chain complexes in an abelian category $\mathscr{A}$
and the category $s\mathscr{A}$ of simplicial objects in $\mathscr{A}$.
In this Appendix, we will describe a monoidal enrichment of this correspondence relating the category of (non-negatively graded) DG $\cP$-algebras to the category of simplicial $\cP$-algebras
for an arbitrary $k$-linear operad $\cP$. For simplicity, we will fix a commutative ring $k$ with unit, and consider only the abelian category
$\mathscr{A} = \Mod_k$.

\subsection{The Dold-Kan correspondence} 
\la{sdk}
To any simplicial $k$-module $ X \in \sMod_k $ we can associate
the chain complex
$$
N(X) = [\,\ldots\to N(X)_{n} \xrightarrow{\partial} N(X)_{n-1}\to \ldots \,]
$$
with $ N(X)_n := \bigcap_{i=1}^n \Ker(d_i:X_n \rar X_{n-1}) $ for $ n \ge 0 $ and the differential $ \partial $ given by  $ d_0$.
The assignment $ X \mapsto  N(X) $ defines a functor $ N: \sMod_k \to \Chp(k) $
from the category of simplicial $k$-modules to the category  of connective chain complexes of $k$-modules.
The functor $N$ is called the {\it normalization functor}. A classical theorem due to Dold and Kan (see \cite[Theorem 8.4.1]{W}) asserts that $ N $ is an equivalence of categories.

For any simplicial $k$-module $ X \in \sMod_k$,  the homology groups of the chain complex $ N(X) $ are naturally isomorphic to the homotopy groups $\pi_{\ast}(|X|)$ of the geometric realization of $X$
(see \cite[Theorem 22.1]{M}).
This justifies the notation $\pi_\ast(X) := {\rm H}_\ast[N(X)]$, which we used throughout the paper.

The inverse $N^{-1} : \Chp(k) \rar \sMod_k$
of the normalization functor is defined as follows.
For any chain complex $V \in \Chp(k)$, the degree $n$ part of the
simplicial $k$-module $N^{-1}(V)$ is given by
\begin{equation}  \label{N_inv}
N^{-1}(V)_n \, = \, \Moplus_{r \geq 0} \, \Moplus_{\sigma : [n] \twoheadrightarrow [r]} \, V_r
\end{equation}

We think of $N^{-1}(V)$ as adjoining to $V$ the degeneracies of all elements in $V$.
We write an element $x \in V_r$ in the summand corresponding to $\sigma$
as $\sigma^*(x)$.
When $\sigma=\id$, we simply write this as $x$, or $\eta(x)$ if we want to emphasize that
we consider $x$ to be an element in $N^{-1}(V)$ rather than $V$.
As suggested by the notation, this determines the degeneracy maps in $N^{-1}(V)$:
 namely, $s_j(\sigma^{\ast}(x)) := (\sigma \circ \sigma^j)^{\ast}(x)$ (recall that $s_j:= [\sigma^j]^{\ast}$).
The face maps in $N^{-1}(V)$ are determined by the requirement that $d_i(\eta(x))=0$ for all $i>0$,
and the canonical map
\begin{equation*}
\eta \, : \, V \rar N[ N^{-1}(V) ], \qquad \qquad \quad x \, \mapsto \, \eta(x) = x
\end{equation*}
commutes with differentials,
{\it i.e.} $d_0(\eta(x)) = \eta(d(x))$.
Since all elements of $N^{-1}(V)$ other than $\eta(x)$ are sums of degenerations of $\eta(x)$,
specifying the face maps on these elements determines all the face maps in $N^{-1}(V)$.
This defines a simplicial $k$-module $N^{-1}(V)$
and hence the functor $N^{-1} : \Chp(k) \rar \sMod_k$ (see \cite{GJ} for more details).
It is easy to check that this functor is indeed the inverse of the normalization functor $N$.

There is an alternative way to define the normalization functor.
For each simplicial $k$-module $X \in \sMod_k$, we can take the chain complex $\overline{N}(X)$ defined by

\begin{equation}   
\label{norm_form2}
\overline{N}(X)_n \, := \,  \frac{ X_n }{ \sum_{j=0}^{n-1} s_j(X_{n-1}) }
\qquad \qquad d = \sum_{i=0}^n (-1)^i d_i \, : \,\overline{N}(X)_n \rar \overline{N}(X)_{n-1}
\end{equation}
Then one can show (see \cite{GJ}) that the canonical map $N(X) \rar \overline{N}(X)$ of chain complexes
given by the composition $N(X)_n \hookrightarrow X_n \twoheadrightarrow \overline{N}(X)_n$
is an isomorphism.

Notice that the inverse \eqref{N_inv} of the normalization functor has an important feature:
the collection of $k$-modules $N^{-1}(V)_n$, as well as the degeneracy maps between them,
depends only on the graded module $V$ and not on its differential.
In other words, \eqref{N_inv} defines a functor
$N^{-1} \, : \, \grMod_k \rar \Mod_k^{ \Delta^{\op}_{\surj} }$
from the category $\grMod_k$ of graded $k$-modules to the category $\Mod_k^{ \Delta^{\op}_{\surj} }$
of  $\Delta^{\op}_{\surj} $-systems of $k$-modules.
Similarly, \eqref{norm_form2} gives a functor
$\overline{N} \, : \, \Mod_k^{ \Delta^{\op}_{\surj} } \rar \grMod_k$.
This will play a role in our construction of the monoidal Dold-Kan correspondence in the next section.

\subsection{Monoidal Dold-Kan correspondence}

It is a classical fact 
that the Dold-Kan normalization functor $N : \sMod_k \rar \Chp(k)$
can be endowed with a symmetric lax monoidal structure.
To describe it, we first introduce some notations.
Given two simplicial modules $X,Y \in \sMod_k$ over a commutative ring $k$,
we denote by $X \sten Y \in \sMod_k$ the result of applying the tensor product levelwise,
{\it i.e.,} $(X \sten Y)_n := X_n \otimes_k Y_n$.
Then, there is a \emph{quasi-isomorphism} of chain complexes
\[ \sh:\, N(X) \otimes N(Y) \, \xrightarrow{\sim} \, N( X \sten Y)  \]
called the \emph{Eilenberg-Zilber shuffle map}, which is natural (in $X$ and $Y$), symmetric, associative and unital in the obvious sense (see, e.g., \cite{M, SS} for details).

This shuffle map allows one to transfer algebraic structures from a simplicial module $A$
to its normalization $N(A)$.
For instance, if $A$ is a simplicial associative algebra, then $N(A)$ is a DG algebra;
if $A$ is a simplicial commutative algebra, then $N(A)$ is a commutative DG algebra, {\it etc}.
In general, for any $k$-linear operad $\cP$,
one can consider the category $\sAlg(\cP)$ of simplicial $\cP$-algebras
as well as the category $\dgAlg(\cP)$ of non-negatively graded DG $\cP$-algebras.
If $A \in \sAlg(\cP)$ is a simplicial $\cP$-algebra,
then each $n$-ary operation $\mu \in \cP(n)$ gives a map
\[ \alpha_A(\mu) \, : \, A \, \sten \, \stackrel{n}{\ldots} \, \sten \, A \rar A \]
One can then use the Eilenberg-Zilber shuffle maps to construct the maps
\[ \alpha_{N(A)}(\mu) \, : \, N(A) \, \otimes \, \stackrel{n}{\ldots} \, \otimes \, N(A)
\xra{\sh} \, N( \, A \, \sten \, \stackrel{n}{\ldots} \, \sten \, A \, ) \xra{ N(\alpha_A(\mu)) } N(A) \]
which form the structure maps for a DG $\cP$-algebra on $N(A)$. This defines a functor
\begin{equation}
\la{monoidalN} N : \sAlg(\cP) \rar \dgAlg(\cP) \ .
\end{equation}
In the special case when $\cP$ is the Lie operad, this last functor has already appeared in \cite{Q}.
Quillen showed that it has a left adjoint in that case. His proof generalizes directly to an arbitrary operad.
\bprop
The functor \eqref{monoidalN} has a left adjoint $N^* : \dgAlg(\cP) \rar \sAlg(\cP)$.
\eprop

\bproof
As in \cite{Q}, for any $A \in \dgAlg(\cP)$, we define $N^*(A)$ as the following (degreewise) coequalizer
of simplicial $\cP$-algebras
\[ N^*(A) \, = \, \coeq \, \left[
\vcenter{ \xymatrix{
T_{\cP}(N^{-1}(T_{\cP}(A)))   \ar@<0.5ex>[r]^-{\alpha_\ast} \ar@<-0.5ex>[r]_-{\sh_\ast}
& T_{\cP}(N^{-1}(A))
}}
 \right]  \]
where $\alpha_\ast$ and $\sh_\ast$ are induced by $N^{-1}(\alpha) \, : \, N^{-1}(T_{\cP}(A)) \rar N^{-1}(A)$
and the Eilenberg-Zilber maps $\sh \, : \, N^{-1}(T_{\cP}(A)) \rar T_{\cP}(N^{-1}(A))$ respectively.
\eproof

Notice that the proof shows that the underlying $\Delta^{\op}_{\surj}$-system of $\cP$-algebra
of $N^*(A)$ depends only on the graded algebra structure of $A$
(see the discussion at the end of the previous subsection).
This observation will allow us to describe the simplicial $\cP$-algebra $N^*(A)$
in the case when $A$ is semi-free.


We first consider the commutative diagrams of functors
\[
\vcenter{
\xymatrix{
\sAlg(\cP)  \ar[r]^-{{\rm forget}}  \ar[d]_{N}  \ar@{}[rd]|{(1)}  &  \sMod_k \ar[d]^{N}  \\
\dgAlg(\cP)  \ar[r]^-{{\rm forget}}           &  \Chp(k)
} }
\qquad \qquad
\vcenter{
\xymatrix{
\sAlg(\cP)  \ar@{}[rd]|{(2)}    &  \sMod_k  \ar[l]_-{T_{\cP}}  \\
\dgAlg(\cP)  \ar[u]^{N^*}          &  \Chp(k) \ar[u]_{N^{-1}}  \ar[l]_-{T_{\cP}}
} }
\]
where we denote by $T_{\cP}$ the free algebra functors in both simplicial and DG contexts.
The square $(1)$ obviously commutes up to isomorphism of functors.
The square $(2)$ is obtained by replacing every functor on the square $(1)$ by its left adjoint.
Therefore, it also commutes up to isomorphism of functors. The commutativity (up to isomorphism) of the square $(2)$ can be written as
\begin{equation}   \label{Nstar_semifree_formula}
N^*(T_{\cP}(V))  \cong T_{\cP} (N^{-1}(V))\ .
\end{equation}
In other words, $N^*$ of a free DG $\cP$-algebra is free. The same is true for semi-free algebras.
Recall that a DG $\cP$-algebra is said to be \emph{semi-free} if its underlying graded algebra is free
over a degreewise free graded $k$-module $V$.
Similarly, a simplicial $\cP$-algebra $A$ is said to be \emph{semi-free}\footnote{By standard definition (cf. \cite{GJ}), a simplicial $\cP$-algebra is called semi-free
if there is a collection of subsets $B_n \subset A_n$, called a \emph{basis},
that is closed under degeneracies and that $A_n = T_{\cP}(B_n)$ for each $n$.
It is clear that our definition implies this.
To see the converse, notice that any basis element
that is not the degeneracy of any other basis element is in fact non-degenerate in the underlying
simplicial set of $A$.
Let $V$ be the graded $k$-module with a basis given by these non-degenerate basis elements.
Then an application of \cite[Lemma I.2.11]{GM} shows that $A = T_{\cP}(N^{-1}(V))$
as a $\Delta^{\op}_{\surj}$-system of $\cP$-algebras.} if its underlying
$\Delta^{\op}_{\surj}$-system of $\cP$-algebras is of the form $A = T_{\cP}(N^{-1}(V))$
for a degreewise free graded $k$-module $V$.

\vspace{1ex}

The above discussion leads to the following
\bprop  \label{Nstar_semifree_prop}
The functor $N^* : \dgAlg(\cP) \rar \sAlg(\cP)$ sends semi-free DG $\cP$-algebras
to semi-free simplicial $\cP$-algebras.
%
\eprop

\bproof
We have seen that $N^*$ sends free algebras to free algebras.
Since the underlying $\Delta^{\op}_{\surj}$-system of $N^*(A)$ depends only on the graded algebra structure of $A$,
the result follows.
\eproof

%
%

Next, we consider the adjunction map $A \rar N(N^*(A))$ in the case when $A = T_{\cP}(V)$ is semi-free
over a graded complex $V$.
To describe this map, we first give a different interpretation of the Eilenberg-Zilber shuffle map.
Namely, we view it a collection of maps that connect two symmetric monoidal structures on the category $\Chp(k)$
of chain complexes on $k$.
We will use the ``quotient" form \eqref{norm_form2} of the normalization functor.
Thus, we consider the equivalence of categories $\overline{N} : \sMod_k \rar \Chp(k)$.
One can use this equivalence to transport the symmetric monoidal structure
$\sten$ on $\sMod_k$ to a symmetric monoidal structure $\uten$ on $\Chp(k)$.
Precisely, we define $V \uten W := \overline{N}(N^{-1}(V) \sten N^{-1}(W))$ for $V,W \in \Chp(k)$.
Then the Eilenberg-Zilber shuffle maps can be written as
\begin{equation}  \label{sh_map2}
\sh : V \otimes W \rar V \uten W, \qquad \qquad x\otimes y \mapsto x \times y := \sh(x,y)
\end{equation}

Now, suppose that a DG $\cP$-algebra $A \in \dgAlg(\cP)$ is semi-free over a graded $k$-module $V$, {\it i.e.,}
\[ A \, = \, T_{\cP}(V) \, := \, \Moplus_{n\geq 0} \cP(n) \otimes_{S_n} V^{\otimes n}  \]
then by \eqref{Nstar_semifree_formula} as well as the discussion that follows,
the DG $\cP$-algebra $\overline{N}(N^*(A)) \in \dgAlg(\cP)$ has a similar description
\[ \overline{N}(N^*(A)) \, = \, \Moplus_{n\geq 0} \cP(n) \uten_{S_n} V^{\uten n}  \]

Moreover, the adjunction map $A \rar \overline{N}(N^*(A))$ is given by
\begin{equation}  \label{adj_map}
\Moplus_{n\geq 0} \cP(n) \otimes_{S_n} V^{\otimes n}
\rar \Moplus_{n\geq 0} \cP(n) \uten_{S_n} V^{\uten n} ,
\qquad
\left( \mu \, , \, x_1 \otimes \ldots \otimes x_n  \right)
\, \mapsto \,
\left( \mu \, , \, x_1 \times \ldots \times x_n  \right)
\end{equation}

This description of the adjunction map will be useful in the next subsection
when we compare the model structures on simplicial $\cP$-algebras and DG $\cP$-algebras.

\subsection{Quillen equivalence}

By \cite[Section II.4, Theorem 4]{Q1}, there is a model structure on the category $\sAlg(\cP)$
of simplicial $\cP$-algebras, where a map $f: A \rar B$ is a weak equivalence (resp., fibration)
if and only if the map of the underlying simplicial sets is a weak equivalence (resp., fibration). Moreover, it is shown in \cite{Hin97} that if $k$ is a field of characteristic $0$, then
the category $\dgAlg(\cP)$ of DG $\cP$-algebras also has a model structure, where a map
$f: A \rar B$ is a weak equivalence (resp., fibration)
if and only if the map of the underlying (connective) chain complexes is a weak equivalence (resp., fibration).

From now on, we assume that $k$ is a field of characteristic $0$,
and the categories $\sAlg(\cP)$ and $\dgAlg(\cP)$ are equipped with the model structures described above.
Then, the normalization functor $N : \sAlg(\cP) \rar \dgAlg(\cP)$ preserves fibrations and
weak equivalences, and therefore the associated adjunction
\begin{equation}   
\label{NNstar_adj}
\xymatrix{ N^* \, : \, \dgAlg(\cP) \ar@<0.5ex>[r]     & \sAlg(\cP) \, : \, N  \ar@<0.5ex>[l] }
\end{equation}
is a Quillen pair.
In fact, we have the following theorem, which is the main result of this Appendix.

\bthm  \label{NNstar_equiv}
The Quillen pair \eqref{NNstar_adj} is a Quillen equivalence.
\ethm

\bproof
It suffices to show that, for any semi-free DG $\cP$-algebra $A = T_{\cP}(V)$, the
unit of the adjunction \eqref{NNstar_adj} is a weak equivalence.
Composing this adjunction map with the isomorphism $N(N^*(A)) \, \rar \, \overline{N}(N^*(A))$,
we can consider the map
$ \label{unitadj} A \rar \overline{N}(N^*(A))$, 
which depends only on the underlying graded $\cP$-algebra structure of $A$,
and is described explicitly by \eqref{adj_map}.

If $A$ is free, i.e. when the differential on $A = T_{\cP}(V)$ is induced by
the differential on a chain complex $V$, then for each $n \geq 0$, the map
\[ \sh : \cP(n) \otimes V^{\otimes n} \ra \cP(n) \uten V^{\uten n} \]
is a quasi-isomorphism as it is induced by the Eilenberg-Zilber shuffle map
(which is always a quasi-isomorphism).
Since $k$ is a field of characteristic $0$, the same is true when we pass to $S_n$-coinvariants
\[ \sh : \cP(n) \otimes_{S_n} V^{\otimes n} \ra \cP(n) \uten_{S_n} V^{\uten n} \]
This shows that the map \eqref{adj_map} is a quasi-isomorphism in the case when $A$ is free.

In the general case, when $A = T_{\cP}(V)$ is semi-free over a graded $k$-module $V$,
choose a homogeneous basis of $V$, and assign a weight grading $\wt(x) \in \mathbb{N}$
for each such basis element $x$.
This induces a grading on $A$, where an element
$(\mu \, , \, x_1 \otimes \ldots \otimes x_n ) \in \cP(n) \otimes_{S_n} V^{\otimes n}$
has weight grading $\wt(x_1) + \cdots + \wt(x_n)$.

The underlying $\Delta^{\op}_{\surj}$-system of $N^*(A)$ is given by $N^*(A) = T_{\cP}(N^{-1}(V))$.
Therefore, its elements in degree $m$ are sums of elements of the form
$(\mu \, , \, \sigma_1^*(x_1) \otimes \ldots \otimes \sigma_n^*(x_n) ) $
where $\phi_i \, : \, [m] \twoheadrightarrow [r_i]$ and $x_i \in V_{r_i}$.
Assign the weight grading $\wt(x_1) + \cdots + \wt(x_n)$ to this element.
Then it is clear that all the degeneracy maps preserve this weight grading.
This induces a weight grading in the normalization $\overline{N}(N^*(A))$.
Moreover, the map \eqref{adj_map} preserves this grading.
We write this grading as
\begin{equation}
 A = \Moplus_{n \geq 0} A^{(n)} \, , \qquad \qquad N^*(A) = \Moplus_{n \geq 0} N^*(A)^{(n)} \, , 
 \qquad \qquad  \overline{N}(N^*(A)) = \Moplus_{n \geq 0} \overline{N}(N^*(A))^{(n)}
\end{equation}

In general, the differentials on both sides of \eqref{adj_map} do not preserve the grading.
However, if we let $F_n(A) = \Moplus_{i \leq n} A^{(i)}$ be the filtration on $A$ induced by the weight grading,
then we can always choose the weight grading on a homogeneous basis of $V$
inductively so that $d(F_n(A)) \subset F_{n-1}(A)$.
Moreover, if we let $G_n = G_n( \overline{N}(N^*(A)) ) = \Moplus_{i \leq n} \overline{N}(N^*(A))^{(i)}$ be the filtration on  $\overline{N}(N^*(A))$
induced by the weight grading on $\overline{N}(N^*(A))$,
then we claim that $d( G_n ) \subset G_n$ for all $n$.

Indeed, consider the filtration  $\tilde{G}_n = \Moplus_{i \leq n} N^*(A)^{(i)}$ on $N^*(A)$ induced by the weight grading.
As we have seen, each graded piece $N^*(A)^{(i)}$ is a $\Delta^{\op}_{\surj}$-system of $k$-modules such that $\overline{N}(N^*(A)^{(i)}) = \overline{N}(N^*(A)^{(i)})$. 
Therefore, to show that $d(G_n) \subset G_n$, 
it suffices to show that $d_i(\tilde{G}_n) \subset \tilde{G}_n$ for all face maps $d_i$.
We will in fact show a more refined statement. 
To express this statement, we recall that the face maps of $N^*(A)$ are determined by the fact that the adjunction map 
$\eta : A \ra N(N^*(A))$ commutes with the differential. 
Indeed, for each homogeneous basis element $x \in V$, considered as an element in $T_{\cP}(V) = A$, 
the requirements $d_0(\eta(x)) = d(\eta(x)) = \eta(d(x))$ and $d_i(\eta(x)) = 0$
specify the values of face maps on the non-degenerate generators $\eta(x)$ of $N^*(A) = T_{\cP}(N^{-1}(V))$.
This in turn specifies the face maps on every other elements by simplicial identites.
Thus, one can write the face maps as $d_i = d_i[d_A]$ to show its dependence on the differential $d_A$ on $A$.
In Lemma~\ref{lemmaA1} below, we will show that, for any differential $d = d_A$ on $A$ such that $d(F_n) \subset F_{n-1}$, the face maps $d_i$
when restricted to homogeneous elements $z \in N^*(A)^{(n)}$, can be decomposed as
$d_i[d_A](z) = d_i'(z) + d_i''[d_A](z)$, where $d_i' : N^*(A)^{(n)} \ra N^*(A)^{(n)}$ and $d_i''[d_A] : N^*(A)^{(n)} \ra \tilde{G}_{n-1}$.
Moreover, $d_i'$ does not depend on the differential $d_A$, and $d_i''[d_A] = 0$ if $d_A=0$.
In particular, we have $d_i' = d_i[0]$.


Assuming this lemma, then we have $d_i(\tilde{G}_n) \subset \tilde{G}_n$, and hence $d(G_n) \subset G_n$.
Therefore, both the domain and target of the map of chain complexes \eqref{adj_map}
admit filtrations by subcomplexes, such that \eqref{adj_map} preserves these subcomplexes.
Since these filtrations are induced by gradings, the graded $k$-modules associated to these filtrations
can be canonically identified with the original graded $k$-modules, i.e. we have
\begin{equation}   \label{assoc_gr}
 {\rm gr}_F(A) \, \cong \, \Moplus_{n \geq 0} A^{(n)} \, = \, A
\ , \qquad 
{\rm gr}_G(\overline{N}(N^*(A))) \, \cong \, \Moplus_{n \geq 0} \overline{N}(N^*(A))^{(n)} \, = \, \overline{N}(N^*(A))
\end{equation}
as graded $k$-modules.
While passing to the associated graded modules does not change the underlying graded $k$-modules, it changes the differentials by 
discarding the part that strictly decrease the grading.
Since we have chosen the differential $d$ on $A$ such that $d(F_n) \subset F_{n-1}$, the induced differential on ${\rm gr}_F(A)$ is zero.
In other words, \eqref{assoc_gr} actually identifies $ {\rm gr}_F(A)$ with the free DG $\cP$ algebra $A' = T_{\cP}(V)$ with trivial differential.
On the other hand, by the above discussion, Lemma~\ref{lemmaA1} gives a description of the differential on the associated 
graded $k$-module ${\rm gr}_G(\overline{N}(N^*(A)))$.
Namely, by discarding the part of the differential on $\overline{N}(N^*(A))$ that strictly decreases the grading, 
one retains precisely the differential in $\overline{N}(N^*(A'))$ where $A'$ is again the free DG $\cP$ algebra $A' = T_{\cP}(V)$ with trivial differential.
In other words, \eqref{assoc_gr} actually identifies ${\rm gr}_G(\overline{N}(N^*(A)))$ with $\overline{N}(N^*(A'))$.
 


%

Therefore, the induced map ${\rm gr}_F(A) \ra {\rm gr}_G(\overline{N}(N^{-1}(A)))$
on the associated graded chain complexes coincides with the adjunction map
$ A' \ra \overline{N} (N^*(A')) $
for the free algebra $A' = T_{\cP}(V)$ with zero differential.
This map is a quasi-isomorphism by our previous argument in the free case.
Since the filtrations $F_\bullet$ and $G_\bullet$ are bounded below and exhaustive, the map
\eqref{adj_map} induces an isomorphism on homology by the Eilenberg-Moore comparison theorem \cite[Theorem~5.5.11]{W}.
\eproof

\blemma
\la{lemmaA1}
For any differential $d = d_A$ on $A$ such that $d(F_n) \subset F_{n-1}$, 
let $d_i = d_i[d_A]$ be the $i$-th face maps on $N^*(A)$.
Then its restriction $d_i|_{N^*(A)^{(n)}}$  to each homogeneous component $N^*(A)^{(n)}$ can be decomposed as
$d_i[d_A](z) = d_i'(z) + d_i''[d_A](z)$, where $d_i' : N^*(A)^{(n)} \ra N^*(A)^{(n)}$ and $d_i''[d_A] : N^*(A)^{(n)} \ra \tilde{G}_{n-1}$.
Moreover, $d_i'$ does not depend on the differential $d_A$, and $d_i''[d_A] = 0$ if $d_A=0$.
\elemma

\bproof
In simplicial degree $m$, the $k$-module $N^*(A)^{(n)}_m$ consists of sums of elements of 
the form 
\[ z = (\mu \, , \, \sigma_1^*(x_1) \otimes \ldots \otimes \sigma_k^*(x_k) ) \]
with  $\wt(x_1) + \cdots + \wt(x_k) = n$,
where $x_j \in V_{r_j}$ and $\sigma_j$ are surjective maps $[m] \twoheadrightarrow [r_j]$.
The image under the face map $d_i$ of this element is given by 
\begin{equation} \label{face_of_z}
 d_i(z) \, = \, (\mu \, , \, d_i(\sigma_1^*(x_1)) \otimes \ldots \otimes d_i(\sigma_k^*(x_k) ))
\end{equation}
Now, for each $j=1,\ldots, k$, the element $d_i(\sigma_j^*(x_j))$ reduces by simplicial identities to either of the two cases:
\begin{enumerate}
 \item[(I)]  $d_i(\sigma_j^*(x_j)) = \sigma_j'^*(x_j)$ for some surjective map $\sigma_j' : [m-1] \twoheadrightarrow [r_j]$ in $\Delta$.
 \item[(II)] $d_i(\sigma_j^*(x_j)) = \sigma_j'^*(d_{i'}(x_j)))$ for some surjective map $\sigma_j' : [m-1] \twoheadrightarrow [r_j-1]$ in $\Delta$, and some $0 \leq i' \leq r_j$.
\end{enumerate}
In case (I), $d_i(\sigma_j^*(x_j))$ has the same weight grading as the term $\sigma_j^*(x_j)$.
We split the case (II) in two subcases:
\begin{enumerate}
 \item[(IIa)] If $i' > 0$, then we have $d_{i'}(x_j)=0$ because by definition $x_j = \eta(x_j)$ is in $N(N^*(A))$.
 \item[(IIb)] If $i' = 0$, then we claim that $d_0(x_j) \in \tilde{G}_{r_j -1}$, and $d_0(x_j) = 0$ if $d_A = 0$.
\end{enumerate}
Indeed, the $0$-th face $d_0(x_j)$ of $x_j = \eta(x_j)$ is uniquely determined by the corresponding differential $d(x_j)$ in the DG $\cP$-algebra $A$.
Namely, since the adjunction map $\eta : A \ra N(N^*(A))$ commutes with differentials, we have $d_0(\eta(x_j)) = \eta(d(x_j))$.
Since we have chosen the weight grading on the generators $x_j$ in such a way that $d(x_j)$ is sum of terms of weight grading strictly less than $x_j$,
we see that $d_0(x_j) \in \tilde{G}_{r_j -1}$ in this case.
The equation $d_0(\eta(x_j)) = \eta(d(x_j))$ also shows that $d_0(x_j) = 0$ if $d_A = 0$.

Thus, to compute $d_i(z)$, one combines the equation \eqref{face_of_z} with the above consideration.
If we are in case (I) or (IIa) for each $1 \leq j \leq k$, then $d_i(z)$ is still in $N^*(A)^{(n)}$.
Thus, we have $d_i(z) = d_i'(z)$ in this case. Moreover, our explicit description shows that $d_i'$ does not depend on $d_A$.
If we are in case (IIb) for some $1 \leq j \leq k$, then we have $d_i(z) \in \tilde{G}_{n-1}$.
Thus, we have $d_i(z) = d_i''(z)$ in this case. Moreover, our description shows that $d_i(z)=0$ in this case if $d_A = 0$.
This completes the proof of the lemma.
\eproof 


\end{document}